\documentclass[a4paper,12pt, twoside]{book}
\usepackage{hyperref}

\usepackage{amssymb,amsmath,graphicx,amsthm}
\usepackage{fancyhdr}
\def\simto{\overset\sim\to}

\def\simleq{\underset\sim<}

\def\simgeq{\underset\sim>}

\def\simle{\underset\sim<}

\def\simge{\underset\sim>}

\def\T{\text}

\def\1#1{\overline{#1}}

\def\2#1{\widetilde{#1}}

\def\3#1{\widehat{#1}}

\def\4#1{\mathbb{#1}}

\def\5#1{\frak{#1}}

\def\6#1{{\mathcal{#1}}}

\def\C{{\4C}}

\def\R{{\4R}}

\def\N{{\4N}}

\def\Z{{\4Z}}
\def\A{\6A}
\def\M{\6M}
\def\N{\6N}
\def\F{\6F}
\def\H{\6H}
\def\S{\6S}
\def\La{\Lambda}

\def\sumK{\underset{|K|=k-1}{{\sum}'}}

\def\sumJ{\underset{|J|=k}{{\sum}'}}

\def\T{\text}
\newcommand{\Om}{\Omega}
\newcommand{\om}{\omega}

\newcommand{\bom}{\bar{\omega}}

\newcommand{\Dom}{\text{Dom} }

\newcommand{\we}{\wedge}

\newcommand{\no}[1]{\|{#1}\|}

\def\R{{\Bbb R}}

\def\C{{\Bbb C}}

\def\Z{{\Bbb Z}}

\def\la{\langle}
\def\ra{\rangle}
\def\di{\partial}
\def\dib{\bar\partial}
\def\Label#1{\label{#1}}

\def\simto{\overset\sim\to}

\def\simleq{\underset\sim<}

\def\simgeq{\underset\sim>}

\def\simle{\underset\sim<}

\def\simge{\underset\sim>}

\def\T{\text}

\def\1#1{\overline{#1}}

\def\2#1{\widetilde{#1}}

\def\3#1{\widehat{#1}}

\def\4#1{\mathbb{#1}}

\def\5#1{\frak{#1}}

\def\6#1{{\mathcal{#1}}}

\def\C{{\4C}}

\def\R{{\4R}}

\def\N{{\6N}}

\def\Z{{\4Z}}

\def\A{\6A}

\def\M{\6M}

\def\La{\Lambda}

\def\sumK{\underset{|K|=k-1}{{\sum}'}}

\def\sumJ{\underset{|J|=k}{{\sum}'}}

\numberwithin{equation}{chapter}

\def\T{\text}

\theoremstyle{plain}

\newtheorem{example}{Example}[chapter]

\newtheorem{theorem}{Theorem}[chapter]

\newtheorem{corollary}[theorem]{Corollary}

\newtheorem{lemma}[theorem]{Lemma}

\newtheorem{proposition}[theorem]{Proposition}

\theoremstyle{definition}

\newtheorem{definition}[theorem]{Definition}

\theoremstyle{remark}

\newtheorem{remark}[theorem]{Remark}

\begin{document}
\pagenumbering{roman}


\begin{titlepage}

\begin{center}
\vspace*{0.5in}
{\LARGE \bf A general method of weights  in\\ the $\bar\partial$-Neumann problem}
\par
\vspace{0.7in}
{\large  \bf Tran Vu Khanh}\\
\par 
\vspace{0.2in}
{\large Supervisor:\hskip0.2cm{ \bf  Giuseppe Zampieri} }
\end{center}
\vspace{0.8in}
\indent

A Thesis submitted for the degree of Doctor of Philosophy in front of the Committee composed by
\begin{itemize}
\item[] {\bf Joseph J. Kohn} (President)
\item[] {\bf Jeffery D. Mc.Neal}
\item[]{\bf Emil J. Straube}
\end{itemize}
\par 
\vspace{0.5in}
\begin{center}
Universit\`{a} degli Studi di Padova\\
\vspace{0.5in}
www.Math.UniPD.it\\
Dipartimento di Matematica Pura ed Applicata
\par
\vspace{0.5in}
December 2009
\end{center}

\end{titlepage}


\chapter*{Abstract}

\addcontentsline{toc}{chapter}{Abstract}


~\\

\indent This thesis deals with Partial Differential Equations in Several Complex Variables and especially focuses on a general  estimate for the $\bar\partial$-Neumann problem on a domain which is  $q$-pseudoconvex or  $q$-pseudoconcave at a  boundary point $z_0$. Generalizing Property ($P$) by \cite{C84}, we define  Property $(f\T-\M\T-P)^k$ at $z_0$.  This property yields the estimate
\begin{equation*}\label{*}
{(f\T-\M)^k} \qquad \no{f(\Lambda)\mathcal M u}^2\le c(\no{\bar\partial u}^2+\no{\bar\partial^*u}^2+\no{u}^2)+C_\M\no{u}^2_{-1}
\end{equation*}
for any  $u\in C^\infty_c(U\cap \bar{\Omega})^k\cap \T{Dom}(\dib^*)$ where $U$ is a neighborhood of $z_0$. We want to point out that under a suitable  choice of $f$ and $\M$, $(f\T-\M)^k$ is the subelliptic, superlogarithmic, compactness and subelliptic multiplier estimate. \\
\vglue9pt

The thesis also aims at exhibiting some relevant classes of domains which enjoy Property $(f\T-\M\T-P)^k$ and at discussing recent literature on the $\bar\partial$-Neumann problem in the framework of this property.

\chapter*{Sunto}
\addcontentsline{toc}{chapter}{Riassunto}

~\\

Questa tesi tratta di Equazioni alle Derivate Parziali in Pi\`{u} Variabili
Complesse e ha come obiettivo principale quello di stabilire una stima
generale per il problema $\bar\partial$-Neumann su un dominio che \`{e}
$q$-pseudoconvesso o $q$-pseudoconcavo in corrispondenza di un punto di
bordo $z_0$. Generalizzando la Propriet\`{a} $(P)$ di \cite{C84}, si introduce la
Propriet\`{a} $(f\T-\M\T-P)^k$ in $z_0$. Essa d\`{a} luogo alla stima
\begin{equation*}
{(f\T-\M)^k} \qquad \no{f(\Lambda)\mathcal M u}^2\le c(\no{\bar\partial
u}^2+\no{\bar\partial^*u}^2+\no{u}^2)+C_\M\no{u}^2_{-1}
\end{equation*}
per ogni  $u\in C^\infty_c(U\cap \bar{\Omega})^k\cap \T{Dom}(\dib^*)$ ove
$U$ \`{e} un intorno di $z_0$. \`{E} il caso di osservare che per opportune
scelte di $f$ e di $\M$, la stima $(f\T-\M\T-P)^k$ coincide con le
principali stime della letteratura quali quelle subellittiche,
superlogaritmiche, di compattezza e infine quelle di moltiplicatore
subellittico.
\vglue9pt
La tesi ha anche l'obiettivo di esibire delle classi rilevanti di domini
che godono della Propriet\`{a} $(f\T-\M\T-P)^k$ e di discutere letteratura
recente sul problema $\bar\partial$-Neumann nel quadro di questa
propriet\`{a}.


\chapter*{Acknowledgements}~\\
\addcontentsline{toc}{chapter}{Acknowledgements}
~\\

When I started the study of the $\bar\partial$-Neumann problem, I would
never have dreamed of having Joseph J. Kohn as President of my Ph.D. Committee. I
knew from the beginning his extraordinary reputation in our community. What I
learned along the way is the continuous presence of his scientific
discoveries in my research work. He is really a master to me.
\\

I am indebted to Emil J. Straube whose research has recently attracted my
 attention. I  have also profited   from his advice and the fruitful
exchanges we have had.
\\

 Finally, Jeffery D. McNeal's work has greatly excited my interest. I thank him
for joining the Committee.\\

\chapter*{Further Acknowledgements}

~\\

 It is a great pleasure to thank my advisor, Giuseppe Zampieri, for his discussions and guidance in directing my research. I am deeply indebted to him for initiating me to the field  of Complex Analysis and CR Geometry. I sincerely appreciate  his support over the years and I feel immensely lucky to have had such an advisor. \\

I would like to thank  Luca Baracco for many fruitful discussions.
I am also indebted to Dmitri Zaitsev for his advice.
\\

 I would like to dedicate my special thanks to the Foundation CARIPARO which offered me a three years scholarship.  In addition, I would like to thank all professors of the Department of Mathematics of the University of Padova for their
interesting lectures.\\

I would like to thank all  friends who have supported me during the time I spent in the Graduate School. \\

Finally, I would like to thank my family for all the support and encouragement they have given me during these past three years. I would especially like to thank my wife, Uyen Phuong. 

\clearpage
\tableofcontents 
\addcontentsline{toc}{chapter}{Contents}

\pagestyle{fancy}
\renewcommand{\chaptermark}[1]{\markboth{Chapter \thechapter. #1}{}}
\renewcommand{\sectionmark}[1]{\markright{ #1}}
\fancyhf{}
\fancyhead[RO]{\bfseries\leftmark}
\fancyhead[LE]{ \bfseries\rightmark}
\fancypagestyle{plain}{%
 \fancyhead{} 
 \renewcommand{\headrulewidth}{0pt} 
}

\addtolength{\headheight}{3pt}
\cfoot{\bf \thepage}


\clearpage
\pagenumbering{arabic}



\chapter{Introduction} 
 The $\bar\partial$-Neumann problem is probably the most important and natural example of a non-elliptic boundary value problem.
It owes its importance to the Cauchy-Riemann system from which it originates. The main tools to prove  regularity of solutions of this problem are various $L_2$-estimates such as subelliptic, superlogarithmic, compactness estimates just to mention a few. In this thesis, we introduce a general estimate, that we call $(f\T-\M)^k$, whose principal related results are Theorems \ref{main1}, \ref{main2} and \ref{main3}.  As an introduction, we give a brief description of the $\dib$-Neumann problem, for whose detailed account we refer to \cite{FK72}.

\section{The $\dib$-Neumann problem}
\indent

 Let $\Omega$  be a bounded domain of $\C^n$ with smooth boundary denoted by $b\Om$.  Let $L_2^{h,k}(\Om)$ \index{Space ! $L_2^{h,k}(\Om)$} be the space of square-integrable $(h, k)$-forms on $\Om$. We have a densely defined complex of operators $\bar\partial$ with adjoint $\bar\partial^*$ \index{Notation !  $\dib^*$}\index{Notation ! $\dib$}
\begin{equation}\begin{split}
L_2^{h,k-1}(\Om)\underset{\dib^*}{\overset{\dib}{\rightleftarrows}} L_2^{h,k}(\Om)\underset{\dib^*}{\overset{\dib}{\rightleftarrows}} L_2^{h,k+1}(\Om).
 \end{split}\end{equation}
 We define $\H^{h,k}\subset L^{h,k}_2(\Om)$\index{Space ! $\H^{h,k}$} by
\begin{equation}\begin{split}
\H^{h,k}=\{u\in \T{Dom}(\dib)\cap\T{Dom}(\dib^*)\big| \dib u=0 \T{ ~~and~~ } \dib^*u=0 \}.
 \end{split}\end{equation}
The {\it  $\dib$-Neumann problem}\index{Problem !$\dib$-Neumann }  for $(h,k)$-forms can  be stated as follows: given $\alpha\in L^{h,k}_2(\Om) $ with $\alpha \perp\H^{h,k}$, does there exist $u\in L^{h,k}(\Om)$ such that
\begin{equation}\begin{split}
\Label{dN}
\begin{cases}(\bar{\partial}\bar{\partial}^*+\bar{\partial}^*\bar{\partial})u=\alpha\\
u\in \T{Dom}(\bar{\partial})\cap \T{Dom}(\bar{\partial}^*)
\\
\bar{\partial} u\in \T{Dom}(\bar{\partial}^*), \bar{\partial}^* u\in \T{Dom}(\bar{\partial})?
\end{cases}
\end{split}\end{equation}

Observe that the solution $u$ to \eqref{dN} under the constraint $u\perp \mathcal H^{h,k}$ is unique, if it exists. We denote  this solution by $N\alpha$. If a solution to \eqref{dN} exists for all $\alpha \perp\H^{h,k}$, then we extend the operator $N$\index{Operator ! $\dib$-Neumann, $N$, $N_k$} to a linear operator on $L_2^{h,k}(\Om)$ by setting 
\begin{equation}\begin{split}
N\alpha =\begin{cases} 0 & \T{ if }  \alpha \in \H^{h,k}\\
 u  & \T{ if }  \alpha \perp\H^{h,k}.
\end{cases}
\end{split}\end{equation}
Then $N$ is bounded self adjoint. Furthermore, if $\dib \alpha =0$, then from \eqref{dN} we obtain $\dib\dib^*\dib N\alpha=0$ and taking inner product with $\dib N\alpha$ we get $\no{\dib^*\dib N\alpha}^2=0$ and hence $\dib^*\dib N\alpha=0$. Thus we see from \eqref{dN} that if $\dib\alpha=0 $ and  $\alpha \perp\H^{h,k}$, then $\alpha=\dib\dib^* N\alpha$. It then follows that $u=\dib^* N\alpha$ is the unique solution to the {\it $\dib$-problem} \index{Problem ! $\dib$}
\begin{equation}\begin{split}
\Label{dE}
\begin{cases}
\bar\partial u=\alpha,
\\
u\T{ is orthogonal to $\T{Ker}\, \bar\partial$}.
\end{cases}
\end{split}\end{equation}

The $\dib$-Neumann problem is a non-elliptic boundary value problem; in fact,   the Laplacian  $\Box=\bar\partial\bar\partial ^*+\bar\partial^*\bar\partial$ \index{Notation ! $\Box$} itself is elliptic but the boundary conditions which are imposed by the membership to $\T{Dom}(\Box)$,\index{Notation ! $\T{Dom}(\Box)$} that is, the second and third line of \eqref{dN}, are not.  The main interest relies in the {\it regularity} at the boundary for this problem, that is, in stating under which condition $u$ inherits from $\alpha$ the smoothness at the boundary $b\Omega$ (it certainly does in the interior). The regularity of the $\dib$-Neumann operator is defined as follows
\begin{definition}\Label{re}
~~
\begin{enumerate}
  \item[(i)] {\it Global regularity :}\index{ Regularity ! global} if $\alpha\in C^\infty(\bar\Om)$ then $N\alpha\in C^\infty(\bar\Om)$.
  \item[(ii)] {\it Local regularity : }\index{ Regularity ! local} if  $\alpha\in C_c^\infty(U\cap  \bar\Om)$ then $N\alpha\in C^\infty(U'\cap \bar\Om)$ where the sets $U'\subset U$ range in a system of  neighborhoods of a point $z_0\in \bar\Om$.
\end{enumerate}  
\end{definition}
One of the main tools used in investigating the local (resp. global ) regularity at the boundary of the solutions of the $\dib$-Neumann problem  consists  in   certain a priori estimates such as subelliptic, superlogarithmic (resp. compactness) estimates. In this thesis we introduce the $(f\T-\M)^k$ estimate which is comprehensive of all these estimates.
 
\section{The $(f\T-\M)^k$ estimate}
\indent 

In order  to define  $(f\T-\M)^k$, some preliminary material is required.  \\

For quantities  $A$ and $B$ we use the notion $A \lesssim B$\index{Notation ! $\lesssim$}\index{Notation ! $\cong$} to mean $A\le cB$ for some constant $c>0$ which is independent of relevant parameters. We write $A\cong B$  to mean $A\lesssim B$ and $B\lesssim A$. For functions $f$ and $g$, we use the notation $f\gg g$ \index{Notation ! $\gg$} to mean that $\frac fg$ is $+\infty$ at $+\infty$.\\
    
Let $\Omega$  be a smooth domain locally defined in a neighborhood of  boundary point $z_0$. Throughout this thesis we assume that $z_0$  is the origin.  For a neighborhood  $U$ of $z_0$, fix a smooth real-valued function $r$ such that 
\begin{equation}\begin{split}
\Om\cap U=\{z\in U: r(z)< 0\}; 
\end{split}\end{equation}
we also require  that $ |\di r|=1$ on $b\Om$. We take a local orthonormal basis of $(1,0)$ forms $\om_1,...,\om_n=\di r$ \index{Notation ! $\om_1,...,\om_n$} and a dual basis of (1,0) vector fields $L_1,...,L_n$;\index{Notation ! $L_1,...,L_n$} thus $L_1,...,L_{n-1}$ generate $T^{1,0}(U\cap b\Om)$. For $\phi \in C^2(U)$, we denote by $\phi_{ij}$ the coefficients of $\di\dib\phi$ in this basis.\index{Notation ! $\phi_{ij}$}\\

 Let $\lambda_1(z)\le ...\le \lambda_{n-1}(z)$ be the eigenvalues of $(r_{jk}(z))^{n-1}_{j,k=1}$. We take a pair of indices $1\le q\le n-1$ and $0\le q_o\le n-1$ such that $q\not=q_o$. We assume that there is a bundle $\mathcal V^{q_o}\in T^{1,0}b\Om$  of rank  $q_o$ with smooth coefficients (that by reordering can be supposed to coincide with $\T{ span }\{ L_1,...,L_{q_o}\}$) such that \index $q_o$
\begin{equation}\begin{split}
\Label{qcon}
\sum^q_{j=1}\lambda_j-\sum_{j=1}^{q_o}r_{jj}\ge 0 \T{ on } U\cap b\Om.
\end{split}\end{equation}
Here we conventionally set $\sum_{j=1}^{q_o}\cdot \equiv 0$ if $q_o=0$.
\begin{definition}
\Label{dqcon}~
\begin{enumerate}
\item[(i)] If $q>q_o$ we say that $\Omega$ is {\it $q$-pseudoconvex } at $z_0$. \index{Domain ! $q$-pseudoconvex}
\item[(ii)] If $q<q_o$ we say that $\Omega$ is {\it $q$-pseudoconcave }  at $z_0$.\index{Domain ! $q$-pseudoconcave}
\end{enumerate}
The $q$-pseudoconvexity/concavity is said to be {\it strong } when \eqref{qcon} holds as strict inequality. 
\end{definition}

\noindent

The notion of $q$-pseudoconvexity was used in \cite{A07} and \cite{Z00}  to prove the existence of $C^\infty(\bar \Omega)^k$ solutions to the equation $\bar\partial u=f$.  Though the notion of $q$-pseudoconcavity is formally symmetric to $q$-pseudoconvexity, it is useless in the existence problem. The reason is intrinsic. Existence is a  ``global" problem
but  bounded domains are never globally $q$-pseudoconcave. Instead, in the present thesis, we deal 
with various  estimates of local type such as subelliptic or superlogarithmic.
Owing to the local nature of these estimates and the related local regularity  of $\dib$-Neumann problem, this is the first occurrence where $q$-pseudoconcavity comes successfully into play. Moreover,  local estimates on a pseudoconcave domain play a leading role in the study the $\dib$ or $\dib_b$-Neumann problem on an annulus or a hypersurface. \\

\begin{remark}{\rm
If $\Om$ is $q$-pseudoconvex at $z_0$, then  $\Om$ is also $k$-pseudoconvex for any $k\ge q$. Similarly,   if $\Om$ is $q$-pseudoconcave at $z_0$, then  $\Om$ is also $k$-pseudoconcave for any $k\le q$. \\}\end{remark}

\begin{remark}{\rm
Definition \ref{dqcon} is a generalization of the  usual pseudoconvexity or pseudoconcavity as well as of  the celebrated condition $Z(q)$ \index{Domain ! condition $Z(q)$} (cf. \cite{FK72}). In fact, for $q_o=0$ and $q=1$, 1-pseudoconvexity is the usual pseudoconvexity. Similarly, for $q_o=n-1$ and $q=n-2$, $(n-1)$-pseudoconcavity is the usual pseudoconcavity. Moreover, if $\Om$ satisfies condition $Z(q)$ at $z_0$, that is, the Levi form has at least $n-q$ positive eigenvalues or at least $q+1$ negative eigenvalues,  then $\Om$ is strongly $q$-pseudoconvex or strongly $q$-pseudoconcave at $z_0$.\\}\index{Domain ! pseudoconvex}\index{Domain ! pseudoconcave}
\end{remark}

 We denote  by $\A^{h,k}$ \index{Space ! $\A^{h,k}$} the space of smooth $(h,k)$-forms in $\bar\Om$. Throughout this thesis we only deal with $(0,k)$-forms since the extension from type $(0,k)$ to type $(h,k)$ is trivial. We use the notation \index{Space ! $C^\infty_c(U\cap \bar\Om)^{k}$}
$$C^\infty_c(U\cap \bar\Om)^{k}=\{u\in \A^{0,k}\big | \T{supp}(u)\subset U\}.$$
 If $u\in C^\infty_c(U\cap \bar\Om)^{k}$, then $u$ can be written as 
\begin{equation}\Label{u}
u=\sumJ u_J\bom_J,
\end{equation}
where $\bom_J=\bom_{j_1}\we...\we\bom_{j_k}$ and where $J=\{j_1,...,j_k\}$ is a ordered multiindex. Also, $\sum'$ denotes summation over strictly increasing index sets. If $J$ decomposes as $J=jK$, then we write $u_{jK}=\epsilon_{J}^{jK}u_J$ where $\epsilon_{J}^{jK}u_J$ is the sign of the permutation $(iK)\overset\sim\to J$. 
Under the choice of a basis of $(1,0)$-forms in which $\omega_n=\partial r$, we
 check readily that  $u\in \T{Dom}(\bar\partial^*)$ if and only if  $u_{nK}|_{b\Om}=0$ for any $K$. \\

We define the {\it function multiplier} \index{Multiplier ! function}
$ \M\in \A^{0,0}$, by
\begin{equation}\Label{Mu1}
\M u := |\M| | u|.
\end{equation}
We define the {\it $1$-form multiplier} \index{Multiplier ! $1$-form}$\M=\sum_j \M_j \bar\om_j\in \A^{0,1}$ \index{Notation ! $\M$} over a  q-pseudoconvex/concave domain $\Om$ at $z_0$ by
\begin{equation}\Label{Mu2}
\M u :=\sqrt{\Big |\sumK |\sum_{j=1}^{n} \M_j u_{jK}|^2-\sum_{j=1}^{q_o}|\M_j|^2 | u|^2\Big|}.
\end{equation}

For $z_0\in b\Om$, we choose the defining function $r$ as the last coordinate and supplement it to a full system $(x,r)\in\R^{2n-1}\times \R$ of so-called {\it special boundary coordinates} \index{Coordinate(s) ! special boundary}  in a neigborhood $U$ of $z_0$. In this situation, the  $x_j$'s  are  the {\it tangential coordinates}\index{Coordinate(s) ! tangential } and $r$ the {\it normal coordinate}\index{Coordinate(s) ! normal }.  Denote by $\xi$ the dual variables of $x$   and define $x\cdot\xi=\sum x_i\xi_j $ and $|\xi|^2=\sum\xi_j^2.$\\

For $\varphi \in C_c^\infty(U\cap\bar\Om)$, we define $\tilde{\varphi }$, the {\it tangential Fourier transform}\index{Notation ! $\tilde{\phi}$} of $\varphi$, by 
$$\tilde{\varphi }(\xi, r)=\int_{\mathbb{R}^{2n-1}}e^{-i x\cdot \xi}\varphi (x,r)dx.$$
We denote by $\Lambda_\xi=(1+|\xi|^2)^{\frac{1}{2}}$ \index{Notation ! $\Lambda_\xi$} the standard ``tangential" elliptic symbol of order 1 and by $\Lambda$ \index{Notation ! $\Lambda$} the operator with symbol $\Lambda_\xi$\index{Operator ! tangential,  $f(\La)$}. We define a class of functions  $\mathcal F$ \index{Notation ! $\mathcal F$} by
$$\mathcal F=\{f\in C^\infty([1,+\infty)) \Big| 0\le f(t)\lesssim t^{\frac{1}{2}}; f'(t)\ge 0 \T{ and } \big |f^{(m)}(t) \big |\lesssim \big |\frac{f^{(m-1)}(t)}{t} \big |, \forall m\in \Z^+ \}.$$
Notice that  for any of the choices $f(t)=1$, $f=(\log t)^s$ or  $f=t^\epsilon$, $\epsilon\le \frac{1}{2}$, we have that $f\in \F$.\\

For $f\in \mathcal F$,  we define the operator  $f(\Lambda)$ by \index{Notation ! $f(\La)$}
\begin{equation}
f(\Lambda)\varphi(x,r) =(2\pi)^{-2n+1}\int_{\R^{2n-1}}e^{ix\cdot \xi} f(\Lambda_\xi)\tilde{\varphi }(\xi,r)d\xi,
\end{equation}
where $\varphi\in C_c^\infty(U\cap \bar\Om)$.\\

We define the energy form $Q$ \index{Notation ! $Q$} on $ C_c^\infty(U\cap \bar\Om)^k\cap \T{Dom}(\dib^*)$ by 
$$Q(u,v)=(\bar{\partial} u,\bar{\partial}v)+(\bar{\partial}^* u,\bar{\partial}^* v)+(u,v).$$
We are now ready to define the $(f\T-\M)^k$ estimate.

\begin{definition}\Label{fMk}Let $\Om$ be $q$-pseudoconvex (resp. $q$-pseudoconcave) at $z_0\in b\Om $.  Then,   the $\bar{\partial}$-Neumann problem is said to satisfy the {\it  $(f\T-\M)^k$ estimate} \index{Estimate ! $(f\T-M)^k$}
 at  $z_0\in b\Om$  if there exist a positive constant $C_\M$   and a neighborhood $U$ of $z_0$ such that 
\begin{equation}\begin{split}
(f\T-\M)^k\qquad ||f(\Lambda)\M u||^2 \lesssim Q(u,u)+C_{\mathcal M}\no{u}^2_{-1}, 
\end{split}\end{equation}
for any $u\in C^\infty_c(U\cap \bar{\Omega})^{k}\cap \T{Dom}(\bar{\partial}^*)$ where $k\ge q$ (resp. $k\le q$). 
\end{definition}

\begin{remark} If $\M$ is a bounded constant, then we can assume that $C_\M=0$.
\end{remark}

\begin{remark} If $\Om$ is pseudoconvex at $z_0$ and $(f\T-\M)^{k_o}$ holds for some $k_o\ge 1$,   then $(f\T-\M)^k$ holds for any $k\ge k_o$. Similarly, If $\Om$ is pseudoconcave at $z_0$ and $(f\T-\M)^{k_o}$ holds for some $k_o\le n-2$, then $(f\T-\M)^k$ holds for any $k\le k_o$.
\end{remark}

We  simply say that $(f\T-\M)^k$ holds when the definition applies.  We remark that  if $f(|\xi|) \cong |\xi|^{\frac{1}{2}}$, then the operator $\M$ has to be bounded  since the best estimate at the boundary is the $\frac{1}{2}$-subelliptic estimate.    
 
\section{How the $(f\T-\M)^k$ estimate mathces the literature}
\indent

We want to point our attention to the choice  of $f$ and  $\M$ in relevant cases and review some results concerning these estimates.\\

1) For $f(|\xi|)=|\xi|^\epsilon$, $0<\epsilon\le \frac{1}{2}$, and for $\M=1$, the $(f\T-
\M)^k$ estimate  becomes the {\it subellipitic estimate} \index{Estimate ! subellipic}
\begin{equation}\begin{split}
\Label{sub}
||| u|||_\epsilon^2 \lesssim Q(u,u).
\end{split}\end{equation}
Here $|||\cdot|||_\epsilon$ \index{Notation ! $|||\cdot|||_\epsilon$} is the {\it tangential} $\epsilon$-Sobolev norm.
When the domain $\Omega$ is pseudoconvex, a great deal of work has been done about subelliptic estimates. The most general results concerning this problem have been  obtained in  \cite{K79} and  \cite{C87}. 
\begin{enumerate} 
  \item[$\bullet$] 
In \cite{K79}, Kohn gave a sufficient condition for subellipticity over pseudoconvex domains with real analytic boundary by introducing a sequence of ideals of subelliptic multipliers.
\item[$\bullet$]In \cite{C87}, Catlin proved, regardless whether $b\Omega$ is real analytic or not, that subelliptic estimates hold  for $k$-forms at  $z_0$ if and only if the D'Angelo type $D_k(z_0)$ is finite.  Catlin applies the method of weight functions used earlier by H\"ormander \cite{H66}. One step in Catlin's proof is the following reduction:
\begin{theorem} \Label{t1.2} Suppose  that $\Om\subset\subset  \C^n$   is a pseudoconvex  domain  defined  by $\Om = \{r < 0\}$, let $z_0 \in  b\Om$, and  let $U$ be a neighborhood of $z_0$.  Suppose that for all $\delta > 0$, there is a smooth real-valued function $\Phi^\delta$  satisfying the properties:
\begin{equation}\begin{split}
\Label{CatP}
\begin{cases}
|\Phi^\delta| \le 1 \T{ on }  U,\\
\Phi^\delta_{ij}   \ge  0      \T{ on } U ,\\
\sum_{ij=1}^{n}\Phi^\delta_{ij}u_i\bar u_j\simge \delta^{-2\epsilon}|u|^2 \T{ on } U\cap \{-\delta<r\le 0\},
\end{cases}
\end{split}\end{equation}
where the matrix $\Phi^\delta_{ij}$ represents the form $\partial\bar\partial \Phi^\delta$ in the chosen basis $\{\omega_j\}_j$.
Then, we have a subelliptic estimate of order $\epsilon$ at $z_0$.
\end{theorem}
 
\end{enumerate}

However, not much is known in the case when the domain is not necessarily pseudoconvex except from the results related to the celebrated $Z(k)$ condition which characterizes the existence of subelliptic estimates of index $\epsilon=\frac12$ 
according to H\"ormander \cite{H65} and Folland-Kohn \cite{FK72}. Some further results, mainly related to the case of forms of top degree $n-1$, are due to Ho \cite{Ho85}. \\

The basic theorem of Kohn and Nirenberg \cite{KN65} shows that local regularity is a consequence of a subelliptic estimate. In fact, if a subelliptic estimate of order $\epsilon$ holds for the $\bar\partial$-Neumann problem on a neighborhood $U$ of $z_0\in b\Omega$, then $\alpha|_U\in H_s(U)^k$ implies $N\alpha|_{U'}\in H_{s+2\epsilon}(U')^k$ for $U'\subset\subset U$; here $H_s(U)^k$ denotes the space of k-forms with coefficients in the $H_s$-Sobolev space of order $s$. \index{Space ! $H_s(U)^k$} \\

2) For $f(|\xi|)=\log|\xi|$ and $\M=\frac{1}{\epsilon}$ for any $\epsilon>0$, we have that  $(f\T-\M)^k$  implies the {\it superlogarithmic estimate } \index{Estimate ! superlogarithmic}
\begin{equation}\begin{split}
\Label{suplog}
 ||\log(\Lambda)u||^2 \lesssim \epsilon( \no{\bar{\partial} u}^2+\no{\bar{\partial}^* u}^2)+C_\epsilon \no{u}^2.
\end{split}\end{equation}
Superlogarithmic estimates were first introduced by Kohn in \cite{K02}. He proved  superlogarithmic estimates for the tangential Kohn-Laplace operator $\Box_b$  on certain pseudoconvex CR manifolds and using them he established local regularity of $\Box_b$ and of the $\bar\partial$-Neumann problem. These estimates are established under the assumption that subellipticity degenerates in certain specified ways.\\

3) For $f(|\xi|)\equiv 1$ and $\M=\frac{1}{\epsilon}$ for any $\epsilon>0$, we have that $(f\T-\M)^k$ implies the {\it compactness estimate } \index{Estimate ! compactness }
\begin{equation}\begin{split}
\Label{compact}
 ||u||^2\lesssim \epsilon( \no{\bar{\partial} u}^2+\no{\bar{\partial}^* u}^2)+C_\epsilon \no{u}_{-1}^2.
\end{split}\end{equation} 
By its definition, the compactness estimate in  \eqref{compact} is a local property. Roughly speaking, the $\bar\partial$-Neumann operator $N$ on $\Omega$ is compact if and only if every boundary point 
has a neighborhood $U$ such that the corresponding $\bar\partial$-Neumann operator on $U\cap \Omega$ is compact. A classical theorem of Kohn and Nirenberg \cite{KN65} asserts that compactness of $N$ (as an operator from $L_2(\Omega)$ into itself) implies global regularity in the sense of preservation of Sobolev regularity.\\

Catlin \cite{C84} introduced Property (P) \index{Property ! (P)} and showed that it implies a compactness estimate for the $\bar\partial$-Neumann problem. A pseudoconvex domain $\Omega$ has Property (P) if for every positive number $M$ there exists a plurisubhamonic function $\Phi^M$ in $C^\infty(\Omega)$, bounded between 0 and 1, whose complex Hessian has all its eigenvalues bounded from below by $M$ on $b\Omega$, that is
\begin{equation}\begin{split}
\Label{Pcom}
\sum_{ij=1}^{n}\frac{\partial^2\Phi^M}{\partial z_i\partial \bar{z}_j}w_i\bar{w}_j\ge M|w|^2,  \T{   for } z\in b\Omega, w\in \C^n.
\end{split}\end{equation}
Compactness  is  completely  understood  on  (bounded)  locally  convexifiable  domains. On such domains, the following are equivalent \cite{FS98}, \cite{FS01} :
\begin{enumerate}
\item[(i)]  $N_k$  is compact,
\item[(ii)] the boundary $b\Omega$ satisfies property (P),
\item[(iii)]  the boundary contains no $k$-dimensional analytic variety.
\end{enumerate}
In general, however, the situation is not understood at all.\\

5) For $f(|\xi|)\equiv 1$ and $\M=\frac{1}{\epsilon}\sum_{j=1}^{n} \overline{r^\epsilon_{z_i\bar z_j}r_{\bar z_i}} d\bar z_j$ where $r^\epsilon$ and $ r$ are defining functions with $|\nabla r^\epsilon|\cong 1$ on $b\Om$, we have that $(f\T-\M)^k$ is equivalent to 
\begin{eqnarray}
\Label{1.8}
\sumK\no{ \sum_{i,j=1}^{n}\overline {r^\epsilon_{z_i\bar z_j}r_{\bar z_i}}u_{jK}}^2\lesssim \epsilon ( \no{\bar{\partial} u}^2+\no{\bar{\partial}^* u}^2)+C_\epsilon\no{u}_{-1}^2,
\end{eqnarray}
for any $u=\sumJ u_Jd\bar z_J\in C^\infty(\bar\Om)^k\cap \T{Dom}(\dib^*)$.\\

 The estimate \eqref{1.8} was introduced by Straube in \cite{S08}. He showed that if \eqref{1.8} holds for all $u\in C^\infty(\Omega)^k\cap \T{Dom}(\bar\partial^*)$, then  the $\bar\partial$-Neumann operator $N_k$  on k-forms is exactly regular in Sobolev norms, that is
$$\no{N_ku}_s\le C_s\no{u}_s$$
for any integer $s\ge 0$ and any $u\in H_s(\Omega)^k$. Notice that  the estimate \eqref{1.8} is weaker than the compactness estimate. \\

6) For $f(|\xi|)=|\xi|^\epsilon$  and $\M\in  \A^{0,0}$ (resp.  $\M\in  \A^{1,0}$), we have that the  $(f\T-\M)^k$-estimate can be written as 
\begin{equation}\begin{split}
\Label{subm1}
||| \M u|||_\epsilon^2 \lesssim Q(u,u).
\end{split}
\end{equation}

One calls $\M$  a {\it subelliptic multiplier} (resp. {\it subelliptic vector-multiplier}). Kohn \cite{K79} introduced the theory of subelliptic multipliers and subelliptic vector-multipliers to obtain  subelliptic estimates on a real analytic domain. To prove subelliptic estimates, he searched for a nonzero constant function belonging to the collection of  multipliers.

\section{The main theorems}
\indent

The first goal of this thesis, we exploit here the full strength of Catlin's method to study  the $(f\T-\M)^k$ estimate on a $q$-pseudoconvex or $q$-pseudoconcave domain. These results are related to  joint work with G. Zampieri in \cite{KZ08a},\cite{KZ08b} and \cite{KZ09a}.\\

For a form $u$, we consider the decomposition $u=(u^\tau,u^\nu)$ where the {\it tangential component} $u^\tau$ collects the coefficients $u_J$ with $n\notin J$ and the {\it normal component} $u^\nu$ those with $n\in J$. \index{Notation ! $u^\tau$, $u^\nu$}
Denote by $S_\delta$ the strip  $\{z\in \Om|-\delta<r<0\}$. Generalizing  conditions \eqref{CatP} in Theorem \ref{t1.2} and Property $(P)$  in \eqref{Pcom}, we define 
\index{Notation ! $S_\delta$}
\begin{definition}\Label{fMP} We say that $\Omega$ satisfies {\it Property $(f\T-\M\T-P)^k$}  \index{Property ! $(f\T-\M\T-P)^k$} at the boundary point $z_0\in b\Omega$, if there is a neighborhood $U$ of $z_0$  and 
a family of real valued $C^2(U)$ weights 
 $\Phi:=\Phi^{\delta,\M}$   such that we have on $S_\delta\cap U$
\begin{eqnarray}(f\T-\M\T-P)^k\qquad
\begin{cases}
|\Phi|\lesssim 1  &    \\
\sumK\underset{ij=1}{\overset{n}\sum}\Phi_{ij}u^\tau_{iK}\bar u^\tau_{jK}&-\underset{j=1}{\overset{q_o}\sum}\Phi_{jj} |u^\tau|^2\\
&\simge f(\delta^{-1})^2|\mathcal M u^\tau|^2 +\underset{j=1}{\overset{q_o}\sum} |L_j(\Phi)u^\tau|^2,
\end{cases}
\end{eqnarray}
for any $u\in C_c^\infty(U\cap\bar\Om)^k$.
\end{definition}

Our result is the following    
\begin{theorem}\Label{main1}Let $\Omega\subset \C^n$ be q-pseudoconvex (resp. q-pseudoconcave) and satisfy Property $(f\T-\M\T-P)^k$ at  $z_0\in b\Omega$; then the $(f\T-\M)^k$-estimate
holds at $z_0$ for $k\ge q$ (resp. $k\le q$).
\end{theorem}

Observe that the condition on the family  $\{\Phi^{\delta,\M}\}$ in our Property $(f\T-\mathcal M\T-P)^k$ is simpler than condition \eqref{CatP} in Theorem \ref{t1.2} and Property $(P)$  in \eqref{Pcom}.\\

The main idea for proving Theorem \ref{main1} follows from  Catlin \cite{C87} combined with some modifications contained in \cite{KZ08a} and \cite{KZ08b}. \\

We remark that our Property $(f\T-\M\T-P)^k$ only deals with the $u^\tau$-tangential component of $u$.  We first get the $(f\T-\M)^k$ estimate for $u$ replaced by $u^\tau$.  On the other hand, for the normal component $u^\nu$ of $u$,  one readily proves the elliptic estimate $\no{u^\nu}_1^2\lesssim Q(u,u)$. Thus the full $u$ enjoys the $(f\T-\M)^k$ estimate.\\ 

On a real hypersurface $M$ of $\C^n$, $\dib$ induces the {\it tangential Cauchy-Riemann} operator $\dib_b$. Let $\dib^*_b$ \index{Operator ! tangential Cauchy-Riemann, $\dib_b$} be  the $L_2$-adjoint of $\dib_b$ and $\Box_b=\dib_b\dib_b^*+\dib_b^*\dib_b$, the induced Kohn-Laplacian. \index{Notation ! $\Box_b$} \index{Notation ! $\dib_b$, $\dib_b^*$} 
Let $(\cdot,\cdot)_b$ denote the inner $L_2$-product on $M$; for $u,\,v\in C^\infty_c(U\cap M)^k$, we denote the tangential energy by
$$Q_b(u,v)=(\dib u,\dib u)_b+(\dib^* u,\dib^*v)_b+(u,v)_b.$$\index{Notation ! $Q_b$}
In the space of frequencies $\xi\in\R^{2n-1}$, we consider a conical smooth partition of the unity 
$$
\psi^+(\xi)+\psi^-(\xi)+\psi^0(\xi)=1,
$$
where $\T{supp}\,\psi^\pm\subset\{\xi :\,\pm\xi_{2n-1}>|(\xi_1,...,\xi_{2n-2})|\}$ and $\T{supp}\,\psi^0\subset\{\xi : \,|\xi_{2n-1}|<2|(\xi_1,...,\xi_{2n-2})|\}$. We consider the corresponding pseudodifferential decomposition of the identity
$$
\Psi^++\Psi^-+\Psi^0=\T{id}.\index{Notation ! $\Psi^+$, $\Psi^-$, $\Psi^0$ }
$$ \index{Notation ! $\psi^+$, $\psi^-$, $\psi^0$ }
Accordingly, for a form $u\in C^\infty_c(U\cap M)^k$, we consider the decomposition
\begin{equation}
\Label{decomposition}
\begin{split}
u&=\zeta\Psi^+u+\zeta\Psi^-u+\zeta\Psi^0u
\\
&=:u^++u^-+u^0,
\end{split}
\end{equation}\index{Notation ! $u^+$, $u^-$, $u^0$} 
where $\zeta$ is a cut off in $C^\infty_c( U'\cap M)$ for $U'\supset\supset U$. 
Our estimates on $M$ are defined as follows 
 \begin{definition} Let $M$ be a hypersurface; then a $(f\T-\M)^k_{b}$ estimate holds for $(\dib_b,\dib_b^*)$  at $z_0\in M$ if, in a neighborhood $U$ of $z_0$, we have
$$(f\T-\M)^k_{b}  \qquad\no{f(\Lambda)\M u}^2_b\lesssim Q_b(u,u)+C_\M\no{u}_{b,-1}^2,$$ \index{Estimate ! $(f\T-\M)^k_b$}
for all $u\in C_c^\infty(U\cap M)^{k}$. And  a $(f\T-\M)^k_{b,+}$ (resp. $(f\T-\M)^k_{b,-}$)  estimate holds for $(\dib_b,\dib_b^*)$  at $z_0$ if the above holds with $u$ replaced by $u^+$ (resp.  $u^-$), that is,  
$$(f\T-\M)^k_{b,+}\qquad \no{f(\Lambda)\M u^+}^2_b\lesssim Q_b(u^+,u^+)+C_\M\no{u^+}_{b,-1}^2$$  \index{Estimate ! $(f\T-\M)^k_{b,+}$}
(resp.
$$(f\T-\M)^k_{b,-}\qquad \no{f(\Lambda)\M u^-}^2_b\lesssim Q_b(u^-,u^-)+C_\M\no{u^-}_{b,-1}^2\quad ).$$  \index{Estimate ! $(f\T-\M)^k_{b,-}$}
 \end{definition}

\begin{definition}
The hypersurface $M$ is said to be $q$-pseudoconvex  if one of two sides $\C^n\setminus M$ is $q$-pseudoconvex. 
\end{definition}

We denote by $\Om^+=\{z\in U | r(z)<0 \}$ the $q$-pseudoconvex side of $\C^n\setminus M$ and  by $\Om^-$ the complementary side. By Remark \ref{convex-concave} which follows, $\Om^-=\{z\in U | -r(z)<0 \}$ is $(n-q-1)$-pseudoconcave. We use the notations   $(f\T-\M)^k_{\Om^\pm}$ in the obvious sense.\\

The second goal of this thesis is to state the equivalence of the  $(f\T-\M)^k$ estimate on a domain and on its boundary.
\begin{theorem}\Label{main2}
Let $M$ be a $q$-pseudoconvex hypersurface at $z_0$
and let $\Omega^\pm$ be  the two sides of $\C^n\setminus M$ locally at $z_0\in M$. 
  Then, for any $k\geq q$, we have the following system of equivalences
\begin{eqnarray}
(f\T-\M)^k_{\Om^+} \Longleftrightarrow(f\T-\M)^k_{b,+}\Longleftrightarrow (f\T-\M)^{n-1-k}_{b,-}\Longleftrightarrow (f\T-\M)^{n-1-k}_{\Om^-}.
\end{eqnarray}
\end{theorem}

Since $Q_b(u^0,u^0)\simge \no{u^0}_{b,1}^2$, Theorem \ref{main2} implies that $(f\T-\M)_b^k$  holds on $M$ at $z_0$ if  Property $(f\T-\M\T-P)^k$ and Property $(f\T-\M\T-P)^{n-1-k}$ hold on $\Om^+$ at $z_0$.   \\

The third goal of this thesis is the analysis of local problems which arise, on a bounded pseudoconvex domain,  from the $(f\T-\M)^k$ estimate when $f\gg \log$ and $\M=1$. Among these we mention local regularity for the $\dib$-Neumann operator, geometric necessary conditions, lower bounds for the Bergmann metric and the Property $(f\T-\M\T-P)^k$.
 
\begin{theorem}\Label{main3}
Let $\Om$ be a bounded  pseudoconvex domain and let $z_0\in b\Om$. Suppose that $(f\T-1)^1$  holds in a neighborhood $U$ of $z_0$ with $f\gg \log$, that is, $\underset{\xi\to \infty}{\lim}\dfrac{f(|\xi|)}{\log |\xi|}=+\infty$. Then, 
\begin{enumerate}
  \item The $\dib$-Neumann operator is locally regular.
  \item If, for an 1-dimensional complex analytic variety $Z$ passing through $z_0$,  $b\Om$ has type $\le F$ along to $Z$ at $z_0$ in the sense that  
 $$|r(z)|\lesssim F(|z-z_0|)\quad \T{for all $z\in Z$ sufficiently close to $z_0$},$$
then, for small $\delta$ and with $F^*$ denoting the inverse function to $F$,
$$
\frac{f(\delta^{-1})}{\log(\delta^{-1})}\simleq F^*(\delta^{-1})^{-1}.
$$ \index{Notation ! $F^*$}
\item Given a constant $\eta>0$, there exists  $U'\subset U$, so that for  any $z\in U'\cap \Om$ and any $X\in T^{1,0}_z\C^n$, the Bergman metric $B(z, X)$ is bounded from below by $|X|  \dfrac f{\log}(\delta^{-1+\eta}(z))$, where $\delta(z)$ is the distance of $z$ to $b\Om$.
\item Given a constant $\eta$, there exist $U'\subset U$ and  constant $\delta_0>0$, so that Property $(\tilde f\T-1\T-P)^1$ holds in $U'$ with $\tilde f(\delta^{-1})=\dfrac{f}{\log^{\frac{3}{2}+\eta}}(\delta^{-1+\eta})$, $0<\delta\le \delta_0.$ 
 \end{enumerate}
\end{theorem}

As for the global problem, it is well known that the  compactness estimate implies global regularity for the $\dib$-Neumann operator. In this thesis we  give an example in which compactness estimate does not hold but, nonetheless, the $\dib$-Neumann operator is globally regular. 
This refines former work by Krantz \cite{Kr88}. 
Moreover, we show that global regularity  follows from an estimate which is weaker than compactness estimate according to Straube \cite{S08}. \\
    
 The thesis also contains the discussion of  Property $(f\T-\M\T-P)^k$ for some classes of domains such as domains satisfying    $Z(k)$ as well as   decoupled or  regular coordinate domains.\\

The thesis is structured as follows.
  In Chapter 2,  we give the background of the $\bar\partial$-Neumann problem
over $q$-pseudoconvex/concave domains following the guidelines of Folland-Kohn \cite{FK72} and Zampieri \cite{Z08}.
In Chapter 3,  we introduce the method of  weights and prove Theorem~\ref{main1}.
 This part is a development of Catlin \cite{C84} and \cite{C87}.
In Chapter 4,  we study, by microlocal tools, the relation between the $\bar\partial$-Neumann system and the tangential system and prove Theorem~\ref{main2}. The main references of this section are Kohn \cite{K86} and Kohn \cite{K02}.
In Chapter 5, we combine the techniques by Kohn \cite{K02} and Catlin \cite{C83} to show that an $(f\T-\M)^k$ estimate implies a local $H_s$ estimate with a 1-parameter family  of cutoff functions. This $H_s$ estimate is the main tool for the proof of Theorem \ref{main3}. The proof of the second (resp. third) part of  Theorem \ref{main3} originates from \cite{C83} (resp. \cite{McN92})  but contains new ideas which are needed to handle domains which are no more of finite type.
Global regularity is discussed in Chapter 6. 
In Chapter 7, we discuss Property $(f\T-\M\T-P)^k$ for decoupled domains and state optimal subelliptic estimates for  regular coordinate domains (for whose detailed analysis we refer to \cite{KZ08c}).



\chapter{Background} 
\label{Chapter2}

In this chapter, we provide most of the background which is needed to introduce the $(f\T-\M)^k$ estimate for  the $\bar\partial$-Neumann problem on a $q$-pseudoconvex/concave domain $\Omega$ at a point $z_0\in b\Omega$.

\section{Terminology and notations}
\noindent

For $z\in \C^n$, we denote by  $\C T_z\C^n$ the complex-valued tangent vectors to $\C^n$ at $z$. We have the direct sum decomposition $\C T_z\C^n=T^{1,0}_z\C^n\oplus T_z^{0,1}\C^n$, where $T^{1,0}_z\C^n$ and $T_z^{0,1}\C^n$ denote the holomorphic and anti-holomorphic tangent vectors at $z$ respectively.\\

Denote by $\mathcal A_z^{0,k}$ the space of $(0,k)$-forms at $z$ and by $\la~,~\ra_z$ the pairing of $\mathcal A^{0,k}_z$ with its dual space. We also denote by $\la~,~\ra_z$ the inner product induced on $\mathcal A_z^{0,k}$ by the hermitian metric and by $|~|_z$ the associated norm.\\

Let $L_1,...,L_n$ be an orthonormal basis of $T^{1,0}\C^n$ in a neighborhood $U$ of $z_0$; then, for $z\in U$, we have 
$\la (L_i)_z,(L_j)_z\ra_z=\delta_{ij}$ where $\delta_{ij}$ is the Kronecker symbol.  \index{Notation ! $\delta_{ij}$}
\\

 Let $\om_1,\cdots, \om_n$ be the dual basis of (1,0)-forms on $U$, so that for each $z\in U$ we have $\la(\om_i)_z,(L_j)_z\ra_z=\delta_{ij}$. We denote  by $\bar L_1,\cdots, \bar L_n$ the conjugates of $L_1,\cdots,L_n$, respectively; these form an orthonormal basis of $T^{0,1}\C^n$ on $U$. Denote by $\bom_1,\cdots,\bom_n$, the conjugates of $\om_1,\cdots,\om_n$ respectively; they are the local basis  of (0,1)-forms on $U$ which is dual to $\bar L_1,\cdots, \bar L_n$. In this basis, for any $\phi\in C^{\infty}(U)$, we can write
$$d\phi=\sum_{j=1}^{n}L_j(\phi)\omega_j+\sum_{j=1}^{n}\bar{L}_j(\phi)\bar\omega_j.$$
Then,   one defines 
$$\partial\phi=\sum_{j=1}^{n}{L}_j(\phi)\omega_j
\T{~~~and~~~~} 
\bar\partial\phi=\sum_{j=1}^{n}\bar{L}_j(\phi)\bar\omega_j.$$\\

We set $\phi_{ij}$ to be the coefficients of $\di\dib \phi$ in the basis $\{\omega_j\}_j$, i.e.
\begin{eqnarray}\Label{coeff}
\di\dib \phi=\sum_{ij}\phi_{ij}\om_i\we\bom_j.\index{Notation ! $\phi_{ij}$}
\end{eqnarray}
For each $k=1,...,n$, let $\bar{c}_{ij}^k$ be the coefficients of the $2$-forms $\partial\bar\omega_k$, that is, 
$$\di\bom_k=\sum_{ij}\bar{c}_{ij}^k\om_i\we\bom_j.$$
Then $\phi_{ij}$ can be calculated as follows 
\begin{eqnarray*}
\di\dib\phi&=&\di\Big(\sum_{k}\bar{L}_k(\phi)\bom_k\Big)\\
&=&\sum_{i,k}L_i\bar{L}_k(\phi)\om^i\we\bom^k+\sum_k\bar{L}_k(\phi)\sum_{i,j}\bar{c}^k_{ij}\om_i\we\bom_j\\
&=&\sum_{i,j}\Big(L_i\bar{L}_j(\phi)+\sum_k\bar{c}^k_{ij}\bar{L}_k(\phi)\Big)\om_i\we\bom_j.
\end{eqnarray*}
From the fact that $\di\dib+\dib\di=0$, we have
\begin{eqnarray}\Label{phiij}
\phi_{ij}=L_i\bar{L}_j(\phi)+\sum_k\bar{c}^k_{ij}\bar{L}_k(\phi)=\bar{L}_jL_i(\phi)+\sum_k{c}^k_{ji}L_k(\phi).
\end{eqnarray} \index{Notation ! $\phi_{ij}$}
From \eqref{phiij} we obtain
\begin{eqnarray}\Label{[]}
[L_i,\bar L_j]=\sum_k{c}^k_{ji}L_k-\sum_k\bar{c}^k_{ij}\bar{L}_k.
\end{eqnarray}
where $[L_i,\bar L_j]=L_i\bar L_j-\bar L_j L_i$ denotes the commutator.\\

Let $\Om\subset \C^n$ be an open subset of $\C^n$ and let $b\Om$ denote the boundary of $\Om$. Throughout this thesis we  restrict ourselves to a domain $\Om$ such that 
$b\Om$ is smooth in the following sense. We assume that in a neighborhood $U$ of $b\Om$ there exist a $C^\infty$ real-valued function $r$ such that $dr\not=0$ in $U$ and $r(z)=0$ if and only if $z\in b\Om$. Without loss of generality, we may assume that $r>0$ outside of $\bar\Om$ and $r<0$ in $\Om$.\\



For $z_0\in b\Om$,  we fix $r$  so that $|\di r|=1$ in a neighborhood $U$ of $z_0$. We choose  $\om_1,\cdots, \om_n$ to be (1,0)-forms on $U$ such that $\om_n=\di r$ and such that $\la \om_i,\om_j\ra=\delta_{ij}$ for $z\in U$. We then define $L_1,\cdots, L_n, \bar L_1,\cdots,\bar L_n, \bom_1,\cdots,\bom_n$ as above. Note that on $U\cap b\Om$, we have 
$$L_j(r)=\bar L_j(r)=\delta_{jn}.$$
Thus,  $L_1,\cdots, L_{n-1}$ and $\bar L_1,\cdots,\bar L_{n-1}$ are local bases of $T^{1,0} b\Om:=\C T b\Om\cap T^{1,0}\C^n$  and $T^{0,1}b\Om:=\C T b\Om\cap T^{0,1}\C^n$ respectively, where $\C T b\Om$ is the space  of complex-valued tangent vectors to $ b\Om$.
We define a vector field $T$ on $U\cap b\Om$ with values in $\C T b\Om$ by 
 $$T=L_n-\bar L_n.$$  \index{Notation ! $T$}

Observe that $L_1,\cdots, L_{n-1}, \bar L_1,\cdots,\bar L_{n-1}, T$ are local basis of $\C T b\Om$ over $U\cap b\Om$. Using integration by parts, we get the proof of following lemma 
\begin{lemma}\Label{l2.1} Let $\varphi, \psi\in C^\infty_c(U\cap \bar\Om)$, then we have
$$(\bar L_j\varphi, \psi)=-(\varphi,  L_j\psi)+\int_{b\Om}\bar L_j(r)\varphi\bar \psi dS+(\varphi,a_j\psi)$$
where $a_j\in C^\infty(U\cap\bar\Om)$ and $dS$ denotes the element of area on $b\Om$. 
\end{lemma}

Let $(r_{ij})$ be the Levi matrix of $r$. 
Substituting $r$ for $\phi$ in \eqref{phiij}, we get $r_{ij}=\bar c^n_{ij}= c^n_{ji}$ and by \eqref{[]} we conclude
\begin{eqnarray}\Label{rij}
[L_i,\bar L_j]= r_{ij}T+\sum_{k=1}^{n-1}c^k_{ji}L_k-\sum_{k=1}^{n-1}\bar c^k_{ij}\bar L_k.
\end{eqnarray}


 Remember that we have  denoted  by $\A^{0,k}$ the space of $(0,k)$-forms in $C^\infty(\bar\Om)$  and  by $C^\infty_c( U\cap\bar\Om)^{k}$
 those which have compact support in $U$. 
  If $u\in C^\infty_c(U\cap \bar\Om)^{k}$, then $u$ can be written as 
\begin{equation}\Label{2u}
u=\sumJ u_J\bom_J,
\end{equation}
where   $\sum'$ denotes summation over strictly increasing indices $J=j_1<...<j_k$ and where $\bar\omega_J$ denotes the wedge product $\bom_J=\bom_{j_1}\we...\we\bom_{j_k}$.  When the multiindices are not ordered, the coefficients are assumed to be alternant. Thus, if $J$ decomposes as $J=jK$, then $u_{jK}=\epsilon^{jK}_{J} u_J$ where $\epsilon^{jK}_{J}$ is the sign if the permutation $jK\simto J$.  \index{Notation ! $\sum'$}

There is a well defined Cauchy-Riemann complex
\begin{equation*}
\A^{0,k-1}\overset{\bar\partial}\to \A^{0,k}\overset{\bar\partial}\to \A^{0,k+1}.
\end{equation*}
The action of $\bar\partial$ on a $(0,k)$ form is \index{Operator ! Cauchy-Riemann operator, $\dib$}
 \begin{equation}\Label{2dib}
\dib u=\sumJ\sum_{j=1}^n \bar L_j u_J\bom_j\we\bom_J+...
 \end{equation}
where the dots refer to terms of order zero in $u$. \index{Notation ! $\dib$}
\\

We  extend this complex to $L_2^{0,k}(\Om)$ the  space of $(0,k)$-forms with $L_2$-coefficients, so that the Hilbert space techniques may be applied. For each $z\in \bar\Om$, denote by $(dV)_z$ the unique positive $(n,n)$-form such that   $|(dV)_z|=1$. We call $dV$ the volume element. For $ u,v \in L_2^{0,k}(\Om)$, we define the inner product and the norm by
$$(u,v)=\int_\Om \la u,v\ra_z (dV)_z \quad\T{ and }  \quad  \no{u}^2=(u,u).$$ \index{Space ! $L_2^{h,k}(\Om)$}
 For each form of degree $(0,k)$, we define
$$\T{Dom}(\dib)=\{ v\in L_2^{0,k}(\Om) : \dib v \T{ (as distribution)}\in L_2^{0,k+1}(\Om)\}.$$
Then the operator $\dib : \T{Dom}(\dib)\to L_2^{0,k+1}(\Om) $ is well-defined and, by noticing that $\A^{0,k}\subset \T{Dom}(\dib)$, we have $\dib :L_2^{0,k}(\Om)\to L^{0,k+1}_2(\Om) $  as a densely defined operator. Thus, the operator $\dib$ has an $L_2$-adjoint $\dib^*$, defined as follows. If $u\in \T{Dom}(\dib^*)$, there must exist $\alpha\in L_2^{0,k-1}$ such that
$$(v,\alpha)=(\dib v, u)   \T{  for all   }  v\in\T{Dom}(\dib).$$ 
In fact, this $\alpha$ is nothing but the form $\bar\partial^*u$ that we are assuming to exist. But we have
\begin{eqnarray}\begin{split}\Label{dom}
(\dib v,u)=&\sumK\sum_{j=1}(\bar L_jv_K,u_{jK})+(v,\dots)\\
=&\sumK\sum_{j=1}\Big(-(v_K,L_ju_{jK})+\delta_{jn} \int_{b\Om}v_K\bar u_{jK}dS \Big)+(v,\dots)\\
=&\sum_{j=1}\left(v,-\sumK L_ju_{jK}\bom_K)+\delta_{jn}\sumK \int_{b\Om}v_K\bar u_{jK}dS\right)+(v,\dots)
\end{split}
\end{eqnarray}
where dots denote an error term in which $u$ is not differentiated. Here the second equality in \eqref{dom} follows from Lemma \ref{l2.1}. If we pretend that this coincides with the action of $\alpha\in L_2^{0,k-1}$, then the boundary integral must vanish. We have thus obtained the proof of the following
\begin{lemma}
\begin{equation}
\Label{2.10*}
u\in \T{Dom}(\dib^*) \T{ if and only if } u_{jK}|_{b\Om}=0 \T{~~for any~~ } K.
\end{equation}\index{Notation ! $\T{Dom}(\dib^*)$}
\end{lemma}
Over a form satisfying \eqref{2.10*}, the action of the Hilbert adjoint of $\dib$, coindides with of its ``formal adjoint" and is therefore expressed by a ``divergence operator"
\begin{eqnarray}\Label{2dib*}
\dib^*u=-{\sumK}\sum_j L_ju_{jK}\bom_K+... 
\end{eqnarray}
\index{Notation ! $\dib^*$}

\section{The basic estimate}
\indent 
\index{Estimate ! basic }
For a real function $\phi$ in class $C^2$, let the weighted $L^\phi_2$-norm be defined by
$$\no{u}_\phi^2=(u,u)_\phi:=\no{ue^{-\frac{\phi}{2}}}^2=\int_\Om\la u,u\ra e^{-\phi}dV.$$
Let $\dib^*_\phi$ be the $L^\phi_2$-adjoint  of $\dib$. It is easy to see that $\T{Dom}(\dib^*)=\T{Dom}(\dib^*_\phi)$ and 
\begin{equation}\Label{2.5}
\begin{split}
\dib_\phi^*u=&- \sumK\sum_{j=1}^{n}\delta_j^\phi u_K\bom_K+\cdots\\
\end{split}
\end{equation}\index{Notation ! $\dib^*_\phi$}
where $\delta^\phi_j \varphi=e^\phi L_j(e^{-\phi}\varphi)$ and where dots denote an error term in which $u$ is not differentiated and $\phi$ does not occur.\\

By developing the equalities \eqref{2dib} and \eqref{2.5},  the key technical result is contained in the following 
\begin{proposition}\Label{kmh}Let $z_0\in b\Om$ and fix an index $q_o$ with $0\le q_o\le n-1$,  then there exists a neighborhood $U$ of $z_0$ and a suitable constant $C$ such that
\begin{eqnarray} \Label{KMH}
\begin{split}
\no{\bar{\partial} u}^2_{\phi}&+\no{\bar{\partial}^*_\phi u}^2_{\phi}+\no{u}^2_{\phi}\\
\simgeq & {\sumK}\sum_{i,j=1}^{n}(\phi_{ij}u_{iK},u_{jK})_\phi-{\sumJ}\sum_{j=1}^{q_o}(\phi_{jj}u_J,u_J)_\phi\\
&+{\sumK}\sum^{n-1}_{i,j=1}\int_{b\Om}e^{-\phi}r_{ij}u_{iK}\bar{u}_{jK}dS-{\sum_{|J|=q}}'\sum^{q_o}_{j=1}\int_{b\Om}e^{-\phi}r_{jj}|u_{J}|^2dS\\
&+\frac{1}{2}\big(\sum^{q_o}_{j=1}\no{\delta_j^{\phi} u}^2_{\phi}+\sum^{n}_{j=q_o+1}\no{\bar{L}_ju}^2_{\phi} \big)
\end{split}
\end{eqnarray}
for any $u\in C^\infty_c(U\cap \bar\Om)^k\cap \T{Dom}(\dib^*)$.\end{proposition}
{\it Proof. }
Let $Au$ denote the sum in \eqref{2dib}; we have
\begin{eqnarray}\Label{2.13}
\no{Au}_\phi^2=\sumJ \sum_{j=1}^{n}\no{\bar{L}_ju_J}^2_\phi-\sumK\sum_{ij}(\bar{L}_iu_{jK},\bar{L}_ju_{iK})_\phi.
\end{eqnarray}
 Let $Bu$ denote the sum  in \eqref{2.5}; we have
\begin{eqnarray}\Label{2.14}
\no{Bu}_\phi^2=\sumK\sum_{ij}(\delta^\phi_iu_{iK},\delta^\phi_ju_{jK})_\phi.
\end{eqnarray}
Remember that $Au$ and $Bu$ differ from $\bar\partial u$ and $\bar\partial^*_\phi u$ by terms of order $0$ which do not depend on $\phi$. We then have
\begin{eqnarray}\Label{2.15}
\begin{split}
\no{\bar{\partial} u}^2_{\phi}&+\no{\bar{\partial}^*_\phi u}^2_{\phi}\\
&= \no{Au}_\phi^2+\no{Bu}^2_\phi+R\\ 
&=\sumJ\sum_{j=1}^n\no{\bar{L}_ju_J}^2_\phi+\sumK\sum_{i,j=1}^n(\delta_i^\phi u_{iK},\delta_j^\phi u_{jK})_\phi-(\bar{L}_j u_{iK},\bar{L}_i u_{jK})_\phi +R,
\end{split}
\end{eqnarray}
where $R$ is an error coming from the scalar product of $0$-order terms with terms $\bar L_ju_J$,  $L_ju_{jK}$ or $u$.\\

We  want to apply now integration by parts to the term $(\delta_i^\phi u_{iK},\delta_j^\phi u_{jK})_\phi$ and $(\bar{L}_j u_{iK},\bar{L}_i u_{jK})_\phi$. Notice that for each $\varphi,\psi \in C_c^1(U\cap \bar{\Omega})$, similarly as in Lemma \ref{l2.1}, we have 
$$\begin{cases}
(\varphi,\delta^\phi_j \psi)_\phi&=-(\bar{L}_j\varphi,\psi)_\phi+(a_j\varphi,\psi)_\phi+\delta_{jn}\int_{b\Omega}e^{-\phi}\varphi\bar{\psi}dS\\
-(\varphi,\bar{L}_i\psi)_\phi&=(\delta^\phi_i\varphi,\psi)_\phi-(b_i\varphi,\psi)_\phi-\delta_{in}\int_{b\Omega}e^{-\phi}\varphi\bar{\psi}dS
\end{cases}$$
and for some $a_j,b_i\in C^1(\bar{\Omega}\cap U)$ independent of $\phi$.

This immediately implies 
\begin{eqnarray}\Label{2.16}
\begin{cases}
(\delta_i^\phi u_{iK},\delta_j^\phi u_{jK})_\phi&=-(\bar{L}_j\delta^\phi_i u_{iK}, u_{jK})_\phi+\delta_{jn}\int_{b\Omega}e^{-\phi}\delta^\phi_i(u_{iK})\bar{u}_{jK}dS+R\\
-(\bar{L}_j u_{iK},\bar{L}_i u_{jK})_\phi&=(\delta^\phi_i\bar{L}_ju_{iK},u_{jK})_\phi-\delta_{in}\int_{b\Omega}e^{-\phi}L_j(u_{iK})\bar{u}_{jK} dS+R.
\end{cases}
\end{eqnarray}
From here on, we denote by $R$ terms involving the product of $u$ by $\delta^\phi_ju$ for $j\le n-1$ or $\bar{L}_ju$ for $j\le n$ but not twice a derivative of $u$. \\

Recall that $ u_{nK}|_{b\Omega}\equiv0$ and $L_j(u_{nK})|_{b\Omega}\equiv0$  if $j\le n-1$.  We thus conclude that the boundary integrals vanish in both equalities of \eqref{2.16}. Now, by taking the sum of two terms in the right side of \eqref{2.16}, after discarding the boundary integrals, we put in evidence the commutator $[\delta^\phi_i, \bar{L}_j]$
\begin{eqnarray}\Label{2.16b}
(\delta_i^\phi u_{iK},\delta_j^\phi u_{jK})_\phi-(\bar{L}_j u_{iK},\bar{L}_i u_{jK})_\phi=([\delta^\phi_i, \bar{L}_j]u_{iK},u_{jK})_\phi+R.
\end{eqnarray}
 \\

Notice that \eqref{2.16} is also true if we replace both $u_{iK}$ and $u_{jK}$ by $u_J$ for indices $i=j\le q_0$. Then we obtain
\begin{eqnarray}\Label{2.17}
\no{\bar{L}_j u_{J}}^2_\phi=\no{\delta^\phi_ju_J}^2_\phi-([\delta^\phi_j, \bar{L}_j]u_J,u_J)_\phi+R.
\end{eqnarray}

Applying \eqref{2.16b} and \eqref{2.17} to the last line in \eqref{2.15}, we have
\begin{equation}\Label{2.19}
\begin{split}
\no{\bar{\partial} u}^2_{\phi}&+\no{\bar{\partial}^*_\phi u}^2_{\phi}\\
&={\sumK}\sum_{i,j=1}^{n}([\delta^\phi_i, \bar{L}_j]u_{iK},u_{jK})_\phi-{\sumJ}\sum_{j=1}^{q_0}([\delta^\phi_j, \bar{L}_j]u_J,u_J)_\phi\\
&+{\sumJ} \Big(\sum_{j=1}^{q_0}\no{\delta^\phi_ju_J}^2_\phi+\sum_{j=q_0+1}^{n}\no{\bar{L}_j u_{J}}^2_\phi\Big)+R.
\end{split}
\end{equation} 
Now we calculate the commutator $[\delta^\phi_i, \bar{L}_j]$,
\begin{equation}\Label{2.20}
\begin{split}
[\delta^\phi_i, \bar{L}_j]&=L_i\bar L_j\phi+[L_i,\bar L_j] \\
&=\phi_{ij}+\sum_k^{n}{c}^k_{ji}\delta^\phi_k-\sum_{j=1}^n\bar{c}^k_{ij}\bar{L}_k\\
&=\phi_{ij}+r_{ij}(\delta^\phi_n-\bar L_n)+\sum_k^{n-1}{c}^k_{ji}\delta^\phi_k-\sum_{j=1}^{n-1}\bar{c}^k_{ij}\bar{L}_k
\end{split}
\end{equation} 
where we have used \eqref{phiij} and  \eqref{rij}.\\

Since $\langle L_n,\partial r\rangle=1$, we have
\begin{equation}\Label{2.21}
(r_{ij}\delta^\phi_nu_{iK}, u_{jK})_\phi=\int_{b\Omega}r_{ij}e^{-\phi}u_{iK}u_{jK}dS+R.
\end{equation} 
Substituting \eqref{2.20} in \eqref{2.19} and combining with \eqref{2.21}, we get
\begin{equation}\Label{2.22}
\begin{split}
\no{\bar{\partial} u}^2_{\phi}&+\no{\bar{\partial}^*_\phi u}^2_{\phi}\\
&={\sumK}\sum_{i,j=1}^{n}(\phi_{ij}u_{iK},u_{jK})_\phi-{\sumJ}\sum_{j=1}^{q_0}(\phi_{jj}u_J,u_J)_\phi\\
&+{\sumK}\sum_{i,j=1}^{n-1}\int_{b\Omega}r_{ij}u_{iK}\bar{u}_{jK}e^{-\phi}dS-{\sumJ}\sum_{j=1}^{q_0}\int_{b\Omega}r_{jj}|u_J|e^{-\phi}dS\\
&+{\sumJ} \Big(\sum_{j=1}^{q_0}\no{\delta^\phi_ju_J}^2_\phi+\sum_{j=q_0+1}^{n}\no{\bar{L}_j u_{J}}^2_\phi\Big)+R.
\end{split}
\end{equation} 

We denote by $S$ the sum in the last line in \eqref{2.22}. To conclude our proof, we only need to prove that 
for a suitable $C$ independent of $\phi$ we have 
\begin{eqnarray}\Label{2.23}
R\le \frac{1}{2}{\sumJ} \Big(\sum_{j=1}^{q_0}\no{\delta^\phi_ju_J}^2_\phi+\sum_{j=q_0+1}^{n}\no{\bar{L}_j u_{J}}^2_\phi\Big)+C\no{u}_\phi^2.
\end{eqnarray}
In fact, if we point our attention at those terms which involve $\delta^\phi_ju$ for $j\le q_o$ or $\bar{L}_ju$ for $q_o+1\le j\le n$, then \eqref{2.23} is clear since $S$ carries the corresponding square $\no{\delta^\phi_ju}^2_\phi$ and $\no{\bar{L}_ju}^2_\phi$. Otherwise, we note that for $j\le n-1$ we may interchange $\bar{L}_j$ and $\delta^\phi_j$ by means of integration by parts: boundary integrals do not occur because $L_j(r)=0$ on $b\Omega$ for $j\le n-1$. As for $\delta^\phi_n$, notice that it only hits coefficients whose index contains $n$ and hence $u_{nK}=0$ on $b\Omega$. So $\delta^\phi_n(u_{nK})\bar{u}$ is also interchangeable with $u_{nK}\bar{L}_n\bar{u}$ by integration by parts. This concludes the proof of Proposition \ref{kmh}.

$\hfill\Box$\\

For the choice $\phi=0$, we can rewrite the estimate \eqref{KMH} as 
\begin{equation}\Label{unweight}\begin{split}
  Q(u,u) \simge &{\sumK}\sum^{n-1}_{i,j=1}\int_{b\Om}u_{iK}\bar{u}_{jK}dS-{\sum_{|J|=q}}'\sum^{q_o}_{j=1}\int_{b\Om}|u_{J}|^2dS\\
&+\sum^{q_o}_{j=1}\no{L_j  u}^2+\sum^{n}_{j=q_o+1}\no{\bar{L}_ju}^2
\end{split}
\end{equation}
for any $u\in C^\infty_c(U\cap \bar\Om)^k\cap\T{Dom}(\dib^*)$.\\

Observe that if $\varphi \in C^\infty_c(U\cap \bar\Om)$ with  $\varphi|_{b\Om}\equiv 0$ on $U\cap b\Om$, then each $\no{L_j \varphi}^2$ can be interchanged with $\no{\bar L_j \varphi}^2+R$ even for $j=n$ due to the vanishing of the boundary integral.
Thus, 
\begin{equation}\Label{index1}
\begin{split}
\sum^{q_o}_{j=1}\no{L_j \varphi}^2+&\sum^{n}_{j=q_o+1}\no{\bar{L}_j \varphi}^2+C\no{\varphi}^2\\
\ge & \frac{1}{2}\sum^{n}_{j=1}\big(\no{L_j\varphi}^2+\no{\bar{L}_j\varphi}^2\big) +\no{\varphi}^2\\
\simge &\no{\varphi}_1^2,
\end{split}
\end{equation}
where $\no{.}_1$ is the Sobolev norm of index 1.\\

In conclusion, we get an estimate which fully expresses  the interior elliptic regularity of the system $(\dib, \dib^*)$  
\begin{eqnarray}\Label{elliptic}
Q(u,u)\simge \no{u}^2_1 
\end{eqnarray}
for any  $u\in C_c^\infty(U\cap \bar\Om)^k$ with $u|_{U\cap b\Om}=0$. 
We keep our choice of $\phi=0$.
 Using observation \eqref{index1} for $\varphi$ replaced by $u_{nK}$ for any $K$, we get
\begin{equation}\begin{split}\Label{2.26}
  Q(u,u) \simge &{\sumK}\sum^{n-1}_{i,j=1}\int_{b\Om}r_{ij}u_{iK}\bar{u}_{jK}dS-{\sum_{|J|=q}}'\sum^{q_o}_{j=1}\int_{b\Om}r_{jj}|u_{J}|^2dS\\
&+\sum^{q_o}_{j=1}\no{L_j u}^2+\sum^{n}_{j=q_o+1}\no{\bar{L}_ju}^2+\sumK\no{u_{nK}}_1^2
\end{split}
\end{equation}
for any $u\in C^\infty_c(U\cap \bar\Om)^k\cap \T{Dom}(\dib^*)$.\\

Conversely, and still under the choice $\phi=0$, we have
\begin{equation}\Label{2.27}
\begin{split}
  Q(u,u) \lesssim &\big|{\sumK}\sum^{n-1}_{i,j=1}\int_{b\Om}r_{ij}u_{iK}\bar{u}_{jK}dS-{\sum_{|J|=q}}'\sum^{q_o}_{j=1}\int_{b\Om}r_{jj}|u_{J}|^2dS\big|\\
&+\sum^{q_o}_{j=1}\no{L_j u}^2+\sum^{n}_{j=q_o+1}\no{\bar{L}_ju}^2+\no{u}^2
\end{split}
\end{equation}
for any $u\in C^\infty_c(U\cap \bar\Om)^k\cap \T{Dom}(\dib^*)$. This inequality is a consequence of the calculation in  Proposition~\ref{kmh} and holds without the assumption of pseudoconvexity.\\

\section{$q$-pseudoconvex/concave domains}
\indent
 
 Let $\lambda_1(z)\le ...\le \lambda_{n-1}(z)$ be the eigenvalues of $(r_{jk}(z))^{n-1}_{j,k=1}$ and denote by $s^+_{b\Om}(z)$, $s^-_{b\Om}(z)$ ,  $s^0_{b\Om}(z) $ their number according to the different sign. \\

 We take a pair of indices $1\le q\le n-1$ and $0\le q_o\le n-1$ such that $q\not=q_o$. We assume that there is a bundle $\mathcal V^{q_o}\in T^{1,0}b\Om$ of rank $q_o$ with smooth coefficients that, by reordering, we may suppose to be the bundle $\mathcal V^{q_o}=\T{ span }\{ L_1,...,L_{q_o}\}$, such that 
\begin{equation}\Label{2qpseu}
\sum_{j=1}^{q}\lambda_j(z)-\sum_{j=1}^{q_o}r_{jj}(z)\ge 0 \qquad z \in  U\cap b\Om.
\end{equation} 
We  define $\Om$ to be $q$-pseudoconvex or $q$-pseudoconcave according to $q>q_o$ or $q<q_o$. 
\begin{lemma}
Condition \eqref{2qpseu} is equivalent to
\begin{equation}
\sumK\sum_{ij=1}^{n-1}r_{ij}u_{iK}\bar u_{jK}-\sum_{j=1}^{q_o}r_{jj}|u|^2\ge 0 \T{ on } U\cap b\Om
\end{equation}
for any $u\in C_c^\infty(U\cap\bar\Om)^q\cap \T{Dom}(\dib^*)$.
\end{lemma}
{\it Proof. } The proof of the Lemma immediately follows from the estimate
$$
\sumK\sum_{ij=1}^{n-1}r_{ij}u_{iK}\bar u_{jK}\geq\underset{j=1}{\overset q\sum}\lambda_j|u|^2,
$$
for any $u\in C_c^\infty(U\cap\bar\Om)^q\cap \T{Dom}(\dib^*)$ with equality for a suitable $u$. In turn, the proof of this estimate is obtained by diagonalizing the matrix $(r_{ij})$ (cf. \cite{H65} and \cite{C87}).

$\hfill\Box$\\

As it has already been noticed, \eqref{2qpseu} for $q>q_o$  implies $\lambda_{q}\ge 0$; hence \eqref{2qpseu} is still true if we replace the first sum $\sum^{q}_{j=1}\cdot$ by $\sum^{k}_{j=1}\cdot$ for any $k$ such that $q\le k\le n-1$. Similarly, if it holds for $q<q_o$, then $\lambda_{q+1}\le0$ and hence it also holds with $q$ replaced by $k\le q$ in the first sum.\\

We notice  that $q$-pseudoconvexity/concavity is invariant under a change of an orthonormal basis but not of an adapted frame. In fact, not only the number, but also the size of the eigenvalues comes into play. Thus, when we say that $b\Om$ is $q$-pseudoconvex/concave, we mean that there is an adapted frame in which \eqref{2qpseu} is fulfilled. Sometimes, it is more convenient to put our calculatiuons in an orthonormal frame. In this case, it is meant that the metric has been changed so that the adapted frame has become orthonormal.

\begin{example}{\rm Let $s^-(z)$ be constant for $z\in b\Om$ close to $z_0$; then \eqref{2qpseu} holds for $q_o=s^-$ and $q=s^-+1$. In fact, we have $\lambda_{s^-}<0\le \lambda_{s^-+1}$, and therefore the negative eigenvectors span a bunlde $\mathcal{V}^{q_o}$ for $q_o=s^-$ that, identified with the span of the first $q_o$ coordinate vector fields, yields $\sum^{q_o+1}_{j=1}\lambda_j(z)\ge \sum^{q_o}_{j=1} r_{jj}(z)$. Note that a pseudoconvex domain is characterized by $s^-(z)\equiv0$, thus, it is 1-pseudoconvex in our terminology. 
In the same way, if $s^+(z)$ is constant for $z\in b\Om$ close to  $z_0$, then $\lambda_{s^-+s^0}\le0<\lambda_{s^-+s^0+1}$. Then, the eigenspace of the eigenvectors $\le0$ is a bundle which, identified to that of the first $q_o=s^-+s^0$ coordinate vector fields yields \eqref{2qpseu} for $q=q_o-1$. In particular a pseudoconcave domain, that is, a domain which satisfies $s^+\equiv0$, is $(n-2)$-pseudoconcave in our terminology.}
\end{example}

\begin{example}
\Label{e2.2}
{\rm Let $\Om$ satisfy $Z(q)$ condition at $z_0$, that is, $s^+(z)\ge n-q$ or $s^{-}(z)\ge q+1$ for $z\in U\cap b\Om$. Thus $\Om$ is strongly $q$-pseudoconvex or strongly $q$-pseudoconcave at $z_0$}
\end{example}
\begin{example}{\rm  Let $\Om$ be a domain in a neighborhood $z_0$ with a defining function 
$$r=\T{2Re}z_n-a(z_1,...,z_{q_o})+b(z_{q_o+1},...,z_{n-1})$$
where $a$ and $b$  are real functions which vanish at $0$ and such that $(a_{z_i\bar z_j})^{q_o}_{ij=1}$ and $(b_{z_i\bar z_j})^{n-1}_{ij=q_o+1}$ are semipositive matrices. We can check that $\Om$ is $(q_o+1)$-pseudoconvex and $(q_o-1)$-pseudoconcave at $z_0$ (cf. Proposition~\ref{p5.3} which follows).}\\
\end{example}

\begin{remark}\Label{convex-concave}{\rm
By the identity  $\sum_{j=1}^{n-1}\lambda_j(z)=\sum^{n-1}_{j=1}r_{jj}(z)$, we get
$$\sum_{j=1}^{q}\lambda_j-\sum_{j=1}^{q_o}r_{jj}=\sum_{j=q+1}^{n-1}(-\lambda_j)-\sum_{j=n-1}^{q_o+1}(-r_{jj}).$$
Therefore, if $\Om$ defined by $r<0$ is $q$-pseudoconvex ($q$-pseudoconcave) at $z_0$, then $\C^n\setminus \bar\Om=\{-r<0\}$ is $(n-q-1)$-pseudoconcave ($(n-q-1)$-pseudoconvex)  at $z_0$. 
}\end{remark}

\section{Tangential pseudodifferential operators}
\noindent

In our study of the $(f\T-\M)^k$-estimates, we  use tangential pseudo-differential operators in a neighborhood of $z_0\in b\Om$. These are expressed in terms of boundary coordinates which are defined as follows.

\begin{definition} A system of real $C^\infty$ coordinates, defined in a neighborhood $U$ of $z_0\in b\Om$ is said a system of {\it boundary coordinates} if one of the coordinates is the defining function $r$. We  denote such a system by $(x,r)=(x_1,...,x_{2n-1},r)\in\R^{2n-1}\times\R$ and call the $x_j$'s the {\it tangential} coordinates and $r$ the {\it normal coordinate}.
\end{definition}

We denote the dual variables of $x$ by $\xi$, and define $x\cdot\xi=\sum x_i\xi_j $; $|\xi|^2=\sum\xi_j^2$. For $\varphi \in C_c^\infty(U\cap\bar\Om)$ we define $\tilde{\varphi }$, the {\it tangential Fourier transform} of $\varphi$, by
$$\tilde{\varphi}(\xi, r)=\int_{\mathbb{R}^{2n-1}}e^{-ix\cdot\xi}\varphi(t,r)dt.$$Denote by $\Lambda_\xi=(1+|\xi|^2)^{\frac{1}{2}}$ the standard ``tangential" elliptic symbol of order 1 and by $\Lambda$ the operator with symbol $\Lambda_\xi$.
For $f\in C^\infty([1,+\infty))$ we define $f(\Lambda)\varphi$ by \index{Notation ! $f(\La)$ }
\begin{equation}
f(\Lambda)\varphi(\xi,r)=(2\pi)^{-(2n-1)}\int e^{ix\cdot\xi}f(\Lambda_\xi)\tilde{\varphi}(\xi,r)d\xi.
\end{equation}
 Hence
\begin{equation}
\no{f(\Lambda)\varphi}^2=\int_{-\infty}^0\int_{\R^{2n-1}}f(\Lambda_\xi)^2|\tilde{\varphi}(\xi,r)|^2 d\xi dr.
\end{equation}
In the case $f(t)=t^s$, $ s\in \R$, we define tangential Sobolev norms  by
\begin{equation}
|||\varphi|||_s=\no{\Lambda^s\varphi}
\end{equation}
\begin{lemma}\Label{l2.7} Let $f,g\in C^\infty([1,+\infty))$ satisfy $f\gg g$.  Then,   for any $\epsilon>0$ and $s\in\R^+ $ there exists a constant $C_{\epsilon,s}$ such that 
$$\no{g(\La)\varphi}^2\le \epsilon \no{f(\La)\varphi}^2+C_{\epsilon,s} |||\varphi|||^2_{-s}$$ 
for any $\varphi\in C_c^\infty(U\cap \bar\Om)$.
\end{lemma}
{\it Proof. } Since $f\gg g$, that is, $\underset{t\to+\infty}{\lim}\frac{f(t)}{g(t)}=+\infty$, then for any $\epsilon>0$,  there exists a constant $t_\epsilon>0$ such that $g(\La_\xi)\le \epsilon f(\La_\xi)$ for  $|\xi| \ge t_\epsilon$. Hence 
\begin{equation*}
\begin{split}
\no{g(\Lambda)\varphi}^2=&\int_{-\infty}^0\int_{|\xi|\ge t_\epsilon}g(\Lambda_\xi)^2|\tilde{\varphi}(\xi,r)|^2dr d\xi+\int_{-\infty}^0\int_{|\xi|\le t_\epsilon }g(\Lambda_\xi)^2|\tilde{\varphi}(\xi,r)|^2dr d\xi\\
\le &\epsilon \int_{-\infty}^0\int_{|\xi|\ge t_\epsilon}f(\Lambda_\xi)^2|\tilde{\varphi}(\xi,r)|^2dr d\xi+\int_{-\infty}^0\int_{|\xi|\le t_\epsilon }g(\Lambda_\xi)^2|\tilde{\varphi}(\xi,r)|^2dr d\xi\\
\le&\epsilon \no{f(\La) \varphi}^2+C_{\epsilon,s} |||u|||^2_{-s}.\\
\end{split}
\end{equation*}
 
 $\hfill\Box$

In $\R^{2n}_{-}=\{(x_1,...,x_{2n})| x_{2n}=r\le 0 \}$, the Schwartz space $\S$ is the space of smooth  functions which decrease rapidly at infinity \index{Space !  $\S$}
$$\S=\{f\in C^\infty(\R^{2n}_{-}) | \, \sup_{x\in \R^{2n}_{-}}|x^\alpha D^\beta f(x)|<\infty \,\T{ for any }  \alpha,\beta\},$$
where $\alpha$ and $\beta$ are $2n$-indices and $D^\beta=\frac{\di^{\beta_1}}{\di x^{\beta_{1}}}\cdots \frac{\di^{\beta_{2n}}}{\di x^{\beta_{2n}}}.$ Recall the class of functions $\mathcal F$  defined by
$$\mathcal F=\{f\in C^\infty([1,+\infty) \Big| f(t)\lesssim t^{\frac{1}{2}}; f'(t)\ge 0 \T{ and }|f^{(m)}(t)| \lesssim |\frac{f^{(m-1)}(t)}{t}|  \}$$
for any $m\in \Z^+$.

\begin{proposition}\Label{p2.8}
Let  $f\in \mathcal F$; $a\in \S(\R_{-}^{2n})$; $s\in \R$; and $S$ be the a vector field with coefficients in $\S(\R_{-}^{2n})$.   Then,   we have
\begin{enumerate}
  \item[(i)]  $|||[f(\La),a]\varphi|||_s\lesssim |||f(\La)\varphi|||_{s-1}$;
  \item[(ii)] $|||[f(\La),S]\varphi|||_s\lesssim |||f(\La)\varphi|||_{s}+|||D_r f(\La) \varphi|||_{s-1}$;\end{enumerate}
for any $\varphi\in C^\infty_c(U\cap \bar\Om)$.
\end{proposition}
{\it Proof. }
Let $\sigma(A)$ denote the symbol of operator $A$. We have
\begin{equation}\Label{2.33}
\sigma([A,B])=\sum_{k>0}\frac{(\di/\di \xi)^k\sigma(A) D_x^k\sigma(B)-(\di/\di \xi)^k\sigma(B) D_x^k\sigma(A)}{k!}.
\end{equation}
(Notice here that the $k=0$ term cancels out.) Using formula \eqref{2.33} for $A=f(\La)$ where $f\in \F$ and $B=a\in \S(\R_{-}^{2n})$, we obtain
$$|\sigma([f(\La), a])|\lesssim \sigma(\La^{-1} f(\La)).$$
The second part of this proposition immediately follows from the first part if we  notice that $S=\sum_k a_k(x,r)\frac{\di}{\di x_k}+b(x,r)\frac{\di}{\di r}$.\\

$\hfill\Box$

\begin{proposition}\Label{1ES} For $f\in \F$, we have the estimate
\begin{eqnarray}
\begin{split}
\no{\La^{-1}f(\Lambda) \varphi}_1^2 \lesssim &\sum_{j=1}^{q_0}\no{L_j\Lambda^{-1}f(\Lambda) \varphi}^2+\sum_{j=q_0+1}^n\no{ \bar L_j \Lambda^{-1}f(\Lambda)\varphi}^2\\
&+\no{\Lambda^{-1}f(\Lambda)\varphi}^2+\no{\Lambda^{-1/2}f(\Lambda) \varphi_b}_b^2
\end{split}
\end{eqnarray}
for any $\varphi\in C_c^\infty(U\cap \bar{\Om})$.
\end{proposition}

The above proposition is a variant of Theorem (2.4.5) of \cite{FK72}. The key in proving this theorem is a passage from functions in $C^\infty(U\cap \bar\Om)$ to $C^\infty(U\cap b\Om)$. This is done by using an extension of $\psi$ from $U\cap b\Om$ to $U\cap \bar\Om$. For $\psi\in C_c^\infty(U\cap b\Om)$, we define $\psi^{(e)}\in C^\infty (\{(x,r)\in \R^{2n}\big| r\le 0\})$ by
$$\psi^{(e)}(x,r)=(2\pi)^{-(2n-1)}\int e^{ix\cdot \xi} e^{r(1+|\xi|^2)^{1/2}}\tilde \psi(\xi,r)d\xi$$
so that $\psi^{(e)}(x,0)=\psi(x)$. \index{Notation ! $\psi^{(e)}$ }
\begin{lemma}\Label{l2.10}
For each $k\in \Z$ with $k\ge 0$, $s\in \R$ and $f\in \F$, we have 
\begin{enumerate}
  \item[(i)] $|||r^kf(\Lambda)\psi^{(e)}|||_s\cong C_{k}\no{f(\Lambda)\psi}_{b, s-k-\frac{1}{2}}$
  \item[(ii)] $|||D_rf(\Lambda)\psi^{(e) }|||_s\cong \no{f(\Lambda)\psi}_{b,s+\frac{1}{2}}$
  \end{enumerate}
for any $\psi\in C_c^\infty(U\cap b\Om)$.
\end{lemma}
{\it Proof.}  {\it (i).} We observe that for every positive integer $k$, we have by integration by parts 
\begin{eqnarray}
\int^{0}_{-\infty}r^{2k}e^{r(1+|\xi|^2)^{1/2}}dr=C'_k (1+|\xi|^2)^{-k-\frac{1}{2}}.
\end{eqnarray}
Hence, 
\begin{eqnarray}
\begin{split}
|||r^k f(\Lambda) \psi^{(e)}|||_s^2\cong&\int_{-\infty}^0\int_{\R^{2n-1}}r^{2k}e^{2r(1+|\xi|^2)^{1/2}}(1+|\xi|^2)^sf((1+|\xi|^2)^{1/2})^2|\psi(\xi)|^2d\xi dr\\
\cong& C'_k\int_{\R^{2n-1}}(1+|\xi|^2)^{s-k-\frac{1}{2}}f((1+|\xi|^2)^{1/2})^2|\psi(\xi)|^2d\xi\\
\cong&C'_k|||f(\Lambda)\psi|||_{b,s-k-\frac{1}{2}}^2.
\end{split}
\end{eqnarray}\\

{\it (ii). } Since the derivative $D_r$ does not affect the variables in which we take the Fourier transform, the two operations commute. Hence
\begin{eqnarray*}
|||D_rf(\Lambda) \psi^{(e)}|||_s^2&=&\int_{\R^{2n-1}}\int_{-\infty}^{0}(1+|\xi|^2)^sf((1+|\xi|^2)^{1/2})^2\\
 &~~&~~~~~~~~~~~~~~~~\times (1+|\xi|^2)e^{2r(1+|\xi|^2)^{1/2}}|\tilde{\psi}(\xi, 0)|^2drd\xi\\
&=&\frac{1}{2}\int_{\R^{2n-1}}f((1+|\xi|^2)^{1/2})^2(1+|\xi|^2)^{s+1/2}|\tilde{\psi}(\xi, 0)|^2d\xi\\
&=&\frac{1}{2}\no{f(\Lambda)\psi}_{b,s+\frac{1}{2}}^2.
\end{eqnarray*}

The lemma is proved.

$\hfill\Box$\\

{\it Proof of Proposition \ref{1ES}. } Instead for the  system $\{L_j\}_{j\le q_o}\cup \{\bar L_j\}_{q_o+1\le j\le n}$, we prove the theorem for a general elliptic system  $\{S_j\}_{1\le j\le n}$ where the $S_j$'s are vector fields in $ \C T \C^n$. This means that  there exist no $\eta\in T^*\C^n\setminus\{0\}$ such that $\la S_j,\eta \ra_x=0$ for all $j=1,. . .,n$.\\

Assume for a moment that $\varphi(x,0)=0$, (i.e, $\varphi|_{b\Om}=0$), which implies that $ \Lambda^{-1}f(\Lambda) \varphi(x,0)=0$ since the boundary condition is invariant under the action of tangential operators. Now,  \eqref{index1} holds for a general elliptic system; this yields
\begin{eqnarray}
\begin{split}\no{\Lambda^{-1}f(\Lambda) \varphi}_1^2\lesssim& \sum_{j=1}^{n}\big(\no{S_j\Lambda^{-1}f(\Lambda) \varphi}^2+\no{\bar S_j\Lambda^{-1}f(\Lambda) \varphi}^2\big)+\no{\Lambda^{-1}f(\Lambda)\varphi}^2\\
\lesssim&
\sum_{j=1}^{n}\no{S_j\Lambda^{-1}f(\Lambda) \varphi}^2+\no{\Lambda^{-1}f(\Lambda)\varphi}^2.
\end{split}
\end{eqnarray}

 Next, consider a general $\varphi$. Let $\varphi_b$ be the restriction of $\varphi$ to boundary, that is, $\varphi_b(x)=\varphi(x,0)$.  We set $\varphi^{(0)}=\varphi-\varphi^{(e)}$. Then,   $\varphi^{(0)}$ vanishes on the boundary so that the previous result applies to $\varphi^{(0)}$. We then have
\begin{eqnarray}
\begin{split}
\no{f(\Lambda) \varphi^{(0)}}^2\lesssim \sum_{j=1}^{n}\no{S_j\Lambda^{-1}f(\Lambda) \varphi^{(0)}}^2&+\no{\Lambda^{-1}f(\Lambda)\varphi^{(0)}}^2.
\end{split}
\end{eqnarray}
Therefore
\begin{eqnarray}\Label{2.38}\begin{split}
\no{\La^{-1}f(\Lambda)\varphi}_1^2\lesssim&\no{f(\Lambda)\varphi^{(0)}}^2+\no{f(\Lambda)\varphi^{(e)}}^2\\
\lesssim&\sum^{n}_{j=1}\no{S_j\Lambda^{-1}f(\Lambda) \varphi^{(0)}}^2+\no{\Lambda^{-1}f(\Lambda)\varphi^{(0)}}^2+\no{f(\Lambda)\varphi^{(e)}}^2\\
\lesssim&\sum^{n}_{j=1}\Big(\no{S_j\Lambda^{-1}f(\Lambda) \varphi }^2+\no{S_j\Lambda^{-1}f(\Lambda) \varphi^{(e)}}^2\Big)\\
&+\no{\Lambda^{-1}f(\Lambda)\varphi}^2+\no{\Lambda^{-1}f(\Lambda)\varphi^{(e)}}^2+\no{f(\Lambda)\varphi^{(e)}}^2\\
\lesssim&\sum^{n}_{j=1}\no{S_j\Lambda^{-1}f(\Lambda)\varphi}^2+\no{D_r\La^{-1}f(\Lambda)\varphi^{(e)}}^2\\
&+\no{\Lambda^{-1}f(\Lambda)\varphi}^2+\no{\Lambda^{-1}f(\Lambda)\varphi^{(e)}}^2+\no{f(\Lambda)\varphi^{(e)}}^2,
\end{split}
\end{eqnarray}
since the $S_j$'s are linear combinations of the tangential derivatives $\frac{\di}{\di x_j}$'s and of $D_r$ . Using Lemma \ref{l2.10}, 
$\no{D_r\Lambda^{-1}f(\Lambda)\varphi^{(e)}}^2$ and $\no{f(\Lambda)\varphi^{(e)}}^2$ are also estimated in the same way.
 This concludes  the proof.\\

$\hfill\Box$\\



\chapter{ $(f\T-\M)^k$-estimate : a sufficient condition} 
\label{Chapter3}

In this chapter,  we give the proof of Theorem \ref{main1}.

\section{Reduction  to the boundary}\noindent

Since  the Property $(f\T-\M\T-P)^k$ only concerns $u^\tau$, the tangential component of $u$,  we first prove the
$(f\T-\M)$-estimate for $u^\tau$, that is,
\begin{eqnarray}\Label{3.1}
\no{f(\La)\M u^\tau}^2\lesssim Q(u^\tau, u^\tau)+C_\M|||u^\tau|||^2_{-1}.
\end{eqnarray}

As a first step, we reduce the estimate \eqref{3.1} to the boundary by the following
\begin{theorem}\Label{t3.1}Let $\Om$ be  $q$-pseudoconvex (resp. $q$-pseudoconcave) at a boundary point $z_0$. Then,    there is a neighborhood $U$ of $z_0$ such that
\begin{eqnarray}\Label{3.3}
\no{f(\Lambda)\M u^\tau}^2\lesssim Q(u^\tau,u^\tau)+C_\M|||u^\tau|||^2_{-1}+||\Lambda^{-1/2}f(\Lambda)\M  u_b^\tau||_b^2,
\end{eqnarray}
 for any $u\in C_c^\infty(U\cap \bar{\Om})^k\cap \T{Dom}(\bar{\partial}^*)$ and $k\ge q$ (resp. $k\le q$ ).
\end{theorem}

Before the proof of Theorem \ref{t3.1}, we need following lemma  
\begin{lemma}\Label{l3.2}
Let $\Om$ be the $q$-pseudoconvex (resp. $q$-pseudoconcave) at boundary point $z_0$. Then,    there is a neighborhood $U$ of $z_0$ such that 
\begin{eqnarray}
\sumJ\Big(\sum_{j=1}^{q_o}\no{L_j u_J}^2+\sum_{j=q_o+1}^{n}\no{\bar L_j  u_J}^2+\no{D_r\La^{-1}u_J}^2+\no{u_J}^2\Big)\lesssim Q(u,u)
\end{eqnarray}
 for any $u\in C_c^\infty(U\cap \bar{\Om})^k\cap \T{Dom}(\bar{\partial}^*)$ with $k\ge q$ (resp. $k\le q$ ).
\end{lemma}
{\it Proof. } We only need to show that
$$\sumJ \no{D_r\La^{-1}u_J}^2\lesssim Q(u,u).$$
Since $D_r=a\bar L_n+bT$ where $a,b\in C^\infty(\bar\Om)$ and $T$ is a tangential operator of order 1, then 
$$\no{D_r\La^{-1}u_J}^2\lesssim \no{\bar L_nu_J}^2+\no{T\La^{-1}u}^2\lesssim \no{\bar L_nu_J}^2+\no{u}^2\lesssim Q(u,u).$$
$\hfill\Box$\\

{\it Proof of Theorem \ref{t3.1}. } 
 Applying Proposition \ref{1ES} for  $\varphi=\M u^\tau\in C_c^\infty(U\cap\bar\Om)$,  we get
\begin{eqnarray}
\begin{split}
\no{ f(\Lambda) \M u^\tau}^2\lesssim &\sum_{j=1}^{q_o}\no{L_j\La^{-1}f(\La)  \M u^\tau}^2+\sum_{j=q_o+1}^{n}\no{\bar L_j \La^{-1}f(\La) \M u^\tau}^2\\
&+\no{\La^{-1}f(\La) \M u^\tau}^2+\no{\La^{-1/2}f(\La)\M  u^\tau_b}_b^2.
\end{split}
\end{eqnarray}
Let $S\in \{L_j\}_{j\le q_0}\cup \{L_j\}_{q_0+1\le j\le n}$.\\

{\bf Case 1.} If $\M=M\rho$, where $\rho\in\A^{0,0}$ and $M$ is either every positive munber or 1, then
$$\M u^\tau :=M|\rho||u^\tau|;$$ it follows
\begin{eqnarray}
\begin{split}
\no{S\La^{-1}f(\La)  \M u^\tau}^2=&M^2 \sumJ\no{S\La^{-1}f(\La)  \rho u^\tau_J}^2\\
\lesssim& M^2 \sumJ\Big(\no{\rho \La^{-1}f(\La)S u_J^\tau}^2\\
&+\no{\rho[S,\La^{-1}f(\La)]u_J^\tau}^2+\no{[S\La^{-1}f(\La),\rho ] u_J^\tau}^2 \Big)\\
\lesssim& M^2 \sumJ\Big(\no{\La^{-1}f(\La) Su_J^\tau}^2+ \no{\La^{-1}f(\La) u_J^\tau}^2+\no{D_r\La^{-2}f(\La) u_J^\tau}^2 \Big)\\
\lesssim &\sumJ \Big(\no{ Su_J^\tau}^2+\no{u_J^\tau}^2+\no{D_r\La^{-1}u^\tau_J}+C_M|||u_J^\tau|||^2_{-1}\Big)\\
\lesssim &Q(u^\tau,u^\tau)+C_M|||u_J^\tau|||^2_{-1},\\
\end{split}
\end{eqnarray}
where the second inequality follows from Proposition \ref{p2.8}, the fourth from Lemma \ref{l2.7} and the last from Lemma \ref{l3.2}.\\

{\bf Case 2.} Let $\M=M\theta$, where $\theta\in\A^{0,1}$ and $M$ is either every positive number or $1$.  Remember  that 
\begin{eqnarray}
\begin{split}
\M u^\tau :=&M\sqrt{\big |\sumK |\sum_{j=1}^{n}\theta_j u^\tau_{jK}|^2-\sumJ\sum_{j=1}^{q_o} |\theta_j|^2 | u^\tau_J|^2\big|}
\end{split}
\end{eqnarray}
In the same way as in case 1, we get 
$$\no{S\La^{-1}f(\La)  \M u^\tau}^2\lesssim Q(u^\tau,u^\tau)+C_M|||u_J^\tau|||^2_{-1}.$$

This completes the proof of Theorem \ref{t3.1}.

$\hfill\Box$\\

\begin{remark} We notice again that we put $C_\M:=C_M$ if $M$ is every positive number and $C_\M:=0$ if $M=1$.
\end{remark}

\section{Estimate on the strip}
In this section, we show the Property $(f\T-\M\T-P)^k$ implies a related  estimate on the $\delta$-strip near the boundary for each $\delta>0$.\\  

Recall that $S_\delta:=\{z\in\C^n: -\delta<r<0\}.$
\begin{theorem}\Label{t3.4} Let $\Om$ be $q$-pseudoconvex (resp. $q$-pseudoconcave) at the boundary point $z_0$.
Assume that property $(f\T-\M\T-P)^k$ holds at $z_0$ with $k\ge q$ (resp. $k\le q$). Then,   there is a neighborhood $U$ of $z_0$ such that  for  any $\delta>0 $
\begin{eqnarray}\Label{strip}
f(\delta^{-1})^2\int_{S_{\delta/2}}| \M u^\tau|^2dV \lesssim Q(u^\tau,u^\tau)
\end{eqnarray} 
holds for any $u\in C^{\infty}_c(U\cap \bar\Om)^k\cap \T{Dom}(\bar{\partial}^*)$.
\end{theorem}
To simplify, define the quadratic form $H_{q_o}^\Phi(u)$ by \index{Notation ! $H^\Phi_{q_o}(u)$}
$$H_{q_o}^\Phi(u)=\sumK\underset{ij=1}{\overset{n}\sum}\phi_{ij}u_{iK}\bar u_{jK}-\underset{j=1}{\overset{q_o}\sum}\phi_{jj} |u|^2$$

{\it Proof of Theorem \ref{t3.4}. }
The proof is divided into two steps. In step 1, we  modify $\Phi^{\delta,\M}$ to $\tilde\Phi^{\delta,\M}$ which has the property that we need not only on the strip but also on the whole $\Omega$. In step 2, we  prove the estimate \eqref{strip}.\\

{\bf Step 1. } By our assumption,  for any $\delta>0$ sufficiently small, there is a  function $\Phi^{\delta,\M}$ such that
\begin{eqnarray}\Label{3.8}
\begin{cases}
H_{q_o}^{\Phi^{\delta,\M}}(u^\tau)\ge c\Big( f(\delta^{-1})^2|\M u^\tau |^2+ \sum_{j=1}^{q_o}|\Phi^{\delta,\M}_ju^\tau|^2\Big)\\
|\Phi^{\delta,\M}|\le1
\end{cases} \T{ on }  U\cap S_\delta
\end{eqnarray}
where $\Phi^{\delta,\M}_j=L_j(\Phi^{\delta,\M})$.\\

Now,  $u^\tau$ has support in  $U$ but the  properties of $\Phi^{\delta,\M}$ only hold on the strip $S_\delta$. Therefore, we have to modify $\Phi^{\delta,\M}$ to $\tilde\Phi^{\delta,\M}$ to get  the properties on the whole $ U$.\\

We define 
\begin{eqnarray}\tilde\Phi^{\delta,\M}:=\Phi^{\delta,\M}\chi(-\frac{r}{\delta}),
\end{eqnarray}
where $\chi$ is a cut off function which satisfies
 $\chi(t)=\begin{cases}1 &\T{ for } t \le \frac{1}{2} \\
0&\T{ for } t\ge 1, \end{cases}$ \\
and $\dot\chi\le 0$. \\

Computation of $\di\dib \tilde\Phi^{\delta,\M} $ shows that 
\begin{eqnarray}
\begin{split}
\di\dib\tilde\Phi^{\delta,\M}=&\chi\di\dib\Phi^{\delta,\M}-\frac{\dot\chi \Phi^{\delta,\M}}{\delta}\di\dib r
-2\T{ Re }\frac{\dot\chi}{\delta}\di\Phi^{\delta,\M}\otimes\dib r+\frac{\Phi^{\delta,\M}\ddot\chi}{\delta^2}\di r\otimes\dib r.
\end{split}
\end{eqnarray}
Notice that, if $i$ or $j$ is $\leq n-1$, then
$$\tilde\Phi_{ij}^{\delta,\M}=\chi \Phi_{ij}^{\delta,\M}-\frac{\dot{\chi}\Phi^{\delta,\M}}{\delta}r_{ij}.$$
We remark that $\langle\partial r,u^\tau\rangle\equiv0$ and therefore
\begin{eqnarray}
\begin{split}\Label{3.11}
H_{q_o}^{\tilde\Phi^{\delta,\M}}(u^\tau)&=\chi H_{q_o}^{\Phi^{\delta,\M}}(u^\tau)-\frac{\dot\chi \Phi^{\delta,\M}}{\delta}H^r_{q_o}(u^\tau).
\end{split}
\end{eqnarray}
 We note now that we can write $r=2Rez_n+h(z_1,. . .,z_{n-1}, y_n)$ for a graphing local defining function. Denote by $z\to z^*$ the projection $\C^n\to b\Om$ in a neigborhood of $z_0$ along the $x_n$-axis. We have   the evident equality $(r_{ij}(z))^{n-1}_{ij}=(r_{ij}(z^*))^{n-1}_{ij=1}$. Thus the second term in right hand side of \eqref{3.11} can be discarded  since $\Om$ is  $q$-pseudoconvex. (Recall here that $\dot\chi\le0$.)\\

Combining with \eqref{3.8},  we have
\begin{eqnarray}\Label{3.12}
\begin{split}
H_{q_o}^{\tilde\Phi^{\delta,\M}}(u^\tau)&\ge \chi H_{q_o}^{\Phi^{\delta,\M}}(u^\tau)\\
&\ge c\chi\Big(f(\delta^{-1})^2|\M u^\tau|^2+\sum_{j=1}^{q_o}|\Phi^{\delta,\M}_ju^\tau|^2\Big)\\
&\ge c\Big(\chi f(\delta^{-1})^2|\M u^\tau|^2+\sum_{j=1}^{q_o}|\tilde\Phi^{\delta,\M}_ju^\tau|^2\Big)
\end{split}
\end{eqnarray}
for $z\in U\cap\Om$. Here last inequality follows from $(\tilde\Phi^{\delta,\M})_j=\chi(\Phi^{\delta,\M})_j$ for $j\le q_o$ and $\chi\ge\chi^2$.\\

{\bf Step 2. } We apply Proposition~\ref{kmh} for $\phi=\psi(\tilde\Phi^{\delta,\M})$ and $u=u^\tau$.
First, we remark that
\begin{eqnarray}
\Label{3.13}
\begin{split}
H_{q_o}^{\phi}(u^\tau) = & \dot\psi H_{q_o}^{\tilde\Phi^{\delta,\M}}(u^\tau)\\
&+\ddot\psi\Big(\sumK|\underset{j=1}{\overset{n-1}\sum}\tilde\Phi^{\delta,\M}_ju^\tau_{jK}|^2-\underset{j=1}{\overset{q_0}\sum}|\tilde\Phi^{\delta,\M}_j|^2|u^\tau|^2\Big).
\end{split}
\end{eqnarray}
We also have
\begin{eqnarray}\Label{3.14}
\begin{split}
\no{\dib^*_{\phi}u^\tau}^2_{\phi}\le& 2\no{\dib^*u^\tau}^2_{\phi}+2\no{\sumK\sum^{n-1}_{j=1} \phi_ju^\tau_{jK}\bom_K}^2_{\phi}\\
=&2\no{\dib^*u^\tau}^2_{\phi}+2\sumK\no{\sum^{n-1}_{j=1} \dot{\psi}\tilde\Phi^{\delta,\M}_ju^\tau_{jK}}^2_{\phi}.
\end{split}
\end{eqnarray}
Thus we get from \eqref{KMH} applied to $u^\tau$, under the choice of the weight $\phi=\psi(\tilde\Phi^{\delta,\M})$  and taking into account \eqref{3.13} and \eqref{3.14} 
\begin{equation}
\Label{3.15}
\begin{split}
2\no{\bar{\partial} u^\tau}^2_{\phi}&+\,4\no{\bar{\partial}^*u^\tau}^2_{\phi}+C\no{u^\tau}^2_{\phi}\\
&\ge\,\int_\Om\dot\psi e^{-\phi}H_{q_o}^{\tilde\Phi^{\delta,\M}}(u^\tau)dV
\\
&+\, \int_\Om(\ddot\psi-4\dot\psi^2) e^{-\phi}\sumK|\underset{j=1}{\overset{n-1}\sum}\tilde\Phi^{\delta,\M}_ju^\tau_{jK}|^2dV\\
&-\int_\Om\ddot\psi  e^{-\phi}\underset{j=1}{\overset{q_o}\sum} |\tilde\Phi^{\delta,\M}_j|^2|u^\tau|^2 dV.
\end{split}
\end{equation}
We now specify our choice of $\psi$. First, we  want $\ddot\psi\ge 4\dot\psi^2$ so that the term in the third line of \eqref{3.15} can be disregarded.  Keeping this condition, 
we need an opposite estimate which assures that the absolute value of the second negative term in the last line of  \eqref{3.15} is controlled by one half of the  term in the second line. In fact, from \eqref{3.12}, we have
\begin{equation}
\begin{split}
\frac{1}{2}\int_\Om\dot\psi e^{-\phi}H_{q_o}^{\tilde\Phi^{\delta,\M}}(u^\tau)dV-\int_\Om\ddot\psi  e^{-\phi}\underset{j=1}{\overset{q_o}\sum} |\tilde\Phi^{\delta,\M}_j|^2|u^\tau|^2 dV\\
\ge \int_\Om e^{-\phi}\Big(\frac{1}{2}c\dot\psi -\ddot\psi\Big)\underset{j=1}{\overset{q_o}\sum} |\tilde\Phi^{\delta,\M}_j|^2|u^\tau|^2 dV.\\
\end{split}
\end{equation}
 The above term is nonnegative as soon as $\ddot\psi\le \frac{c}{2}\dot\psi$ .  If we then set $\psi:=\frac12e^{\frac c{2}(t-1)}$ then both requests are satisfied.
 Thus the inequality of \eqref{3.15} continues as
\begin{equation}
\begin{split}
\Label{3.17}
\quad\quad&\ge \,\frac{1}{2}\int_\Om\dot\psi e^{-\phi}H_{q_o}^{\tilde\Phi^{\delta,\M}}(u^\tau)dV\\
&\ge\,\int_\Om\frac{1}{2}\dot\psi e^{-\phi}cf(\delta^{-1})^2\chi(-\frac{r}{\delta})|\M u^\tau|^2dV\\
&\ge\frac{c}{2}f(\delta^{-1})^2\int_{S_\delta/2}\dot\psi e^{-\phi}|\M u^\tau|^2dV.
\end{split}
\end{equation} 
Here  the first inequality comes from \eqref{3.12} and the last equality follows from $\chi(\frac{-r}{\delta})=1$ on  $S_{\delta/2}$.

Now, we want to remove the weight from the resulting inequality. The weight in the first line of \eqref{3.15} can be handled owing to 
 $e^{-\phi}\le 1$ on $\bar \Om\cap U$. Furthermore, since $|\tilde\Phi^{\delta,\M}|<1$  on $U$, then the term $\dot \psi  e^{-\phi}$ in the last line of \eqref{3.17} is bounded from below by  a positive constant. We end up with the unweighted estimate
\begin{eqnarray}
\Label{3.18}
\no{\bar{\partial} u^\tau}^2+\no{\bar{\partial}^* u^\tau}^2+\no{u^\tau}^2
\simge f(\delta^{-1})^2\int_{S_{\delta/2}}|\M u^\tau|^2dV.
\end{eqnarray}

This concludes the proof of the theorem.\\

$\hfill\Box$\\

\section{Proof of Theorem 1.10}

In this section,  we give the proof of   Theorem \ref{main1}.  We first prove it for $u$ replaced by $u^\tau$.
\begin{theorem}\Label{t3.5}
Let $\Om$ be $q$-pseudoconvex (resp. $q$-pseudoconcave) at $z_0\in b\Om$ and
assume that Property $(f\T-\M\T-P)^k$ holds at $z_0$. Then,  we have
\begin{eqnarray}\Label{3.19}
\no{f(\La)\M u^\tau}^2 \lesssim Q(u^\tau,u^\tau)+C_\M |||u^\tau|||_{-1}^2
\end{eqnarray} 
for any $u^\tau\in C^{\infty}_c( U\cap\bar\Om)^k\cap \T{Dom}(\bar{\partial}^*)$ with $k\ge q$ (resp. $k\le q$).
\end{theorem}

For the proof of Theorem~\ref{t3.5}, we use the method by \cite{C87}.\\
Let $\{p_k\}$ with ${k=0,1,\dots }$ be a sequence of cutoff functions with the properties
\begin{enumerate}
  \item $\sum^\infty_{k=0}p^2_k(t)\cong1$; for any $t\ge 0$
  \item $p_k(t)\equiv 0$ if $t\not\in(2^{k-1},2^{k+1})$ with $k\ge 1$ and $p_0(t)\equiv 0$ if $ t\ge 2$.
  \end{enumerate}
We can also choose $p_k$ so that $p_k'(t)\lesssim 2^{-k}$.

Let $P_k$ denote the operator defined by 
$$(\widetilde{P_k\varphi})(\xi,r)=p_k(|\xi|){\varphi}(\xi,r)$$ for any $\varphi\in C^\infty_c$. To prove inequality \eqref{3.19}, we need the following 
\begin{lemma}\Label{l3.6} Let $z_0\in b\Om$, $f\in \mathcal{F}$, $a\in S(\R_{-}^{2n})$ and let $S$ be a vector field. Then, for any smooth $\varphi$ with support in a neighborhood of $z_0$, we have
\begin{enumerate}
  \item[(i)] $\no{f(\Lambda)\varphi}^2\cong\sum_{k=0}^{\infty}f(2^{k})^2\no{P_k\varphi }^2;$
  \item[(ii)] $\sum^\infty_{k=0}f(2^{k})^2||[P_k,a]\varphi||^2\lesssim ||\Lambda^{-1}f(\Lambda)\varphi||^2;$
  \item[(iii)] $\sum^{\infty}_{k=0} \no{[P_k,S]\varphi}^2\lesssim \no{D_r\La^{-1} \varphi}^2+\no{\varphi}.$
\end{enumerate}
\end{lemma}
{\it Proof.}
{\it (i).}
  We have
\begin{eqnarray*}
\no{f(\Lambda)\varphi}^2&=&\int_{-\infty}^{0}\int_{\R^{2n-1}}f(\Lambda_\xi)^2|\tilde{\varphi}(\xi, r)|^2d\xi dr\\
&=&\int_{-\infty}^{0}\int_{\R^{2n-1}}f(\Lambda_\xi)^2\Big(\sum_{k=0}^{\infty}p_k^2(|\xi|)\Big)|\tilde{\varphi}(\xi, r)|^2d\xi dr\\
\end{eqnarray*}
since $\sum_{k=0}^\infty p_k^2=1$. We notice that $\Lambda_\xi=(1+|\xi|^2)^{1/2}\cong 2^{k}$ as long as $|\xi|$ is in the support of $p_k$. Thus, it follows
 \begin{eqnarray*}
\no{f(\Lambda)u}^2&\cong&\sum_{k=0}^{\infty}\int_{-\infty}^{0}\int_{\R^{2n-1}}f(2^{k})^2|p_k(\tau)\tilde{\varphi}(\xi, r)|^2d\xi dr\\
&=&\sum_{k=0}^{\infty}f(2^{k})^2\no{P_ku}^2.
\end{eqnarray*}

{\it (ii). }  We can choose another sequence of cutoff functions $\{q_k\}_{k=0}^{\infty}$ such that \begin{enumerate}
  \item $q_k\equiv1$ on supp$(p_k)$ ;
  \item $q_k(t)\equiv 0$ if $t\not\in(2^{k-2},2^{k+3})$ with $k\ge 1$ and $q_0(t)\equiv 0, t\ge 4$.\end{enumerate}
Then \begin{eqnarray}\Label{3.20}
 |p_k(x)-p_k(y)|\lesssim 2^{-k}q_k(y)|x-y| 
\end{eqnarray}
 for any $x,y\ge 0$.
Observe that 
\begin{eqnarray}
\begin{split}
\widetilde{([P_k,a]\varphi)}(\xi,r)=\int_{\R^{2n-1}}\big(p_k(|\xi|)-p_k(|\tau|)\big)\tilde{a}(\xi-\tau,r)\tilde \varphi(\tau,r)d\tau.
\end{split}
\end{eqnarray}
Using \eqref{3.20}, and Plancherel Theorem together with Young's Inequality, we get
\begin{eqnarray}\Label{3.22}
\begin{split}
\no{[P_k,a]\varphi}^2 \lesssim 2^{-2k}\no{Q_ku}^2
\end{split}
\end{eqnarray}
where $Q_k$ is the operator with symbol $q_k(|\xi|)$.
Multiplying \eqref{3.22} by $f(2^k)^2$, taking summation over $k=0,1,...$ and using  the result of (i) for $f(|\xi|)$ replaced by $|\xi|^{-1}f(|\xi|)$ and $P_k$ by $Q_k$, we get the conclusion of (ii).\\

{\it (iii).}
The proof of (iii) follows immediately from (ii) by observing that $S=\sum a_j \frac{\di}{\di x_j}+b \frac{\di}{\di r}$.

$\hfill\Box$\\

\noindent
{\it Proof of Theorem~\ref{t3.5}. } By Theorem~\ref{t3.1}, we only need to estimate $\no{\Lambda^{-1/2}f(\Lambda)(\M u^\tau)_b}^2_b$.
Let $\chi_k\in C^\infty_c(-2^{-k},0]$  with $0\leq\chi_k\leq1$ and $\chi_k(0)=1$. \\ 

We have  the elementary inequality 
\begin{equation}
\Label{new*}
|g(0)|^2\le \frac{2^{k}}{\eta }\int^0_{-2^{-k}}|g(r)|^2dr+2^{-k}\eta\int^0_{-2^{-k}}|g'(r)|^2dr,
\end{equation}
which holds for any $g$ such that $g(-2^{-k})=0$.
If we apply it for $g(r)=\no{\chi_k(r)P_k \M u^\tau(\cdot,r)}_b$, we get
\begin{eqnarray*}
\no{\Lambda^{-1/2} f(\Lambda)(\M u^\tau)_b}_b^2&\simeq&\sum_{k=0}^{\infty}f(2^{k})^2 2^{-k}\no{\chi_k(0)P_k\M u^\tau(\cdot,0)}^2_b\\
&\le&\eta^{-1}\underbrace{ \sum_{k=0}^{\infty} f(2^{k})^2\int_{-2^{-k}}^0\no{\chi_k P_k\M u^\tau(.,r)}^2dr}_{I}\\
&&+\eta\underbrace{\sum_{k=0}^{\infty}f(2^{k})^2 2^{-2k}\int_{-2^{-k}}^0 \no{D_r\Big(\chi_k P_k\M u^\tau(.,r)\Big)}^2dr}_{II},\\
\end{eqnarray*}
where the first ``$\simeq$" follows from (i) of Lemma~\ref{l3.6} and the subsequent ``$\leq$" from \eqref{new*}.
Observe that $\chi_k\le 1$ and recall Theorem~\ref{t3.4}  that we apply for $P_k\M u^\tau$ and $\delta=2^{-k}$.  Thus the first sum (I) above can be estimated by
\begin{eqnarray}
\begin{split}
(I) &\lesssim\sum_{k=0}^{\infty} Q(P_ku^\tau,P_ku^\tau)\\
&\lesssim\sum_{k=0}^{\infty} \no{P_k\dib  u^\tau}^2+\no{P_k\dib^*u^\tau}^2+\no{[P_k,\dib]u^\tau}^2+\no{[P_k,\dib^*]u^\tau}^2\\
 &\lesssim Q(u^\tau,u^\tau)+\no{\Lambda^{-1} D_ru^\tau}^2,
\end{split}
\end{eqnarray}
where the estimate of the commutators  follows from Lemma~\ref{l3.6}. We remark that $D_r u^\tau$ can be expressed as a linear combination of $\bar{L}_n u^\tau$ and $Tu^\tau$ for some tangential vector field $T$. Then, similarly as in the proof of Lemma~\ref{l3.2}, we have
\begin{eqnarray*}
 ||\Lambda^{-1}D_r(u^\tau)||^2&\lesssim& ||\Lambda^{-1}\bar{L}_n u^\tau||^2+||\Lambda^{-1} T u^\tau||^2\\
&\lesssim& \no{\bar{L}_nu^\tau}^2+||u^\tau||^2\\
&\lesssim& Q(u^\tau,u^\tau).
\end{eqnarray*}Therefore, $$(I)\lesssim Q(u^\tau,u^\tau).$$
 
We now estimate (II). Since $D_r(\chi_k)\le 2^k$, $D_rP_k=P_kD_r$ and $\chi_k\le 1$, we get\begin{eqnarray}\Label{3.24}
\begin{split}
(II)&\le\sum_{k=0}^{\infty}f(2^{k})^2 2^{-2k}\left(\int_{-2^{-k}}^0 \no{D_r(\chi_k)P_k \M u^\tau(.,r)}^2dr\right.\\
&\left.+\int_{-2^{-k}}^0 \no{\chi_k D_r(P_k\M u^\tau(.,r))}^2dr\right)\\
&\le\sum_{k=0}^{\infty}f(2^{k})^2\int_{-2^{-k}}^0 \no{P_k\M u^\tau(.,r)}^2dr\\
&+\sum_{k=0}^{\infty}f(2^{k})^2 2^{-2k}\int_{-2^{-k}}^0 \no{ P_kD_r(\M u^\tau(.,r))}^2dr\\
&\le\no{f(\La) \M u^\tau}^2+\no{f(\La)\La^{-1} D_r(\M u^\tau)}^2,
\end{split}
\end{eqnarray}
where the last inequality follows from Lemma \ref{l3.6}. We now estimate the second term in the last line of \eqref{3.24}. Since $D_r=a\bar L_n+bT$, as before, we obtain
\begin{eqnarray*}
 \no{f(\La)\La^{-1} D_r(\M u^\tau)}^2 &\lesssim& 
||f(\Lambda)\Lambda^{-1}\bar{L}_n \M u^\tau||^2+||f(\Lambda)\Lambda^{-1} T\M u^\tau||^2\\
&\lesssim& \no{\bar{L}_nu^\tau}+C_\M|||u^\tau|||^2_{-1}+||f(\Lambda)\M u^\tau||^2\\
&\lesssim& Q(u^\tau,u^\tau)+C_\M|||u^\tau|||^2_{-1}+||f(\Lambda)\M u^\tau||^2.
\end{eqnarray*}

  Combining all our estimates of $\no{f(\Lambda)\Lambda^{-1/2}(\M u^\tau)_b}_b^2$, we obtain
$$\no{f(\Lambda)\Lambda^{-1/2}\M u^\tau_b}_b^2\lesssim \eta^{-1}\Big(Q(u^\tau,u^\tau)+\no{u^\tau}^2+C_\M|||u^\tau|||^2_{-1}\Big)+\eta ||f(\Lambda)\M u^\tau||^2.$$
Summarizing up, we have shown that
$$||f(\Lambda)\M u^\tau||\lesssim \eta^{-1}Q(u^\tau,u^\tau)+C_\M\no{u^\tau}^2_{-1}+\eta||f(\Lambda)\M u^\tau||^2.$$
Choosing $\eta>0$ sufficiently small, we can move the term $\eta||f(\Lambda)\M u^\tau||^2$ into the left-hand-side and get $$||f(\Lambda)\M u^\tau||^2\lesssim Q(u^\tau,u^\tau)+C_\M |||u^\tau|||^2_{-1}.$$ The proof is complete.\\

$\hfill\Box$

The proof of Theorem~\ref{main1} follows from Theorem~\ref{t3.5}  and the following 

\begin{lemma}\Label{l3.7}
Let $\Om$ be  $q$-pseudoconvex (resp. $q$-pseudoconcave ) at $z_0$ and $U$ be a neigborhood of $z_0$. For each $u\in C^\infty(U\cap\bar\Om)^k\cap \T{Dom}(\dib^*)$ with $k\ge q$ (resp. $k\le q$), assume that $(f\T-\M)^k$  holds for $u^\tau$. Then,   $(f\T-\M)^k$ holds for $u$. 
\end{lemma}
{\it Proof. }  Observe that $u^\nu|_{b\Om}\equiv0$; from \eqref{2.26} it follows  
$$Q(u^\nu,u^\nu)\lesssim \no{u^\nu}_1^2\lesssim Q(u,u).$$

On the other hand, 
$$\begin{cases}
\no{\dib u^\tau}=\no{\dib (u-u^\nu)}\le \no{\dib u}+\no{\dib u^\nu}\\
\no{\dib^* u^\tau}=\no{\dib^* (u-u^\nu)}\le \no{\dib^* u}+\no{\dib^* u^\nu}.\\
\end{cases}$$
Hence 
\begin{eqnarray}\Label{3.25}
Q(u^\tau,u^\tau)\le Q(u,u)+Q(u^\nu,u^\nu)\lesssim Q(u,u).
\end{eqnarray}
Therefore, 
\begin{eqnarray}
\begin{split}
\no{f(\La)\M u}\lesssim&\no{f(\La)\M u^\tau}+\no{f(\La)\M u^\nu}\\
\lesssim& Q(u^\tau,u^\tau)+C_\M|||u^\tau|||_{-1}^2+\no{u^\nu}^2_1+C_\M|||u^\nu|||^2_{-1}\\
\lesssim& Q(u,u)+C_\M|||u|||^2_{-1}.\\
\end{split}
\end{eqnarray}
Finally, we remark that
\begin{eqnarray}
\begin{split}
C_\M|||u|||^2_{-1}\lesssim&\no{D_r\La^{-1}u}^2\tilde{C}_{\M}\no{u}_{-1}^2\\
\lesssim& Q(u,u)+\tilde C_\M\no{u}_{-1}^2.
\end{split}
\end{eqnarray}

This completes the proof of Lemma \ref{l3.7}.

$\hfill\Box$

\section{Some remarks about $(f\T-\M)^k$ }
In this section we make some remarks about $(f\T-\M)^k$.\\

We recall that we write $f\gg g$ when $\underset{\xi\to\infty}\lim\frac{f(|\xi|)}{g(|\xi|)}=+\infty$.
\begin{lemma}
 Assume that $(f\T-\M)^k$ holds for $f\gg g$. Then for any $\epsilon>0$,  $(g\T-\frac{1}{\epsilon}M)^k$ also  holds.
\end{lemma}
The proof of the lemma follows from Lemma~\ref{l2.7}.\\

For example, if $f\gg \log$, then  $(f\T-1)^k$ estimate implies the superlogarithmic estimate. Similarly, if $f\gg 1$, then the $(f\T-1)^k$  estimate implies  the  compactness estimate.

\begin{lemma}
Let $\M' =\M$ on $b\Om$. Assume that $(f\T-\M)^k$ holds; then $(f\T-\M')^k$ also  holds.
\end{lemma}
 Using \eqref{elliptic}, we get the proof of this lemma.\\

\begin{lemma}\Label{l3.38}
If the $(f\T-\M)^k$ estimate holds, then we have
$$\no{D_r\Lambda^{-1}f(\Lambda)(\M u)}^2\lesssim Q(u,u)+C_\M |||u|||_{-1}^2,$$
for any $u\in C^\infty_c(U\cap \bar\Om)^k\cap \T{Dom}(\dib^*)$. 
\end{lemma}
{\it Proof. }
There are functions $a$ and $b$ such that 
$$D_r=a\bar L_n+ a T,$$
where $T$ is a tangential operator of order one. Therefore
\begin{eqnarray}
\begin{split}
\no{\Lambda^{-1}f(\Lambda)\frac{\di}{\di r}\M u}^2\lesssim &\no{\Lambda^{-1}f(\Lambda)\bar L_n (\M u)}^2+\no{\Lambda^{-1}f(\Lambda)T \M u}^2\\
\lesssim &\no{\bar L_n u}^2+\no{f(\Lambda) \M u}^2\\
\lesssim& Q(u,u)+C_\M\no{u}_{-1}^2.
\end{split}
\end{eqnarray}
This completes  the proof of the lemma.

$\hfill\Box$\\




It is interesting to remark that when $\Om$ is pseudoconvex (resp. pseudoconcave), then   $(f\T-\M)^k$ implies $(f\T-\M)^{k+1}$ (resp. $(f\T-\M)^{k-1}$).
Notice that the similar result is not clear when $\Om$ is $q$-pseudoconvex (resp. $q$-pseudoconcave) in case $q_o\neq0$ (resp. $q_o\neq n-1$). 

\begin{lemma}\Label{l3.11}
Let $\Om$ be  pseudoconvex (resp. pseudoconcave) at  $z_0\in b\Om$ and assume that $(f\T-\M)^k$ holds. Then,   $(f\T-\M)^{k+1}$  (resp. $(f\T-\M)^{k-1}$) also holds.
\end{lemma}

{\it Proof.} {We first discuss the pseudoconvex case.}  Let 
$u=\underset{|L|=k+1}{{\sum}'}u_L\bar \om_L$ have degree $k+1$.
We rewrite $u$ as a non-ordered sum
$$u=\frac{1}{(k+1)!}\underset{|L|=k+1}{{\sum}'}u_L\bar \om_L=\frac{(-1)^k}{k+1}\sum_{l=1}^n\Big(\frac{1}{k!}{{\sum}'}_{|J|=k}u_{lJ}\bom_J \Big)\we\bom_l.$$
 For $l=1,\dots, n$, we define a set of $k$-forms $v_l$ by $v_l:=\underset{|J|=k}{\sum'}u_{lJ}\bar\om_{J}$. It is easy to see that  $\overset{n}{\underset{l=1}{\sum}} |v_l|^2=(k+1)|u|^2$  and $\overset{n}{\underset{l=1}{\sum}}\sumK r_{ij}(v_l)_{iK}\overline{({v_l})}_{jK}=k\sumJ r_{ij}u_{iJ}\bar u_{jJ}$.\\ 

Using  \eqref{2.27} we have 

\begin{eqnarray}
\begin{split}
\sum_{l=1}^n Q(v_l,v_l)&\simeq \sum_{l=1}^n\Big( \sum _{j=1}^n \no{\bar L_j v_l}^2+\sum_{ij}\sumK\int_{b\Om}r_{ij}(v_l)_{iK}\overline{(v_l)}_{jK}dS+\no{vl}^2\Big )\\
&= (k+1)\sum _{j=1}^n \no{\bar L_j u}^2+k\sum_{ij}\sumJ \int_{b\Om}r_{ij}u_{iJ}\overline{u}_{jJ}dS+(k+1)\no{u}^2\\
&\simeq(k+1)Q(u,u).
\end{split}
\end{eqnarray}
\\
If $\M\in \A^{0,0}$, then
$$\sum_{l=1}^n\no{f(\La) \M v_l}^2=(k+1)\no{f(\La) \M u}^2.$$
If $\M\in \A^{1,0}$, we have
$$\sum_{l=1}^n| \M v_l|^2=\sum_{l=1}^n\sumK |\M_j (v_l)_{jK} |^2=k\sumJ|\M_j u_{jJ}|^2=k|\M u|^2;$$
then
$$\sum_{l=1}^n\no{f(\La) \M v_l}^2=k\no{f(\La)\M u}^2.$$
We discuss now the pseudoconcave case. Let  $u=\sumK u_K\bom_K$ have degree $k-1$. For l=1,...,n, define $v_l=\sumK u_K\bom_K\we \bom_l$; this has degree $k$. Then the proof follows the same lines as in the pseudoconvex case.

$\hfill\Box$

So far, only local estimates have been treated. If one goes to the full boundary, one has to consider local estimates over a covering. In fact, only on each local patch, the pseudodifferential operator $f(\Lambda)$ makes sense. Instead, the multiplier $\mathcal M$ has global meaning and thus a ``global $\M$ estimate" is meaningful. One readily checks that this comes as the collection of the local estimates.

\begin{theorem}\Label{t3.12}
Let $\Om$ be a bounded pseudoconvex domain. Assume that for multiplier $\M$  there exists a  family $\{\Phi^\M\}$ such that 
\begin{eqnarray}
\Label{supernew}
\begin{cases}
|\Phi^\M|\lesssim 1\\
\di\dib\Phi^\M(u^\tau,u^\tau)\ge  |\M u^\tau|^2
\end{cases}\T{~~~~on ~~}b\Om.
\end{eqnarray} 
Then, 
\begin{eqnarray}\Label{3.32}
\no{\M u}^2\lesssim \no{\dib u}^2+\no{\dib^* u}^2+C_\M\no{u}^2_{-1}
\end{eqnarray}
for any $u\in C^\infty(\bar\Om)^k\cap \T{Dom}(\dib^*)$.
\end{theorem}

The proof of Theorem \ref{t3.12} follows from Theorem \ref{main1} and elliptic estimate for the interior of domain.
In fact, from \eqref{supernew} one gets, when $z$ leaves $b\Om$ and enters in $\Om$, an error term proportional to $-\no{ru^\tau}^2$. Since $\no{ru^\tau}^2_1\le Q(ru^\tau,ru^\tau)$ (by ellipticity in the interior),  the error term can be absorbed.

We notice that this can be obtained directly (without using Theorem \ref{main1}) by the global method such as in Catlin \cite{C84}.
\begin{remark}
The inequality \eqref{3.32} is trivial when $\M$ is bounded.
\end{remark}



\chapter{ The $(f\T-\M)^k$ estimate on the boundary } 
\label{Chapter4}


In this chapter we  study the behavior of the boundary value of forms associated to the $\dib$-Neumann problem. Precisely, we establish a relation between the $(f\T-\M)^k$-estimate on $\Om$ and $b\Om$.

\section{Definitions and notations}\noindent

Let $M$ be a smooth real hypersurface in $\C^n$. We start by denoting  by  $\A_b^{0,k}$  the space of restriction of element of $\A^{0,k}\cap \T{Dom}(\dib^*)$ to the boundary $b\Om$. Then $\A^{0,k}_b$ is the space of smooth sections of the vector bundle $\big (T^{0,1}(M)^*\big)^k$ on $M$.  \\

    The tangential Cauchy-Riemann operator $\dib_b : \A_b^{0,k}\to \A_b^{0,k+1} $ is defined as follows. If $u \in \A_b^{0,k}$,  let $u'$ be a (0,k)-form whose restriction to $M$  equals $u$. Then $\dib_b u$ is the restriction of $\bar\partial u'$ to $M$. \\

Let $z_0\in M$ and $U$ be a neigbourhood of $z_0$; we fix a defining function $r$ of $M$ such that $|\di r|=1$ on $U\cap M$. We assume that $\Om$ is one of the two sides of $\C^n\setminus M$ in a neighborhood of $z_0$. Let $L_1,...,L_n$ be the local basis for (1,0) vector fields defined in $U$, as defined in Chapter 1.\\

We can define a Hermitian inner product on $\A^{0,k}_b$ by 
$$(u,v)_b=\int_M\la u,v \ra dS,$$    
where $dS$ is the surface element on $M$. The inner product gives rise to an $L_2$-norm $\no{\cdot}_b$.\\

In analogy with Chapter 2,  we define $\dib^*_b$ to be the $L^2$-adjoint of $\dib_b$ in the standard way. Thus $\dib^*_b:\A_b^{0,k+1}\to \A_b^{0,k} $  for $k\ge 0$. The Kohn-Laplacian is defined by 
$$\Box_b=\dib_b\dib^*_b+\dib^*_b\dib_b.$$ \index{Notation ! $\Box_b$ }
Remember that we have already defined 
$$Q_b(u,v)=(\dib_b u,\dib_b v)_b+(\dib^*_b u,\dib^*_b v)_b+( u,v)_b.$$

\index{Notation ! $Q_b$ }

 Denote $C_c^\infty(U\cap M)^k$ the forms of $\A^{0,k}_b$ with  compact support in $U$. Let
$u=\sumJ u_J\bom_J\in C^\infty_c(U\cap M)^k;$
then on $M$, the operator $\dib_b$ and $\dib^*_b$ are expressed by
\begin{equation}
\dib_b u=\sumJ\sum_{j=1}^{n-1} \bar L_j u_J\bom_j\we\bom_J+...
 \end{equation}
and 
\begin{eqnarray}
\dib_b^*u=-{\sumK}\sum_j^{n-1} L_ju_{jK}\bom_K+... 
\end{eqnarray}
where dots refer the error terms in which $u$ is not differentiated.  \\
\index{Space ! $\A^{0,k}_b$ }\index{Space ! $C_c^\infty(U\cap M)^k$ }\index{Notation ! $\dib_b$ }\index{Notation ! $\dib_b^*$ }

In $U$, we choose special boundary coordinate $(x_1,...,x_{2n-1},r)$. Let $\xi=$ $(\xi_1,...,\xi_{2n-1})$ $=(\xi',\xi_{2n-1})$ be the dual coordinates to $\{x_1,...,x_{2n-1}\}$. 
We also decompose $(x_1,...,x_{2n-1})$ $=(x',x_{2n-1})$ so that $T^\C_{z_0}M$ is defined by $x_{2n-1}=0$ in $T_{z_0}M$.
Let $\psi^++ \psi^-+ \psi^0=1$ be a $C^\infty$ partition of the unity in the sphere $|\xi|=1$ such that
$\psi^\pm$ are 1 at the ploles $(0,...,\pm 1)$ and $\psi^0$ at the equator, that is, at  $\xi_{2n-1}=0$. We extend these functions to $\R^{2n-1}\setminus\{0\}$ as homogeneous functions of degree 0. We may assume that the supports of the functions $\psi^+$, $\psi^-$ and $\psi^0$ are contained in the cones 
 \begin{eqnarray}\begin{split}
\mathcal C^+=&\{\xi \big |\xi_{2n-1}> \frac{1}{2} |\xi'| \};\\
\mathcal C^-=&\{\xi\big|-\xi_{2n-1}> \frac{1}{2} |\xi'|  \};\\
\mathcal C^0=&\{\xi \big | | \xi_{2n-1}|<|\xi'| \}.
\end{split}
\end{eqnarray}\index{Notation ! $\mathcal C^+, \mathcal C^-,\mathcal C^0$ } 
Then $\T{supp}{\psi^+}\subset\subset  \mathcal C^+$, $\T{supp}{\psi^-}\subset\subset  \mathcal C^-$ and  $\T{supp}{\psi^0}\subset\subset \mathcal C^0$.\\

The operators $\Psi=\Psi^{\overset\pm0}$ with symbols $\psi=\psi^{\overset\pm0}$ are defined by 
$$\widetilde{\Psi \varphi}(\xi)=\psi(\xi)\tilde{\varphi}(\xi) \qquad \T{for }\quad \varphi\in C^\infty_c(U\cap M);$$
$$\widetilde{\Psi \varphi}(\xi,r)=\psi(\xi)\tilde{\varphi}(\xi,r) \qquad \T{for }\quad \varphi\in C^\infty_c(U\cap \Om).$$
The microlocal decomposition  $\varphi=\varphi^++\varphi^-+\varphi^0$ of a function $\varphi\in C^\infty_c(U\cap M)$ is defined by \index{Notation ! $\Psi^+, \Psi ^-,\Psi^0$ } \index{Notation ! $u^+, u^-, u^0$ }
$$\varphi=\zeta \Psi^+\varphi+\zeta \Psi^-\varphi+\zeta \Psi^0\varphi,$$
 where $\zeta\in C^\infty(U'), \bar U\subset \subset U'$ and $\zeta\equiv1$ on $U$.  \\
 For a form $u$, the microlocal decomposition $u=u^++u^-+u^0$ is accordingly defined coefficientwise.
We recall some definitions from Chapter 1. 
\begin{definition} If $M$ is a hypersurface and $z_0\in M$, then a $(f\T-\M)^k_{b}$ estimate  for $(\dib_b,\dib_b^*)$ holds at $z_0$ if there exists a neighborhood $U$ of $z_0$ such that 
$$(f\T-\M)^k_{b}  \qquad\no{f(\Lambda)\M u}^2_b\le cQ_b(u,u)+C_\M\no{u}_{b,-1}^2$$ 
for all $u\in C_c^\infty(U\cap M)^{k}$. And  an $(f\T-\M)^k_{b,+}$- (resp. $(f\T-\M)^k_{b,-}$-) estimate  for $(\dib_b,\dib_b^*)$  holds at $z_0$ if $(f\T-\M)^k_{b}$ holds with $u$ replaced by $u^+$ (resp.  $u^-$), that is,  
$$(f\T-\M)^k_{b,+}\qquad \no{f(\Lambda)\M u^+}^2_b\le cQ_b(u^+,u^+)+C_\M\no{u^+}_{b,-1}^2$$ 
(resp.
$$(f\T-\M)^k_{b,-}\qquad \no{f(\Lambda)\M u^-}^2_b\le cQ_b(u^-,u^-)+C_\M\no{u^-}_{b,-1}^2.\quad )$$ 
 \end{definition}
\begin{definition}
The hypersurface $M$ is said to be $q$-pseudoconvex at $z_0$ if either of the two components of $\C^n\setminus M$  is $q$-pseudoconvex at $z_0$. 
\end{definition}

Denote by $\Om^+=\{z\in U | r(z)<0 \}$ the $q$-pseudoconvex side of $M$ and by $\Om^-$ the other one. By Remark \ref{convex-concave}, $\Om^-=\{z\in U | -r(z)<0 \}$ is $(n-q-1)$-pseudoconcave at $z_0$.  Remember that  we are choosing $\om_1,\dots,\om_{n-1},\om_n= \partial r$ as an orthonormal basis of (1,0)-forms. We also use the notation $\omega_n^\pm=\pm\partial r$ for the exterior conormal vectors to $\Om^\pm$.  We define $T=\frac{1}{2}(L_n-\bar L_n)$ and $\frac{\di}{\di r}=\frac{1}{2}(L_n+\bar L_n)$. It follows 
\begin{eqnarray}\Label{4.4}
\begin{split}
L_n=\frac{\di}{\di r}+T \T{~~~~and~~~~ }\bar L_n=\frac{\di}{\di r}-T.
\end{split}
\end{eqnarray}
\index{Notation ! $T$ }
\section{Basic microlocal  estimates on $b\Om$}
In this section we prove the basic microlocal estimates on hypersurface.\\
\index{Estimate ! basic microlocal}

In a  similar way as in Proposition~\ref{kmh}, we get
\begin{lemma} \Label{l4.1}For two indices  $q_1,q_2$ ; $(1\le q_1\le q_2\le n-1)$,  there is a constant $C$ such that 
\begin{eqnarray}\begin{split}\Label{4.3}
\no{\dib_b u}_b^2+&\no{\dib^*_b u}_b^2+ \no{ u}_b^2\\
\simgeq& \sumK\sum_{ij=1}^{n-1}(r_{ij}Tu_{iK},u_{iK})_b-\sumJ\sum_{j=q_1}^{q_2}(r_{jj}Tu_J,u_J)_b^2\\
&+\frac{1}{2}\sumJ\Big(\sum_{j=1}^{q_1-1}\no{\bar L_j u_J}^2_b+\sum_{j=q_2+1}^{n-1}\no{\bar L_j u_J}^2_b+\sum_{j=q_1}^{q_2}\no{ L_j u_J}^2_b\Big).
\end{split}
\end{eqnarray}
\end{lemma}
Note that, conversely, we have
\begin{eqnarray}\begin{split}\Label{4.4}
\no{\dib_b u}_b^2+&\no{\dib^*_b u}_b^2\\
\simleq& \left|\sumK\sum_{ij=1}^{n-1}(r_{ij}Tu_{iK},u_{iK})_b-\sumJ\sum_{j=q_1}^{q_2}(r_{jj}Tu_J,u_J)_b^2\right|\\
&+\sumJ\Big(\sum_{j=1}^{q_1-1}\no{\bar L_j u_J}^2_b+\sum_{j=q_2+1}^{n-1}\no{\bar L_j u_J}^2_b+\sum_{j=q_1}^{q_2}\no{ L_j u_J}^2_b\Big)+\no{u}_b^2,
\end{split}
\end{eqnarray}
 for all $u\in C_c^\infty(U\cap M)^k$ for any $k$.

\begin{lemma}  \Label{l4.4}
 Let $M$ be a  hypersurface and $z_0$ a point of $M$. Then there is a  neighborhood $U$ of $z_0$ such that
  $$ Q_b(u^0,u^0) \cong \no{ u^0}_{b,1}^2, $$
for all $u\in C_c^\infty(U\cap M)^k$ with any $k$. 
\end{lemma}
{\it Proof.} Using twice the inequality \eqref{4.3}  for $q_1=q_2=0$ and $q_1=0, q_2=n-1$ and taking summation, we get
\begin{eqnarray}\Label{4.5}
\begin{split}
\no{\dib_b u}_b^2+&\no{\dib^*_b u}_b^2+ \no{ u}_b^2\\
\simgeq& \sumK\sum_{ij=1}^{n-1}(r_{ij}Tu_{iK},u_{iK})_b-\sumJ\sum_{j=1}^{n-1}(r_{jj}Tu_J,u_J)_b^2\\
&+\frac{1}{2}\Big(\sumJ\sum_{j=1}^{n-1}\no{L_j u_J}^2_b+\sumJ\sum_{j=1}^{n-1}\no{\bar L_j u_J}^2_b\Big)\\
\simgeq& \no{\Lambda' u}^2_b-\big(\epsilon+\T{diam}(U)\big)\no{Tu}^2_b-C_\epsilon \no{u}^2_b,
\end{split}
\end{eqnarray}
where $\La'$ is the pseudodifferential operator of order 1 whose symbol is  $(1+\sum_{j=1}^{2n-2}|\xi_j|^2)^\frac{1}{2}$. 
To explain the last estimate in \eqref{4.5}, observe that $L_j|_{z_0}=\sqrt2\partial_{z_j},\,\,j=1,...,n-1$ and that the coefficients of the $L_j$'s are $C^1$;  therefore, the third line of \eqref{4.5} is bounded from below by $||\Lambda'u||^2-\T{diam}(U)\no{Tu}^2_b$.
Apply \eqref{4.5} for $u^0$ and notice that $\no{\Lambda'u^0}^2_b\simgeq \no{Tu^0}^2_b$. Taking $U$ and $\epsilon$ suitably small, we conclude
$$Q_b(u^0,u^0)\simge \no{\La' u^0}^2_b\simge \no{\La u^0}^2_b.$$ 
On the other hand, the converse inequality is always true. 

$\hfill\Box$

\begin{lemma}  \Label{l4.5}
Let $M$ be a $q$-pseudoconvex hypersurface at $z_0$. Then, for a neighborhood $U$ of $z_0$, and for $\zeta'\equiv1$ over $\T{supp}(\zeta)$ and $\psi^{\pm\,\prime}\equiv1$ on $\T{supp}\,\psi^\pm$, we have 
\begin{eqnarray}{\it (i) } \begin{split}  
 Q_b(u^+,u^+) &+\no{\zeta'\Psi^+u}_{b,-\infty}^2\\
\cong& \sumK\sum_{ij=1}^{n-1}(r_{ij}\zeta' R^+ u^+_{iK},\zeta'R^+u^+_{iK})_b\\
&-\sumJ\sum_{j=1}^{q_0}(r_{jj}\zeta'R^+u^+_J,\zeta'R^+u^+_J)_b^2\\
&+\sumJ\sum_{j=1}^{q_0}\no{L_j u_J^+}^2_b+\sumJ\sum_{j=q_0+1}^{n-1}\no{\bar L_j u^+_J}^2_b+\no{u^+}_b^2+\no{\Psi^+u}_{b,-\infty}^2,
\end{split}
\end{eqnarray}
for all $u\in C_c^\infty(U\cap M)^k$ with  $k\ge q$, where $R^+$ is the pseudodifferential operator of order 1 whose symbol is $\xi_{2n-1}^{\frac{1}{2}}\psi^{+\,\prime}(\xi)$. Similarly, we have 
 \begin{eqnarray}{\it  (ii) }\begin{split}  
 Q_b(u^-,u^-)&+\no{\zeta'\Psi^-u}_{b,-1}^2\\
 \cong&\sumJ\sum_{j=q_0+1}^{n-1}(r_{jj}\zeta' R^-u^-_J,\zeta' R^-u^-_J)_b \\
&- \sumK\sum_{ij=1}^{n-1}(r_{ij} \zeta' R^- u^-_{iK},\zeta' R^-u^-_{iK})_b\\
&+\sumJ\sum_{j=1}^{q_0}\no{\bar L_j u^-_J}^2_b+\sumJ\sum_{j=q_0+1}^{n-1}\no{ L_j u^-_J}^2_b+\no{u^-}_b^2+\no{\Psi^-u}_{b,-\infty}^2,
\end{split}
\end{eqnarray}
for any  $u\in C_c^\infty(U\cap M)^k$ with  $k\le n-1- q$, where $ R^-$ is the pseudodifferential operator of order 1 whose symbol is $(-\xi_{2n-1})^{\frac{1}{2}}\psi^{-\,\prime}(\xi)$.
\end{lemma}

 Recall that, since $M$ is $q$-pseudoconvex at $z_0$, then there is a defining function of $M$ which satisfies on $M$
$$\sumK\sum_{ij=1}^{n-1}r_{ij}u_{iK}\bar u_{jK}-\sum_{j=1}^{q_o}r_{jj}|u|^2\ge 0,$$
for any $u\in C_c^\infty(U\cap M)^k$ with $k\ge q$, and also satisfies on $M$
$$\sum_{j=q_o+1}^{n-1}r_{jj}|u|^2-\sumK\sum_{ij=1}^{n-1}r_{ij}u_{iK}\bar u_{jK}\ge 0,$$
for any $u\in C_c^\infty(U\cap M)^k$ with $k\le n-q-1$,
 where $U$ is a neighborhood of $z_0$.\\
{\it Proof.}  {\it (i). } 
We have
$$\varphi^+=\zeta \Psi^{+}\varphi=\zeta (\Psi^{+\,\prime})^2 \Psi^+\varphi = (\Psi^{+\,\prime})^2\zeta \Psi^+\varphi+[\zeta, (\Psi^{+\,\prime})^2]\Psi^+\varphi.$$
Since the supports of symbols of $\Psi^+$ and $[\zeta, (\Psi^{+\,\prime})^2]$ are disjoint, the operator $[\zeta, (   \Psi^{+\,\prime})^2]\Psi^+$ is of order $-\infty$ and we have
\begin{eqnarray}\Label{4.8}\begin{split}
(r_{ij}T\varphi^+, \varphi^+)_b=&(r_{ij}T\zeta \Psi^+\varphi, \zeta \Psi^+\varphi)_b\\
&=(r_{ij}T (\Psi^{+\,\prime})^2\zeta \Psi^+\varphi, \zeta \Psi^+\varphi)_b+O(\no{\Psi^+\varphi}^2_{b,-\infty})\\
&=((\zeta')^2 r_{ij}R^{+\,*}R^+\zeta \Psi^+\varphi, \zeta \Psi^+\varphi)_b+O(\no{\Psi^+\varphi}^2_{b,-\infty})\\
&=(r_{ij} \zeta' R^+\zeta \Psi^+\varphi,\zeta' R^+\zeta \Psi^+\varphi)_b\\
&+([(\zeta')^2 r_{ij},R^{+\,*}]R^+\zeta \Psi^+\varphi, \zeta \Psi^+\varphi)_b+O(\no{\Psi^+\varphi}^2_{b,-\infty}).
\end{split}\end{eqnarray}
From the pseudodifferential operator calculus we get
 \begin{eqnarray}\begin{split}
([\zeta^{\prime2} r_{ij},R^{+*}]R^+\zeta \Psi^+\varphi, \zeta \Psi^+\varphi)_b|\lesssim \no{\varphi^+}^2_b
\end{split}\end{eqnarray}
Substituting $u$ or $u_J$ for $\varphi$ in \eqref{4.8},   we obtain
\begin{eqnarray}\Label{4.10}\begin{split}
&\sumK\sum_{ij=1}^{n-1}(r_{ij}T u^+_{iK},u^+_{iK})_b-\sumJ\sum_{j=1}^{q_0}(r_{jj}Tu^+_J, u^+_J)_b^2\\
=& \sumK\sum_{ij=1}^{n-1}(r_{ij}\zeta' R^+ u^+_{iK},\zeta'R^+u^+_{iK})_b-\sumJ\sum_{j=1}^{q_0}(r_{jj}\zeta'R^+u^+_J,\zeta'R^+u^+_J)_b^2\\
&+O(\no{u^+}^2_b)+O(\no{\Psi^+u}^2_{b,-\infty}).
\end{split}\end{eqnarray}
Since $M$ is $q$-pseudoconvex, then  the sum in second line in \eqref{4.10} is nonnegative if $k\ge q$.  Thus the first part of Lemma \ref{l4.5}  is proven  by applying Lemma \ref{l4.1} to $u^+$ with $q_1=0$ and $ q_2=q_0$.\\

{\it (ii). } The proof of the second part is similar. We have to notice that
\begin{eqnarray}\begin{split}
&\sumK\sum_{ij=1}^{n-1}(r_{ij}T u^-_{iK},u^-_{iK})_b-\sumJ\sum_{j=q_0+1}^{n-1}(r_{jj}Tu^-_J, u^-_J)_b^2\\
=&-\sumK\sum_{ij=1}^{n-1}(r_{ij}\zeta'  R^- u^-_{iK},\zeta' R^-u^-_{iK})_b+\sumJ\sum_{j=q_0+1}^{n-1}(r_{jj}\zeta' R^-u^-_J,\zeta' R^-u^-_J)_b^2\\
&+O(\no{u^-}^2_b)+O(\no{\Psi^-u}^2_{b,-\infty}),
\end{split}\end{eqnarray}
By Remark \eqref{convex-concave}, the second line is nonnegative for any $k$-form $u$ with $k\le n-q-1$ and this concludes the proof. 

$\hfill\Box$\\

\section{Basic  microlocal estimates on $\Om^+$ and $\Om^-$}
\index{Estimate ! basic microlocal }
In this section,  we prove  the basic microlocal estimates on $\Om^+$ and $\Om^-$. We begin by introducing  the harmonic extension of a form from $b\Om$ to $\Om$ following  Kohn  \cite{K86} and \cite{K02}. \\

In terms of special boundary coordinate $(x,r)$, the operator $L_j$ can be written as
$$L_j=\delta_{jn}\frac{\di}{\di r}+\sum_k a_j^k(x,r) \frac{\di}{\di x_k}$$
for $j=1,. . ., n$. We define the tangential symbols of $L_j$, $1\le j\le n-1$, by
$$\sigma(L_j)((x,r),\xi)=-i \sum_k a_j^k(x,r) \xi_k,$$ 
 and 
$$\sigma(T)((x,r),\xi)=\frac{-i}{2} \sum_k \big(a_n^k(x,r)-\bar a_n^k(x,r)\big) \xi_k.$$
Note that $\sigma(T)$ is real. We set $$\sigma_b(L_j)(x,\xi)=\sigma({L_j})((x,0),\xi) \T{~~and~~ } \sigma_b(T)(x,\xi)=\sigma(T)((x,0),\xi)$$ 
and $$\mu(x,\xi)=\sqrt{\sum_j|\sigma_b(L_j)(x,\xi)|^2+ |\sigma_b(T)(x,\xi)|^2+1}.$$
Remember the notation $\Lambda_\xi=(1+|\xi|^2)^{\frac12}$;
in a neighborhood of $z_0$, we have $\mu(x,\xi)\sim\Lambda^1_\xi$. 
\\

Harmonic extension is defined as follows. Let $\varphi\in C_c^\infty(U\cap M)$; define $\varphi^{(h)}\in C^\infty(\{(x,r)\in \R^{2n}| r\le 0\})$ by 

  $$\varphi^{(h)} (x, r)=(2\pi)^{-2n+1}\int_{\R^{2n-1}}e^{i x\cdot\xi} e^{r\mu(x,\xi)}\tilde\varphi(\xi)d\xi,$$
  so that $\varphi^{(h)}(x,0)=\varphi(x)$.  
This extension is called  ``harmonic"  since $\triangle \varphi^{(h)} (x, r)$ has order 1 on $M$. In fact, we have 
\begin{eqnarray}\begin{split}
\triangle=&-\sum_{j=1}^n \frac{\di^2}{\di z_j \di \bar z_j}\\
=&-\sum_{j=1}^n L_j\bar L_j+ \sum_{k=1}^{2n-1}a^k(x,r)\frac{\di}{\di x_k}+a(x,r)\frac{\di}{\di r}\\
=&-\frac{\di^2}{\di^2 r}+T^2-\sum_{j=1}^{n-1} L_j\bar L_j+ \sum_{k=1}^{2n-1}b^k(x,r)\frac{\di}{\di x_k}+b(x,r)\frac{\di}{\di r}
\end{split}
\end{eqnarray}
since \eqref{4.4} implies that $L_n\bar L_n=\frac{\di^2}{\di^2 r}-T^2+D$, where $D$ is a first order operator. Hence if $(x,r)\in U\cap \bar\Om^+, $
 \begin{eqnarray}\begin{split}\Label{4.15b}
\triangle (\varphi^{(h)})(x,r)=(2\pi)^{(2n-1)}\int e^{ix\cdot \xi}e^{r\mu(x,\xi)}\big(p^1(x,r,\xi)+rp^2(x,r,\xi)\tilde\varphi(\xi,0)\big)d\xi.
\end{split}
\end{eqnarray}
For future use, we prepare the notation $P^1+rP^2$ for the pseudodifferential operator with symbol $p^1+rp^2$ which appears in the right of \eqref{4.15b}.
Along with \eqref{4.15b} we have 
$$L_j\varphi^{(h)}(x,r)=(L_j\varphi)^h(x,r)+E_j\varphi(x,r)$$
where 
$$E_j\varphi (x,r)= (2\pi)^{(2n-1)}\int e^{ix\cdot \xi}e^{r\mu(x,\xi)}\big(p_j^0(x,r,\xi)+rp_j^1(x,r,\xi)\tilde\varphi(\xi)\big)d\xi$$
and 
$$\bar L_j\varphi^{(h)}(x,r)=(\bar L_j\varphi)^h(x,r)+\bar E_j\varphi(x,r)$$
for $j=1,...,n-1$.

\begin{lemma}\Label{l4.10}
For any $k\in \Z$ with $k\ge 0$, $s\in \R$, and $f\in \F$ we have
\begin{enumerate}
  \item[(i)] $|||r^kf(\Lambda)\varphi^{(h)}|||_s\simle \no{f(\Lambda)\varphi}_{b, s-k-\frac{1}{2}}$,
  \item[(ii)] $|||D_rf(\Lambda)\varphi^{(h) }|||_s\simle \no{f(\Lambda)\varphi}_{b,s+\frac{1}{2}}$
  \end{enumerate}
for any $\varphi\in C^\infty_c(U\cap \bar\Om^+).$
\end{lemma}
{\it Proof. } We notice again that $\mu(x,\xi)\cong (1+|\xi|^2)^{\frac{1}{2}}$ over a small neighborhood of $z_0$, and then the proof of this lemma is similar to the proof  of Lemma \ref{l2.10}.\\

$\hfill\Box$\\

We define $\varphi_b$ to be the restriction of $\varphi\in C^\infty_c( U\cap \bar\Om^+)$ to the boundary. We have the elementary estimate
\begin{eqnarray}\Label{4.16b}
\no{\varphi_b}_{b, s}^2\lesssim |||\varphi|||_{s+\frac{1}{2}}^2+|||D_r\varphi |||^2_{s-\frac{1}{2}}.
\end{eqnarray}

The following lemma states the basic microlocal estimates on $\Om^+$.
\begin{lemma}\Label{l4.6} Let $\Om^+$ be  $q$-pseudoconvex  at  $z_0$. If U is  a sufficiently small  neighborhood of $z_0$, then
we have the three estimates which follow
 \begin{multline}\Label{4.12}|||\Psi^0\varphi|||^2_{1}\lesssim \sum_{j=1}^{q_o}\no{L_j \Psi^0\varphi}^2+\sum_{j=q_o+1}^{n}\no{\bar L_j\Psi^0\varphi}+\no{\Psi^0\varphi}^2\\
 \T{for any } \varphi \in C_c^\infty(U\cap \bar\Om^+),
 \end{multline}
\begin{multline}\Label{4.13}|||\Psi^-\varphi|||^2_{1}\lesssim \sum_{j=1}^{q_o}\no{L_j\Psi^- \varphi}^2+\sum_{j=q_o+1}^{n}\no{\bar L_j\Psi^-\varphi}+\no{\Psi^-\varphi}^2\\
\T{for any } \varphi \in C_c^\infty(U\cap \bar\Om^+),
\end{multline}
\begin{multline}\Label{4.14}|||\bar L_n\Psi^+\varphi^{(h)}|||^2_{\frac{1}{2}}\lesssim \sum_{j=1}^{q_o}\no{L_j \Psi^+\varphi}_b^2+\sum_{j=q_o+1}^{n-1}\no{\bar L_j\Psi^+\varphi}_b^2+\no{\Psi^+\varphi}_b^2\\
\T{for any } \varphi\in C_c^\infty(U\cap b\Om^+).
\end{multline} 
\end{lemma}

{\it Proof. } We start from \eqref{4.12}. Since $\T{supp }\psi^0\subset \mathcal C^0$, then we have 
\begin{eqnarray}\begin{split}
(1+|\xi'|^2) |\psi^0(\xi)|^2&\lesssim \left(1+\sum_{j=1}^{n-1}|\sigma(L_j)|^2((0,0),\xi)\right)|\psi^0(\xi)|^2
\\
&\simleq \left(1+\sum_{j=1}^{n-1}|\sigma(L_j)|^2((x,r),\xi)\right)|\psi^0(\xi)|^2\\
&+\left(\sum_{j=1}^{n-1}|\sigma(L_j)|^2((0,0),\xi)-\sum_{j=1}^{n-1}|\sigma(L_j)|^2((x,r),\xi)\right)|\psi^0(\xi)|^2
\\
&\simleq \left(1+\sum_{j=1}^{n-1}|\sigma(L_j)|^2((x,r),\xi)\right)|\psi^0(\xi)|^2\\
&+\T{diam}(\bar\Om^+\cap U)\left(1+\sum_{j=1}^{n-1}|\xi_j|^2\right) |\psi^0(\xi)|^2.
\end{split}
\end{eqnarray}
Hence
\begin{eqnarray}\Label{4.16}
\begin{split}
\no{\Lambda \Psi^0 \varphi }^2_1\lesssim& \sum_{j=1}^{q_0}\no{L_j\Psi^0\varphi}^2+\sum_{j=q_0+1}^{n}\no{\bar L_j\Psi^0\varphi}^2\\
&+\no{\Psi^0\varphi}^2+\T{diam}(U\cap \bar\Om^+) \sum_{j=1}^{2n-1} \no{D_j\Psi^0\varphi }^2.
\end{split}
\end{eqnarray}
 The estimate \eqref{4.12} follows from \eqref{4.16} by taking $U$ sufficiently small so that the last term is absorbed in the left hand side of the estimate.\\

We pass to prove \eqref{4.13}. For all $\varphi\in C_c^\infty(U\cap \bar \Om^+)$,  let $\varphi^{(h)}$ be the harmonic extension of $\varphi_b=\varphi|_{U\cap b\Om^+}$. We have
\begin{eqnarray}\Label{4.17}
|||\Psi^-\varphi|||_1^2
\lesssim \no{\Psi^-(\varphi - \varphi^{(h)})}_1^2+|||\Psi^-\varphi^{(h)}|||^2_1.
\end{eqnarray}

We estimate now $|||\Psi^-\varphi^{(h)}|||_1^2$; we have 
\begin{eqnarray}
\begin{split}
\bar L_n \Psi^-\varphi^{(h)} (x,r)=&(2\pi)^{(2n-1)}\int e^{ix\cdot \xi} e^{r (\mu(x,\xi))}\Big(\mu(x,\xi)-\sigma_b(T)(x,\xi)\\
&~~~~~~~~~~~~~~~~~~~~~~~~~~~~~~~~~~~+rp^1(x,\xi)\Big)\psi^-(\xi)\tilde{\varphi}(\xi,0)d\xi  
\end{split}
\end{eqnarray}
where $p^1(x,\xi)$ is the symbol which appears in \eqref{4.15b} and whose associated operator we have denoted by $P^1$.
Choosing $U$ sufficiently small we have  $\sigma(T)_b(x,\xi)\le 0$ when $\xi \in\T{supp}(\psi^-)\subset \mathcal C^-$.  Then, 
$$\mu(x,\xi)-\sigma_b(T)(x,\xi)\simge |\xi|+1.$$
It follows
\begin{eqnarray}\Label{4.19}
\begin{split}
|||\Psi^-\varphi^{(h)}|||_1^2\lesssim& \no{\bar L_n\Psi^-\varphi^{(h)}}^2+\no{rP^1\Psi^-\varphi^{(h)}}^2.\\
\end{split}
\end{eqnarray}
Applying Lemma \ref{l4.10} and inequality \eqref{4.16b} to  the second term in \eqref{4.19}, we get
\begin{eqnarray}\Label{4.20}
\begin{split}
\no{rP_1\Psi^-\varphi^{(h)}}^2 
\lesssim& \no{\Lambda^{-1/2}\Psi^-\varphi}_b^2\\
\lesssim& \no{\Lambda^{-1} D_r\Psi^-\varphi}^2+\no{\Psi^-\varphi}^2\\
\lesssim& \no{\bar L_n\Psi^-\varphi}^2+\no{\Psi^-\varphi}^2.
\end{split}
\end{eqnarray}
For the first term in \eqref{4.19}, we have 
\begin{eqnarray}\Label{4.21}
\begin{split}
\no{\bar L_n\Psi^-\varphi^{(h)}}^2
\lesssim& \no{\bar L_n\Psi^-(\varphi - \varphi^{(h)})}^2+\no{\bar L_n\Psi^-\varphi}^2\\
\lesssim& \no{\Psi^-(\varphi - \varphi^{(h)})}_1^2+\no{\bar L_n\Psi^-\varphi}^2.
\end{split}
\end{eqnarray}\\

Combining \eqref{4.17}, \eqref{4.19}, \eqref{4.20} and \eqref{4.21}, we get
\begin{eqnarray}
\begin{split}
|||\Psi^-\varphi|||_1^2\lesssim \no{\Psi^-(\varphi-\varphi^{(h)})}_1^2+\no{\bar L_n \Psi^-\varphi}^2+\no{\Psi^-\varphi}^2.
\end{split}
\end{eqnarray}
Finally, we estimate $\no{\Psi^-(\varphi-\varphi^{(h)})}_1^2$. Since $\Psi^-(\varphi-\varphi^{(h)})=0$ on $U\cap b\Om$, then
\begin{eqnarray}\Label{4.23}
\begin{split}
\no{\Psi^-(\varphi-\varphi^{(h)})}^2_1\lesssim&\no{\Delta\Psi^-(\varphi-\varphi^{(h)})}^2_{-1}\\
\lesssim &\no{\Delta\Psi^- \varphi}^2_{-1}+\no{\Delta\Psi^-\varphi^{(h)}}^2_{-1}\\
\lesssim &\sum_{j=1}^{q_0}\no{ \bar L_j L_j\Psi^-\varphi}^2_{-1}+\sum_{j=q_0+1}^{n}\no{ L_j\bar  L_j\Psi^-\varphi}^2_{-1}\\&+\no{P^1\Psi^-\varphi}_{-1}^2+ \no{(rP^2+P^1)\Psi^-\varphi^{(h)}}^2_{-1}\\
\lesssim &\sum_{j=1}^{q_0}\no{  L_j\Psi^-\varphi}^2+\sum_{j=q_0+1}^{n}\no{\bar  L_j\Psi^-\varphi}^2+\no{\Psi^-\varphi}^2.\\
\end{split}
\end{eqnarray}
Here the third inequality in \eqref{4.23} follows from \eqref{4.15b}.  This completes the proof of \eqref{4.17}.\\

We prove now \eqref{4.14}. For any $\varphi\in C_c^\infty(U\cap b\Om^+)$, we have 
\begin{equation}
\Label{4.23bis}
\bar L_n \Psi^+\varphi^{(h)} (x,r)=(2\pi)^{(2n-1)}\int e^{ix\cdot \xi} e^{r (\mu(x,\xi))}\Big(\mu(x,\xi)-\sigma_b(T)(x,\xi)+rp^1(x,\xi)\Big)\psi^+(\xi)\tilde{\varphi}(\xi,0)d\xi.
\end{equation}

Choosing $U$ sufficiently small we have $\sigma_b(T)(x,\xi)>0$ when $\xi\in \T{supp}\psi^+\subset\mathcal C^+$. So that 
$$\mu-\sigma_b(T)=\sum_{j=1}^{q_o}\frac{\sigma_b(\bar L_j)}{\mu+\sigma_b(T)}\sigma_b(L_j)+\sum_{j=q_o+1}^{n-1}\frac{\sigma_b( L_j)}{\mu+\sigma_b(T)}\sigma_b(\bar L_j).$$
Since the symbols $$\Big\{\frac{\sigma_b(\bar L_j)}{\mu+\sigma_b(T)}\Big\}_{j\le q_o}\T{~~~  and ~~~~}\Big\{\frac{\sigma_b( L_j)}{\mu+\sigma_b(T)}\Big\}_{ q_o+1\le j\le n-1}$$ are absolutely bounded, then
\begin{equation}
\Label{4.23ter}
|\mu-\sigma_b(T)|\simleq |\sum_{j=1}^{q_o}\sigma_b(L_j)|+|\sum_{j=q_o+1}^{n-1}\sigma_b(\bar L_j)|.
\end{equation}
Hence, from \eqref{4.23bis} and \eqref{4.23ter} we get

\begin{eqnarray}
\begin{split}
|||\bar L_n \Psi^+ \varphi^{(h)}|||_{\frac{1}{2}}^2\lesssim &\sum_{j=1}^{q_o}|||(L_j\Psi^+\varphi)^h|||_{\frac{1}{2}}+\sum_{j=q_o+1}^{n-1}|||(L_j\Psi^+\varphi)^h|||_{\frac{1}{2}}+|||rP_1\Psi^+\varphi^{(h)}|||_{\frac{1}{2}}^2\\
\lesssim &\sum_{j=1}^{q_o}||L_j\Psi^+\varphi||_b^2+\sum_{j=q_o+1}^{n-1}||L_j\Psi^+\varphi||_b^2+\no{\Psi^+\varphi}^2_b.
\end{split}
\end{eqnarray}

$\hfill\Box$

Using Lemma~\ref{l4.6} for coefficients of forms, we obtain 

\begin{lemma}\Label{l4.7} 
Let $\Om^+$ be a $q$-pseudoconvex at  $z_0$. Then, for a suitable neighborhood $U$ of $z_0$ and for any $u\in C_c^\infty(U\cap \bar\Om^+)^k\cap \T{Dom}(\dib^*)$ with  $k\ge q$, we have
  \item[(i)] 
  $$|||\Psi^0u|||_1^2+|||\Psi^-u|||_1^2\lesssim Q(u,u).$$
Moreover, 
for any $u\in C_c^\infty(U\cap b\Om^+)^k$ with  $k\ge q$, we have
$$|||\bar L_n \Psi^+ (u^{+})^{(h)}|||_{\frac{1}{2}}^2\lesssim Q_b(u^+,u^+).$$
\end{lemma}

Similarly, we get the basic microlocal estimates for $\Om^-$.

\begin{lemma}\Label{l4.8} 
Let $\Om^-$ be $(n-1-q)$-pseudoconcave at $z_0$. Then, for a suitable neighborhood $U$ of $z_0$ and for any $u\in C_c^\infty(U\cap \bar\Om^-)^k\cap \T{Dom}(\dib^*)$ with $k\le  n-1-q$, we have
  $$|||\Psi^0u|||_1^2+|||\Psi^+u|||_1^2\lesssim Q(u,u).$$
Moreover, for any  $u\in C_c^\infty(U\cap b\Om^-)^k$  with $k\le n-1-q$, we have
$$|||\bar L_n \Psi^- (u^{-})^{(h)}|||_{\frac{1}{2}}^2\lesssim Q_b(u^-,u^-).$$
\end{lemma}

\section{The equivalence of $(f\T-\M)^k$ estimate on $\Om$ and $b\Om$}
In this section, we give the proof  of Theorem \ref{main2}.  This is a consequence of the three theorems which follow, that is, Theorem \ref{t4.10}, \ref{t4.11} and \ref{t4.12}.
\begin{theorem}\Label{t4.10}
Let $\Om^+\subset\C^n$ be a smooth $q$-pseudoconvex domain   with boundary $M=b\Om$ at $z_0\in b\Om$. Then
\begin{enumerate}
  \item[(i)] $(f\T-\M)^k_{\Om^+}$ implies $(f\T-\N)^k_{b,+}$  where $\N$ is the restriction of  $\M$ to $M$.
  \item[(ii)] $(f\T-\N)^k_{b,+}$ implies $(f\T-\M)^k_{\Om+}$ where $\M$ is any extension of $\N$ from $M$ to $\Om$, that is, $\M\big|_M=\N$.
\end{enumerate}
\end{theorem}

{\it Proof.} {\it (i). } We need to show that over a neighborhood $U$ of $z_0$ we have
$$\no{f(\La)\N u^+}_b^2\lesssim Q_b(u^+,u^+)+C_\N \no{u^+}^2_{b,-1}+\no{\Psi^+u}^2_{b,-\infty}$$
for any $u\in C^\infty_c(U\cap M)^k$.   Let $\chi=\chi(r)$ be a cut off function  with $\chi(0)=1$.   Applying inequality \eqref{4.16b}, we have
\begin{eqnarray}\Label{4.24}
\begin{split}
\no{f(\Lambda)\N u^+}^2_b\lesssim& \no{\Lambda^{\frac{1}{2}}\chi f(\Lambda)\M u^{+\,(h)}}^2 +\no{\Lambda^{-\frac{1}{2}}  D_r( \chi f(\Lambda)\M u^{+\,(h)})}^2\\
\lesssim&\no{f(\Lambda)\M\chi\zeta'  R^+  u^{+\,(h)}}^2\\
& +\no{\Lambda^{-1} f(\Lambda) D_r(\M\chi\zeta'R^+  u^{+\,(h)})}^2+error,
\end{split}
\end{eqnarray}
where $\zeta'=1$ on supp$(u^+)$ and $\T{supp}(\chi\zeta')\subset\subset U'$. Here, the error term is estimated by
$$error\lesssim \no{\La^{\frac{1}{2}} u^{+\,(h)}}^2+C_\M\no{\La^{\frac{1}{2}} u^{+\,(h)}}^2_{-1}+\no{\Psi^+u^{(h)}}^2_{-\infty}.$$
Notice that $\chi\zeta'  R^+  u^{+\,(h)}\in \T{Dom}(\dib^*)$. Using the hypothesis of the theorem to estimate the second  line of \eqref{4.24} and applying Lemma \ref{l3.38} to the second term in the last line of \eqref{4.24}, we have that \eqref{4.24} can be continued by
\begin{eqnarray}
\begin{split}
\lesssim& Q(\chi\zeta'  R^+  u^{+\,(h)}, \chi\zeta'  R^+  u^{+\,(h)})+C_\M\no{\chi\zeta'  R^+  u^{+\,(h)}}_{-1}^2+error\\
\lesssim& \sumK\sum_{ij}^{n-1}(r_{ij}\zeta'  R^+  u^{+}_{iK}, \zeta'  R^+  u^{+}_{jK})_b-\sumJ(r_{jj}\zeta'  R^+  u_J^{+}, \zeta'  R^+  u_J^{+})_b\\
&+\sumJ\Big(\sum_{j=1}^{q_0}\no{L_j \chi\zeta'  R^+  u_J^{(h)+}}^2+ \sum_{j=1}^{n}\no{\bar L_j\chi\zeta'  R^+  u_J^{(h)+}}^2\Big)\\
&+\no{\chi\zeta'  R^+  u^{+\,(h)}}^2+C_\M\no{\chi\zeta'  R^+  u^{+\,(h)}}_{-1}^2+error\\
\lesssim& Q_b(u^+,u^+)+C_\M\no{u^{+}}_{b,-1}^2+\no{\Psi^+u}^2_{b,-\infty}\\
&+\sumJ\Big(\sum_{j=1}^{q_0}\no{\Lambda^{\frac{1}{2}}(L_j u_J^{+})^{(h)}}^2+ \sum_{j=1}^{n-1}\no{\Lambda^{\frac{1}{2}} (\bar L_j u_J^{+})^{(h)}}^2\Big)+\no{\La^{\frac{1}{2}} \bar L_n \Psi^+(u^{+})^{(h)}}^2\\
\lesssim& Q_b(u^+,u^+)+C_\M\no{u^{+}}_{b,-1}^2+\no{\Psi^+u}^2_{b,-\infty},
\end{split}
\end{eqnarray}
where the second inequality follows from \eqref{2.27}, the third from Lemma \ref{l4.5} and the last from Lemma \ref{l4.7}.\\

{\it  (ii). }  For any $u\in C^\infty_c(U\cap \bar\Om)^k\cap \T{Dom}(\dib^*)$, we decompose $u=u^\tau+u^\nu$ and  $u^\tau=u^{\tau\,+   }+u^{\tau\,-   }+u^{\tau\,0   }$. Since $u^\nu$ satisfies elliptic estimates and on account of  Lemma \ref{l4.7}, we have
\begin{eqnarray}
\begin{split}
\no{f(\Lambda)\M u^\nu}^2\le |||u^\nu|||_1^2+C_\M\no{u^\nu}_{-1}^2&\lesssim Q(u,u)+C_\M\no{u}_{-1}^2,\\
\no{f(\Lambda)\M u^{\tau\,0   }}^2\le |||u^{\tau\,0   }|||_1^2+C_\M\no{u^{\tau\,0   }}_{-1}^2&\lesssim Q(u,u)+C_\M\no{u}_{-1}^2,\\  
\no{f(\Lambda)\M u^{\tau\,-   }}^2\le |||u^{\tau\,-   }|||_1^2+C_\M\no{u^{\tau\,-   }}_{-1}^2&\lesssim Q(u,u)+C_\M\no{u}_{-1}^2.\\    
\end{split}
\end{eqnarray}
Moreover, by Theorem \ref{t3.1} and \eqref{3.25}, we have
\begin{eqnarray}
\begin{split}
\no{f(\Lambda)\M u^{\tau\,+   }}^2\lesssim& \no{\Lambda^{-\frac{1}{2}}f(\Lambda)\N u^{\tau\,+   }_b}^2_b+Q(u^\tau,u^\tau)+C_\M\no{u^\tau}_{-1}^2\\    
\lesssim& \no{\Lambda^{-\frac{1}{2}}f(\Lambda)\N u^{\tau\,+   }_b}^2_b+Q(u,u)+C_\M\no{u}_{-1}^2.\\    
\end{split}
\end{eqnarray}
Thus, we obtain
\begin{eqnarray}
\begin{split}
\no{f(\Lambda)\M u}^2\lesssim &\no{f(\Lambda)\M u^{\tau\,+   }}^2+\no{f(\Lambda)\M u^{\tau\,-   }}^2+\no{f(\Lambda)\M u^{\tau\,0   }}^2+\no{f(\Lambda)\M u^\nu}^2\\
\lesssim& \no{\Lambda^{-\frac{1}{2}}f(\Lambda)\N u^{\tau\,+   }_b}^2_b+Q(u,u)+C_\M\no{u}_{-1}^2.\\    
\end{split}
\end{eqnarray}

Hence,  we only need to estimate $\no{\Lambda^{-\frac{1}{2}}f(\Lambda)\N u^{\tau\,+   }_b}^2_b$.
We begin by noticing that
\begin{eqnarray}\Label{4.30}
\begin{split}
\no{\Lambda^{-\frac{1}{2}}f(\Lambda)\N u^{\tau\,+   }_b}^2_b\lesssim \no{f(\Lambda)\N \zeta'\Lambda^{-\frac{1}{2}}u^{\tau\,+   }_b}^2_b+error,
\end{split}
\end{eqnarray}
where $\zeta'\equiv1$ over $\T{supp}u^{\tau\,+}$.
Using the hypothesis we can continue  \eqref{4.30} by
\begin{eqnarray}
\begin{split}
\lesssim&Q_b(\zeta'\Lambda^{-\frac{1}{2}}u^{\tau\,+   }_b,\zeta'\Lambda^{-\frac{1}{2}}u^{\tau\,+   }_b)+C_\M\no{\zeta'\Lambda^{-\frac{1}{2}}u^{\tau\,+   }_b}^2_{b,-1}+error\\
\lesssim&\sumK\sum_{ij=1}^{n-1}(r_{ij}\zeta'R^+(\zeta'\Lambda^{-\frac{1}{2}}u^{\tau\,+   }_b)_{iK},\zeta'R^+(\zeta'\Lambda^{-\frac{1}{2}}u^{\tau\,+   }_b)_{jK})_b\\
&\sumJ\Big(\sum_{j=1}^{q_0}\no{L_j(\zeta'\Lambda^{-\frac{1}{2}}u^{\tau\,+   }_b)_J}_b^2+ \sum_{j=1}^{n-1}\no{\bar L_j( \zeta'\Lambda^{-\frac{1}{2}}u^{\tau\,+   }_b)_J}_b^2\Big)\\
&+\no{\zeta'\Lambda^{-\frac{1}{2}}u^{\tau\,+   }_b}^2_b+C_\M\no{\zeta'\Lambda^{-\frac{1}{2}}u^{\tau\,+   }_b}^2_{b,-1}+\no{\Psi^+\Lambda^{-\frac{1}{2}}u^{\tau}_b}^2_{b,-\infty}+error\\
\lesssim&Q(\zeta''R^+\zeta'\Lambda^{-\frac{1}{2}}\zeta \Psi^+u^{\tau},\zeta''R^+\zeta'\Lambda^{-\frac{1}{2}}\zeta \Psi^+ u^{\tau})\\
&\sumJ\Big(\sum_{j=1}^{q_0}\no{\zeta' L_j \Lambda^{-\frac{1}{2}}(u^{\tau\,+   }_b)_J}_b^2+ \sum_{j=1}^{n-1}\no{\zeta' \bar L_j \Lambda^{-\frac{1}{2}}(u^{\tau\,+   }_b)_J}_b^2\Big)\\
&+\no{\zeta' \Lambda^{-\frac{1}{2}}u^{\tau\,+   }_b}^2_b+\tilde C_\M\no{\Lambda^{-\frac{1}{2}}u^{\tau\,+   }_b}^2_{b,-1}.\\
\end{split}
\end{eqnarray}
Since $\zeta''R^+\zeta'\Lambda^{-\frac{1}{2}}\zeta \Psi^+$ is a tangential pseudodifferential operators of order zero, then 
 $$Q(\zeta''R^+\zeta'\Lambda^{-\frac{1}{2}}\zeta \Psi^+u^{\tau},\zeta''R^+\zeta'\Lambda^{-\frac{1}{2}}\zeta \Psi^+u^{\tau})\lesssim Q(u^{\tau},u^{\tau}).$$

To estimate the last two lines  we proceed as follows. For $j\le q_0$,  since $\zeta' L_j \Lambda^{-\frac{1}{2}}u^{+\,\tau}\in C^\infty_c(U\cap \bar\Om)$, then using inequality \eqref{4.16b}, we have
\begin{eqnarray}
\begin{split}
\no{\zeta' L_j \Lambda^{-\frac{1}{2}}u^{+\,\tau}}_b^2\lesssim&\no{\Lambda^{\frac{1}{2}}\zeta' L_j \Lambda^{-\frac{1}{2}}u^{+\,\tau}}^2+\no{\Lambda^{-\frac{1}{2}}\frac{\di}{\di r}\zeta' L_j \Lambda^{-\frac{1}{2}}u^{+\,\tau}}^2\\
\lesssim&\no{L_j u^{+\,\tau}}^2+\no{\frac{\di}{\di r}\Lambda^{-1}L_j u^{+\,\tau}}^2+error\\
\lesssim&\no{L_j u^{+\,\tau}}^2+\no{T\Lambda^{-1}L_j u^{+\,\tau}}^2+\no{\bar L_n\Lambda^{-1}L_j u^{+\,\tau}}^2+error\\
\lesssim&\no{L_j u^{+\,\tau}}^2+\no{\bar L_n u^{+\,\tau}}^2+error\\
\lesssim&Q(u^{+\,\tau},u^{+\,\tau})\lesssim Q(u^\tau,u^\tau).
\end{split}
\end{eqnarray}
The same estimate holds for each term $\no{\zeta'  \bar L_j \Lambda^{-\frac{1}{2}}u^{+\,\tau}}_b^2$ for $q_0+1\le j\le n-1$.
This concludes the proof of Theorem \ref{t4.10}.

$\hfill\Box$

Similarly, we get the equivalence of $(f\T-\M)^k$ on $\Om^-$ and $M$

\begin{theorem}\Label{t4.11}
Let $\Om^-$ be a smooth $q$-pseudoconcave domain at $z_0\in b\Om$. 
Then $(f\T-\M)^k_{\Om^-}$ is equivalent to $(f\T-\N)^k_{b,-}$ for $\M|_M=\N$.
\end{theorem}
Let $M$ be $q$-pseudoconvex. Recall that this means that one of the two sides of $\C^n\setminus M$ is $q$-pseudoconvex (or equivalently, the complementary side is $(n-1-q)$-pseudoconcave).
\begin{theorem}\Label{t4.12}
Let $M$ be a $q$-pseudoconvex hypersurface at $z_0$.  Then $(f\T-\N)^k_{b,+}$ holds if and only if $(f\T-\N)^{n-1-k}_{b,-}$ holds.
\end{theorem}

{\it Proof.} We define the local conjugate-linear duality map $F^k : \A^{0,k}_b\to \A^{0,n-1-k}_b$ as follows. If $u=\sumJ u_J\bom_J$ then 
$$F^{k}u=\sum \epsilon^{\{J ,J'\}}_{\{1,..., n-1\}} \bar u_J\bom_{J'} , $$
where $J'$ denotes the strictly increasing $(n-k-1)$-tuple consisting of all integers in $[1,n-1]$ which do not belong to $J$ and $\epsilon^{J,J'}_{\{1,n-1\}}$ is the sign of the permutation $\{J, J'\}\simto \{1,\dots, n-1\}$. \\

Since $\overline{(\varphi^+)}=(\overline\varphi)^-$,
then
$$F^{k}u^+=\sum \epsilon^{\{J ,J'\}}_{\{1,..., n-1\}} (\bar u)^-_J\bom_{J'}. $$

Also,
 $F^{n-1-k}F^ku^+=u^+$, $\no{F^ku^+}=\no{u^-}$ and finally
$$\dib_b F^k u^+=F^{k-1}\dib^*_b u^+ +\cdots$$
and 
$$\dib^*_b F^k u^+=F^{k+1}\dib_b u^+ +\cdots$$
where dots refers the term in which $u$ is not differentiated. Hence
$$Q_b(F^k u^+, F^k u^+) \cong Q_b(u^+,u^+).$$
On the other hand, we also have $\no{f(\Lambda) \N F^ku^+}^2=\no{f(\Lambda) \N u^-}^2$.

$\hfill\Box$

\begin{corollary}\Label{c4.12}
Let $M$ be a pseudoconvex hypersurface at $z_0$, let $\N \in C^\infty(M)$
and let $\M$ be any extension of $\N$ from $M$ to $\C^n$.
 Then,  for $1\le k\le n-2$, the estimate $(f\T-\N)^k_{b}$ holds if  one of the following conditions is satisfied
\begin{enumerate}
  \item $(f\T-\N)^k_{b,+}$ and $(f\T-\N)^k_{b,-}$ 
  
\item $(f\T-\N)^k_{b,+}$ and $(f\T-\N)^{n-1-k}_{b,+}$ 

  \item $(f\T-\N)^l_{b,+}$ for $l\le \min\{k,n-1-k\}$

  \item $(f\T-\N)^l_{b,-}$  for $l\ge \max\{k,n-1-k\}$

  \item $(f\T-\M)^k_{\Om^+}$ and $(f\T-\M)^k_{\Om^-}$ 

  \item $(f\T-\M)^k_{\Om^+}$ and $(f\T-\M)^{n-1-k}_{\Om^+}$ 

  \item $(f\T-\M)^l_{\Om^+}$ for $l\le \min(k,n-1-k)$

  \item $(f\T-\M)^l_{\Om^-}$  for $l\ge \max(k,n-1-k)$.
\end{enumerate}

\end{corollary}

{\it Proof. }
The equivalence of $(f\T-\N)^k_{b}$ to (1), (2), (3) and (4) follows from Theorem~\ref{t4.12} combined with the fact that $\no{u^0}^2_{1\,b}\simleq Q_b(u,u)$. The equivalence of these conditions to (5), (6), (7) and (8) follows from Theorems \ref{t4.10} and \ref{t4.11}. 

$\hfill\Box$



\chapter{Geometric and analytic consequences of $(f\T-M)^k$ estimate }
\label{Chapter6}

The chapter is devoted to the proof of Theorem~\ref{main3}.

\section{Estimates with a $1$-parameter family of cutoff functions}
Let $\Om$ be a domain in $\C^n$ and let $U$ be a neighborhood of a given point $z_0\in b\Om$. Choose real coordinates $(x_1,...,x_{2n-1},r)$ in $U$.  Let $\chi_0 ,\chi_1,\chi_2$ be cutoff functions supported in $U$ such that $\chi_j\equiv 1$ in a neighborhood of the support of $\chi_{j-1}$ for $j=1,2$. For every $t\in (0,1]$, we define  two new cutoff functions $\chi^t_j$ such that $\chi_j^t(z)=\chi_j(\frac{z}{t})$ for $j=0,1,2$. The main tool in the proof of Theorem~\ref{main3} is  the following 
\begin{theorem}
\Label{t6.1}
Let $\Om$, $U$, $\chi^t_j,$ $j=0,1,2$ be as above. Assume that $(f\T-1)^k$ holds in a neighborhood $U$ with $f\gg \log$. Then, for any positive integer $s$,  we have
\begin{eqnarray}
\Label{6.1}
\begin{split}
\no{\chi_0^t u}_s^2 \lesssim& t^{-2s}\no{\chi^t_1\Box u}_s^2+(g^*(t^{-1}) )^{2(s+1)}\no{\chi^t_2 u}^2, 
\end{split}
\end{eqnarray}
for any $u\in \A^{0,k}\cap\T{Dom}(\Box)$,  where $g=\frac{f}{\log }$ and $g^*$ is the inverse function to $g$.  
\end{theorem}
\begin{remark}
In \cite{C83}, Catlin proves the same statement for the particular choice $f(|\xi|)=|\xi|^\epsilon$ ending up with a $g$ in \eqref{6.1} which is $f$ itself. In fact, starting from subelliptic estimates, \eqref{6.1} is obtained by induction over $j$ such that $j\epsilon\ge s$. For us, who use Kohn method of \cite{K02}, a loss of sharpness, that is,  $g=\frac f\log$ instead of $f$, is unavoidable.
\end{remark}
{\it Proof of Theorem~\ref{t6.1}. } For each integer $s\ge 0$,  we interpolate two sequences of cutoff functions $\{\zeta_m\}_{m=0}^s$ and $\{\sigma_m\}_{m=1}^s$ with support on $U$ and such that 
$\zeta_j\equiv1$ on $\T{supp}\,\zeta_{j+1}$ and $\sigma_j\equiv1$ on $\T{supp}\,\zeta_j$. 
 We define two new sequences $\{\zeta^t_m\}$ and $\{\sigma^t_m\}$ by $\zeta_m^t(z)=\zeta_m(\frac{z}{t})$ and $\sigma_m^t(z)=\sigma_m(\frac{z}{t})$.
 In our conclusion, we take $\chi_1$ as $\zeta_0$ and $\chi_0$ as $\zeta_s$. 
 \\

We  also need a  pseudodifferential partition of the unity. Let $\lambda_1(|\xi|)$ and $\lambda_2(|\xi|)$ be real valued $C^\infty$ functions such that $\lambda_1+\lambda_2\equiv 1$ and 
$$\lambda_1(|\xi|)=\begin{cases}1 & \T{ if } |\xi|\le 1\\
0 & \T{ if } |\xi|\ge  2. \end{cases}$$
Recall that $\La^m$ is the tangential pseudodifferential operator of order $m$. Denote by $\La_t^m$ the pseudodifferential operator with symbol $\lambda_2(t|\xi|)(1+|\xi|^2)^{\frac{m}{2}}$  and by $E_t$ the operator with symbol $\lambda_1(t|\xi|)$. Note that 
\begin{eqnarray}\Label{6.2b}
\no{ \La^m \zeta^t_m u}^2\lesssim \no{\La_t^m \zeta^t_m u}^2+t^{-2m}\no{\zeta_m^tu}^2.
\end{eqnarray}
In this estimate, it is understood that  $t^{-1}\le g^*(t^{-1})$. 
\\

For $m=1,2,\dots$, we define the pseudodifferential operator $R^m_t$ by 
$$R^m_t\varphi(x,r)=(2\pi)^{-2(n-1)}\int_{\R^{2n-1}}e^{ix\cdot \xi} \lambda_2(t|\xi|)(1+|\xi|^2)^{\frac{m\sigma^t_m(x,r)}{2}}\tilde{\varphi}(\xi,r)d\xi$$
for $\varphi\in C^\infty_c(U\cap \bar\Om)$. Since the symbol of $(\La_t^m-R_t^m)\zeta_{m}^t$ is of order zero, then
\begin{eqnarray}\label{6.3b}
\begin{split}
||\La_t^m\zeta^t_{m}u||^2\lesssim &\no{R_t^m\zeta^t_{m}u}^2+\no{\zeta^t_{m}u}^2\\
\lesssim &\no{\zeta^t_mR_t^m\zeta^t_{m-1}u}^2+\no{[R_t^m, \zeta^t_{m}]\zeta^t_{m-1}u}^2+\no{\zeta^t_{m}u}^2\\
\lesssim &\no{f(\La) \zeta^t_{m-1}R_t^m\zeta^t_{m-1}u}^2+\no{[R_t^m, \zeta^t_{m}]\zeta^t_{m-1}u}^2+\no{\zeta^t_{m}u}^2
\end{split}
\end{eqnarray}
By Proposition \ref{p6.2} which follows, the commutator  in the last line of  \eqref{6.3b} is dominated by $\sum_{j=1}^m t^{-2j}||| \zeta^t_{m-j}u |||_{m-j}^2$.  From \eqref{6.2b} and \eqref{6.3b},  we get
\begin{eqnarray}\Label{6.4b}
\begin{split}
|||\zeta^t_{m}u|||^2_m \lesssim \no{f(\La)\zeta^t_{m-1}R_t^m\zeta^t_{m-1}u}^2+\sum_{j=1}^m t^{-2j}||| \zeta^t_{m-j}u |||_{m-j}^2.
\end{split}
\end{eqnarray}
Similarly, 
\begin{eqnarray}\Label{6.5b}
\begin{split}
|||D_r\La^{-1}\zeta^t_{m}u|||^2_{m} \lesssim &\no{D_r\La^{-1} f(\La) \zeta^t_{m-1}R_t^m\zeta^t_{m-1}u}^2\\
&+\sum_{j=1}^m t^{-2j}|||D_r\La^{-1} \zeta^t_{m-j}u |||_{m-j}^2.
\end{split}
\end{eqnarray}
Denote $A_t^m=\zeta^t_{m-1}R_t^m\zeta^t_{m-1}$ and remark that $A^m_t$ is self-adjoint; also, we have $A_t^mu\in C_c^\infty(U\cap \bar\Om)^k\cap \T{Dom}(\dib^*)$ if $u\in \A^{0,k}\cap\T{Dom}(\dib^*)$. Using the hypothesis and Lemma \ref{l3.38}, we obtain
\begin{eqnarray}\Label{6.6b}
\begin{split}
\no{f(\La) A^m_tu}^2+ \no{D_r\La^{-1} f(\La) A^m_tu}^2 \lesssim Q(A_t^mu, A_t^m u)
\end{split}
\end{eqnarray}
\\

Next, we estimate $Q(A_t^mu,A_t^mu)$. We  have 
\begin{eqnarray}\Label{6.7b}
\no{\dib A_t^m u}^2=(A_t^m \dib u, \dib A_t^m u) +([\dib, A_t^m]u, \dib A_t^m u).
\end{eqnarray}
and 
\begin{multline}
\Label{6.7c}
(A_t^m \dib u, \dib A_t^m u)
=(A_t^m \dib^* \dib u,  A_t^m u) +([\dib, A_t^m]^*u, \dib^* A_t^m u)\\
+(f(\La)^{-1}[[A_t^m, \dib^*]^*, \dib]u, f(\La) A_t^mu).
\end{multline}
Combination of \eqref{6.7b} and \eqref{6.7c} yields
\begin{multline}
\Label{6.7d}
\no{\dib A_t^m u}^2=(A_t^m \dib^* \dib u,  A_t^m u) +([\dib, A_t^m]^*u, \dib^* A_t^m u)\\
+(f(\La)^{-1}[[A_t^m, \dib^*]^*, \dib]u, f(\La) A_t^mu)+([\dib, A_t^m]u, \dib A_t^m u).
\end{multline}
Similarly, 
\begin{eqnarray}\Label{6.8b}
\begin{split}
\no{\dib^* A_t^m u}^2=&(A_t^m \dib \dib^* u,  A_t^m u) +([\dib^*, A_t^m]^*u, \dib A_t^m u) \\
&+(f(\La)^{-1}[[A_t^m, \dib]^*, \dib^*]u, f(\La) A_t^mu)+ ([\dib^*, A_t^m]u, \dib^* A_t^m u)
\end{split}
\end{eqnarray}
Taking summation  of \eqref{6.7d} and \eqref{6.8b}, and using the ``small constant - large constant" inequality, we obtain
\begin{eqnarray}\Label{6.9b}
\begin{split}
Q(A_t^mu,A_t^mu)\lesssim &\no{ A_t^m \Box u}^2+error \\
\lesssim &|||\zeta^t_{m-1} \Box u|||^2_m+error ,\\
\end{split}
\end{eqnarray}
where 
\begin{eqnarray}\Label{6.10b}
\begin{split}
error= &\no{[\dib, A_t^m]u}^2+\no{[\dib^*, A_t^m]u}^2+ \no{[\dib, A_t^m]^*u}+\no{[\dib^*, A_t^m]^*u}\\
&+\no{f(\La)^{-1}[A_t^m, \dib]^*, \dib^*]u}^2+\no{f(\La)^{-1}[A_t^m, \dib^*]^*, \dib]u}^2+\no{A_t^mu}^2.
\end{split}
\end{eqnarray}

 Using Proposition \ref{p6.2} (see below), the error terms are dominated by 
$$\epsilon Q(A_t^mu, A_t^m u)+C_\epsilon (g^{*}(t^{-1}) )^{2(m+1)}\no{\chi_2^t u}^2+\sum_{j=1}^m t^{-2j}\no{ \zeta^t_{m-j}u }_{m-j}^2. $$

Therefore
\begin{eqnarray}\Label{6.11b}
\begin{split}
Q(A_t^mu,A_t^mu)\lesssim |||\zeta^t_{m-1} \Box u|||^2_m+\sum_{j=1}^m t^{-2j}\no{ \zeta^t_{m-j}u }_{m-j}^2+ (g^{*}(t^{-1}) )^{2(m+1)}\no{\chi^t_2 u}^2. 
\end{split}
\end{eqnarray}

Combining \eqref{6.4b}, \eqref{6.5b}, \eqref{6.6b}  and \eqref{6.11b}, we obtain
\begin{eqnarray}\Label{6.12b}
\begin{split}
|||\zeta^t_{m}u|||^2_m+|||D_r\La^{-1}\zeta^t_{m}u|||^2_{m}\lesssim& |||\zeta^t_{m-1} \Box u|||^2_m\\
&+\sum_{j=1}^m t^{-2j}|| \zeta^t_{m-j}u ||_{m-j}^2+(g^{*}(t^{-1}) )^{2(m+1)}\no{\chi_2^tu}^2. 
\end{split}
\end{eqnarray}

Since the operator $\Box$ is elliptic, and therefore non-characteristic with respect to the boundary, we have for $m\ge2$
\begin{eqnarray}\Label{6.13b}
\no{\zeta^t_{m}u}_m^2\lesssim \no{\Box\zeta^t_{m}  u}_{m-2}^2+ |||\zeta^t_{m}u|||^2_m+|||D_r\La^{-1}\zeta^t_{m}u|||^2_{m}. 
\end{eqnarray}
Replace the first term in the right of \eqref{6.13b} by $\no{\zeta^t_m\Box u}^2_{m-2}+\no{[\Box,\zeta^t_m]u}^2_{m-2}$ and observe that the commutator is estimated by $t^{-2}\no{\zeta^t_{m-1}u}^2_{m-1}+t^{-4}\no{\zeta^t_{m-1}u}^2_{m-2}$. Using also \eqref{6.12b}, we get
\begin{eqnarray}\Label{6.14b}
\begin{split}
\no{\zeta^t_{m}u}_m^2\lesssim & \no{\zeta^t_{m-1} \Box u}^2_m+\sum_{j=1}^m t^{-2j}|| \zeta^t_{m-j}u ||_{m-j}^2+(g^{*}(t^{-1}) )^{2(m+1)}\no{\chi^t_2u}^2
\end{split}
\end{eqnarray}
for $m=2,...,s$ (the cases $m=1$ being trivial).

Iterated use of  \eqref{6.14b} to estimate the terms of type $\zeta^t_{m-j}u $ by those of type  $\zeta^t_{m-1} \Box u$ in the right side yields
\begin{eqnarray}
\begin{split}
\no{\zeta^t_{s}u}_s^2\lesssim &\sum_{m=0}^s t^{-2m}\no{ \zeta^t_{s-m}\Box u}_{s-m}^2+(g^{*}(t^{-1}) )^{2(s+1)}\no{\chi^t_2u}^2\\
\lesssim &t^{-2s}\no{ \zeta^t_{0}\Box u}_s^2+(g^{*}(t^{-1}) )^{2(s+1)}\no{\chi^t_2u}^2.
\end{split}
\end{eqnarray}
Choose $\chi^t_0=\zeta^t_s$ and $\chi^t_1=\zeta^t_0$; we then conclude
\begin{eqnarray}
\begin{split}
\no{\chi^t_0u}_s^2\lesssim &t^{-2s}\no{ \chi^t_1\Box u}_s^2+(g^{*}(t^{-1}) )^{2(s+1)}\no{\chi^t_2u}^2,
\end{split}
\end{eqnarray}
for any $u\in \A^{0,k}\cap \T{Dom}(\Box)$.

$\hfill\Box$

\begin{proposition}\Label{p6.2} 
\begin{enumerate}
\item[(i)] $\no{[R_t^m, \zeta^t_{m}]\zeta^t_{m-1}u}^2\lesssim \sum_{j=1}^m t^{-2j}||| \zeta^t_{m-j}u |||_{m-j}^2$
\item[(ii)] Assume that $(f\T-1)^k$ holds with $f\gg \log$, then for any $\epsilon\ge 0$, there is a constant $C_\epsilon$ such that the error term in \eqref{6.9b} is dominated by   $$\epsilon Q(A_t^mu, A_t^m u)+C_\epsilon (g^{*}(t^{-1}) )^{2(m+1)}\no{\chi^t_2 u}^2+\sum_{j=1}^m t^{-2j}\no{ \zeta^t_{m-j}u }_{m-j}^2,$$
 \end{enumerate}
where $g:=\frac{f}{\log }$ and $g^*$ is the inverse function to $g$.
\end{proposition}

{\it Proof.  (i).}
Using formula \eqref{2.33}, the proof of (i) is straightforward.

{\it (ii). }  First, we show 
\begin{eqnarray}\Label{6.19d}
\no{[\dib, A_t^m]u}\le \epsilon Q(A_t^mu, A_t^m u)+C_\epsilon (g^{*}(t^{-1}) )^{2(m+1)}\no{\chi_2^tu}^2+\sum_{j=1}^m t^{-2j}\no{ \zeta^t_{m-j}u }_{m-j}^2.
\end{eqnarray}
By the Jacobi identity,  
\begin{eqnarray}\Label{6.20c}
\begin{split}
[\bar{\partial},A^m_t]=&[\dib, \zeta^t_{m-1}R^m_t\zeta^t_{m-1}]\\
=&   [\bar{\partial},\zeta^t_{m-1}]R^m_t\zeta^t_{m-1}+\zeta^t_{m-1}[\bar{\partial},R^m_t]\zeta^t_{m-1}+\zeta^t_{m-1}R^m_t[\bar{\partial}, \zeta^t_{m-1}].
\end{split}
\end{eqnarray}
Since the support of the derivative of $\zeta_{m-1}^t$ is disjoint from the
support of $\sigma_m^t$, the  first and third terms in the second line of \eqref{6.20c} are bounded by $\frac{C}{t}$ in $L_2$. The middle term in \eqref{6.20c} is treated  as follows. Let  $a$ belong to $\S$ and $D$ be $\frac{\di}{\di x_j}$ or $D_r$; we have 
\begin{eqnarray}\Label{6.21d}
[aD, R^m_t]=[a, R^m_t]D+a[D,R^m_t].
\end{eqnarray}
Now, if $D=\frac{\di}{\di x_j}$, the term first term of  \eqref{6.21d}  is bounded  by $C R^m_t$; if, instead, $D=D_r$, we decompose $D_r=\bar L_n +Tan$, so that $[a, R^m_t]D$  is bounded by $C\La^{-1}\bar L_n R^m_t +C R^m_t $. As for the second term,  we have $[D, R^m_t ]=m D(\sigma_m^t)\log( \La) R^m_t$; in particular, $[D, R^m_t]$ is bounded by $\frac{C}{t}\log(\La) R^m_t$. Therefore, 
\begin{eqnarray}\Label{6.22d}
\begin{split}
\no{[\dib, A_t^m]u}^2
\lesssim & \frac{1}{t}\no{\log(\La) A_t^m u}^2+\epsilon \no{\bar L_nA^m_t u}^2+C_\epsilon \no{\chi_2^t u}^2+ \sum_{j=1}^m t^{-2j}\no{ \zeta^t_{m-j}u }_{m-j}^2.
\end{split}
\end{eqnarray}
To estimate the first term in \eqref{6.22d}, we  check that 
$$ \frac{1}{t}\log \La_\xi \le \epsilon f(\La_\xi)~~~~\T{on support of }~ \lambda_2\Big (g^*\big((\epsilon t)^{-1}\big)^{-1}\La_\xi\Big) $$
and hence
\begin{eqnarray}
\frac{1}{t}\log\La_\xi\lesssim \epsilon f(\La_\xi)+t^{-1}\lambda_1\Big (g^*\big((\epsilon t)^{-1}\big)^{-1}\La_\xi\Big)\log \La_\xi.
\end{eqnarray}
It follows 
\begin{eqnarray}\begin{split}
\frac{1}{t}\no{\log\La A^m_t u}^2 \le& \epsilon \no{f(\La)A^m_t u}^2+t^{-2}g^*\big((\epsilon t)^{-1}\big)^{2m}\log^2\Big( g^*\big((\epsilon t)^{-1}\big)\Big )\no{\chi_2^t u}^2\\
\le& \epsilon \no{f(\La)A^m_t u}^2+C_\epsilon g^*(t^{-1})^{2(m+1)}\no{\chi_2^t u}^2.
\end{split}
\end{eqnarray}

Since we are supposing that $(f\T-1)^k$ holds, we get the proof of the inequality \eqref{6.19d}. 
By a similar argument, we can estimate all error terms in \eqref{6.9b} and obtain the conclusion of the proof of Theorem~\ref{t6.1}.

$\hfill\Box$

\section{Kohn's proof of  local regularity}
In \cite{K02}, Kohn proved that superlogarithmic estimates of the system $(\dib_b, \dib_b^*)$ on a hypersurface $M$ imply local regularity of the operator  $\Box_b^{-1}$; and superlogarithmic estimates of the system $(\dib_b, \dib_b^*)$ 
on microlocalized forms $u^+$
imply local regularity of the $\dib$-Neumann operator on the side $\Om^+$. The purpose of this section is to give a direct proof, without passing through the tangential system, of the local regularity of the $\dib$-Neumann operator and the Bergman projection and of the smoothness of the Bergman kernel. Also, our argument applies to general $q$-pseudoconvex domains. 

\begin{theorem}\Label{t6.4g}
Let $\Om\subset\C^n$ be a bounded $q$-pseudoconvex domain with $C^\infty$ boundary.
Suppose that $\Box$ is invertible, modulo $\mathcal H^{0,k}$ for $k\ge q$, and denote by $N_k$ its inverse over $k$-forms. 
 Suppose further that superlogarithmic estimates hold in degree $k\ge q$, that is, $(f\T-\M)^k$ estimates hold for $f=\log$ and   $\M=\frac{1}{\epsilon}$ for any $\epsilon>0$. Let $\chi_0$ and $\chi_1$ be smooth cutoff functions supported in $U$ with $\chi_1\equiv1$ in a neighborhood of the support of $\chi_0$. Then the $\dib$-Neumann operator $N_k$ as well as  the Bergman projection $P_{k-1}$ satisfy for any $s$ and for any $\alpha\in H_s(\Om)^{0,k}$ the estimates
\begin{eqnarray}
\no{\chi_0 N_k \alpha}_s^2&\lesssim &\no{ \chi_1\alpha }_s^2+\no{\alpha}^2,\\
\no{\chi_0 \dib^* N_k \alpha}_s^2+\no{\chi_0 \dib N_k \alpha}_s^2&\lesssim &\no{ \chi_1\alpha }_s^2+\no{\alpha}^2,\\
\no{\chi_0 P_{k-1} \alpha}_s^2&\lesssim &\no{ \chi_1\alpha }_s^2+\no{\alpha}^2.
\end{eqnarray}
Moreover, if $(z_0,w_0)$ is a point in $(b\Om\times b\Om)\setminus\T{\rm (Diagonal)}$ such that $(f\T-\M)^k$ estimates hold at both points $z_0$ and $w_0$, then the Bergman kernel $K(z,w)$ extends smoothly to $\bar\Om\times\bar\Om$ locally in a neighborhood of $(z_0,w_0)$.

\end{theorem}

{\it Proof. }
The first statement, about $N_k$, $\dib^* N_k$, $\dib N_k$ and $P_{k-1}$ follows from Theorem~\ref{t6.1} combined with the standard method of elliptic regularization by Kohn-Nirenberg \cite{KN65}. As for the last statement concerning $K(z,w)$, it follows from the previous one combined with the results of Kerzmann \cite{Ke72}.

$\hfill\Box$

 \section{Catlin's proof of the necessary condition}
In this section we derive geometric conditions for $b\Om$ if the  $(f\T-1)^k$ estimate  holds. We restate and prove now the first half part of Theorem~\ref{main3}.

\begin{theorem}
\Label{t6.3b}
Let $\Om$ be a pseudoconvex domain  with a local defining function $r$ at $z_0$  and suppose further that the $(f\T-1)^k$ estimate holds. Let $Z$ be a $k$-dimensional complex analytic variety passing through $z_0$ and suppose that $b\Om$ has type $\le F$  along $Z$ at $z_0$ in the sense that \index{Type ! $\le F$}
\begin{equation}
\Label{6.100}
|r(z)|\lesssim F(|z-z_0|)\quad\T{for $z\in Z$ close to $z_0$.}
\end{equation}
 Then,  for small $\delta$, $$\frac{f(\delta^{-1})}{\log(\delta^{-1})} \lesssim F^*(\delta )^{-1},$$  where $F^*$ is the inverse function to $F$.
\end{theorem}
To prove Theorem~\ref{t6.3b} we have to recall the following theorem from \cite{C83} which holds for any domain $\Om\subset\C^n$ not necessarily $q$-pseudoconvex nor $q$-pseudoconcave. The result of Catlin is proved for domains of {\it finite type}, that is, for $F(\delta)=\delta^m$, but it holds in full generality of $F$.

\begin{theorem}
\Label{t6.4b}
Let $\Om$ be a domain in $\C^n$ with smooth boundary and assume that that there is a function $F$ and a $k$-dimensional complex-analytic variety $Z$ passing through $z_0$ such that \eqref{6.100} is satisfied  for $z\in Z$, $z$ close to $z_0$. Then, in any neighborhood $U$ of $z_0$, there is a family $\{Z_t, t\in T\}$ of $k$-dimensional complex manifolds contained in $U$ of diameter comparable to $t$ such that 
$$\sup_{z\in Z_t}|r(z)|\lesssim F(t).$$
\end{theorem}

For the proof of Theorem~\ref{t6.3b} we also  need the next lemma which is contained in \cite{C83}.

\begin{lemma}
\Label{l6.5b}
Let $\Om$ be a bounded pseudoconvex domain in $\C^n$.  Let $z_0\in b\Om$ and assume that $b\Om$ is smooth near $z_0$ with a smooth boundary-defining function $r$ normalized by $\frac{\di r}{\di x_n}(z_0)>0$.  Then there exists a neighborhood $U$ of $z_0$ such that for each $w\in U\cap\Om$, there is a function $G$ holomorphic and $L_2$ on $\Om$  with the following properties   
\begin{enumerate}
  \item \quad$\no{G}^2\lesssim 1$
  \item \quad$\Big |\frac{\di^k G}{\di z_n}(w)\Big|\simge |r(w)|^{-(k+\frac{1}{2})}$ for all $k\ge 0$.
 \end{enumerate} 
\end{lemma}

{\it Proof of Theorem~\ref{t6.3b}.}  By Theorem~\ref{t6.4b}, for each $t\in T$, there is a point $\gamma_t\in Z_t$, which satisfies $|r(\gamma_t)|\lesssim F(t)$ and by Lemma~\ref{l6.5b} there is a function $G_t\in A(\Om)\cap L_2(\Om)$ such that  
$$\no{G_t}\le 1$$
and 
$$\Big|\frac{\di^m G_t}{\di z_n^m}(\gamma_t)\Big |\simge F(t)^{-(m+\frac{1}{2})}.$$
We parametrize $Z_t$ over $\C^k\times\{0\}$ by
$$
z'\mapsto (z',h_t(z')).
$$

Let $\phi $ be a cutoff function on $\R^+$ such that $\phi=1$ on $[0,1)$ and $\phi=0$ on $[2,+\infty)$. Let, for some $a$ to be chosen later 
$$\psi_t(z')=\phi\Big(\frac{8|z'-\gamma_t'|}{at}\Big),$$
where $z'=(z_1,...,z_k)$.
Choose the datum $\alpha_t$ as  
$$\alpha_t=\psi_t(z')G_t(z)d\bar z_1\we...\we d\bar z_k$$
Clearly the form $\alpha_t$ is $\dib$-closed and its coefficient belongs to $L_2$.  Let $v_t$ be the canonical solution of $\dib$ with datum $\alpha_t$, that is,  $v_t=\dib^*u_t$ where $\Box u_t=\alpha_t$.
We have 
$$\dib v_t=\alpha_t.$$
Define the action of $\frac{\di^m}{\di z^m_n}$ over a form  by applying  $\frac{\di^m}{\di z^m_n}$ to each coefficient. Since $\frac{\di^m}{\di z^m_n}$ and $\dib$ commute, we then have
\begin{equation}
\Label{com*}
\frac{\di^m}{\di z^m_n}\alpha_t=\frac{\di^m}{\di z^m_n}\dib v_t=\dib \frac{\di^m}{\di z^m_n}v_t.
\end{equation}
 \\

On each polydisc $B_{at}^k(\gamma'_t)$ of center $\gamma'_t$ and radius $at$, we put the standard Euclidean metric. This gives rise to an inner product on differential forms supported in $B^k_{at}(\gamma_t')$ that we denote by 
\begin{equation}
\Label{inner*}
\int_{B^k_{at}(\gamma'_t)}\la,\ra dV,
\end{equation}
where $dV$ is the volume element in $\C^k$.\\
We also define a $k$-form $w_t$ by
$$w_t=\phi(\frac{8|z'-\gamma'_t|}{3at})d\bar z_1\we...\we d\bar z_k.$$
 We use the notation
\begin{equation}
\Label{new}
\mathcal K^m_t=\int_{B_{at}^k(\gamma'_t)}\la \frac{\di^m }{\di z_n^m}\alpha_t(z',h_t(z')),w_t \ra dV.
\end{equation}
We use  polar coordinates to evaluate \eqref{new} and get
\begin{equation}
\Label{newbis}
\mathcal K_t^m=2^k\int_{0}^{ta/4}\phi(\frac{8s}{ta})\phi(\frac{8s}{3ta})\int_{|w'|=1} \frac{\di^m }{\di z_n^m}\psi_tG_t(w',h_t(\gamma'_t+sw'))s^{2k-1}dSds.
\end{equation}
Since $\frac{\di^m }{\di z_n^m}G_t(z',h_t(z'))$ is holomorphic, it satisfies the Mean Value Equality. Thus

$$\mathcal K_t^m=2^k\Big(\int_{0}^{ta/4}\chi(\frac{8s}{ta})s^{2k-1}ds\Big)\frac{\di^m }{\di z_n^m}\psi_tG_t(w',h_t(\gamma'_t+sw'))\int_{|w'|=1}dS.$$

Recalling that $\psi_tG_t(w',h_t(\gamma'_t+sw'))=G_t(\gamma_t)$, then, on account of \eqref{new}, we have a lower bound for $\mathcal K_t^m$ 
\begin{equation}
\Label{lower*}\mathcal K_t^m\simge t^{2k}F(t)^{-(m+\frac{1}{2})}.
\end{equation}

On the other hand, from \eqref{com*} we obtain
$$\mathcal K_t^m=\int_{B_{at}^k(\gamma'_t)}\la \dib \frac{\di^m }{\di z_n^m} v_t(h_t),w_t \ra dV=\int_{B_{at}^k(\gamma'_t)}\la  \frac{\di^m }{\di z_n^m} v_t(h_t),\vartheta w_t \ra dV$$
where $\vartheta$ is the adjoint of $\dib$ with respect to the inner product \eqref{inner*}. We define a set $$S_t=\{z'\in \C^k : \frac{3at}{8}\le |z'-\gamma_t'|\lesssim \frac{6at}{8} \}.$$ Since  $\vartheta w_t$ is supported in $S_t$ and  $|\vartheta w_t|\lesssim t^{-1}$, then 
\begin{equation}
\Label{upper*}
\mathcal K_t^m\lesssim t^{2k-1}\sup_{S_t} | \frac{\di^m }{\di z_n^m} v_t(h_t) |\lesssim t^{2k-1}\sup_{h_t(S_t)} | \frac{\di^m }{\di z_n^m} \dib^* u_t |\lesssim t^{2k-1}\sup_{h_t(S_t)} | \underset{|\beta|=m+1}{D^\beta} u_t |.
\end{equation}
Before completing the proof of Theorem~\ref{t6.3b}, we need
the proposition which follows which is dedicated to state an upper bound for $\mathcal K_t^m$. 
\begin{proposition}\Label{p6.6}
 \begin{equation}
 \Label{upperbis*}
 \sup_{h_t(S_t)} | \underset{|\beta|=m+1}{D^\beta} u_t |\lesssim g^*(t^{-1})^{m+n+3}.
 \end{equation}
\end{proposition}
{\it Proof.}  Since the set $S_t$ has  diameter $O(at)$ and the function $h_t$ satisfies $|dh_t(z')|\le C$ for $z'\in B^k_{at}(\gamma'_t)$, then the set $h_t(S_t)$ has diameter of size $O(at)$. Moreover, by  construction, there exists a constant $c$ such that 
$$\inf\{|z_1-z_2| : z_1\in \T{supp}\alpha_t, \quad   z_2\in (\T{id}\times h_t)(S_t)\}>2ct.$$
Therefore, we may choose $\chi_0$ and $ \chi_1$ such that if we set $\chi_k^t(z)=\chi_k (\frac{z}{ct})$ for  $k=0,1$, we have the properties
\begin{enumerate}
  \item[(1)]\qquad$\chi^t_0=1 \T{~~~~on~~~~} (\T{id}\times h_t)(S_t)$
\item[(2)] \qquad$\alpha_t=0 \T{~~~~on~~~~} \T{supp}\chi^t_1 .$
\end{enumerate}

Hence
\begin{eqnarray}
\begin{split}
\sup_{(\T{id}\times h_t)(S_t)} | \underset{|\beta|=m+1}{D^\beta} u_t |\lesssim&\sup_{\Om\cap Z_t}| \underset{|\beta|=m+1}{D^\beta} \chi_0^t u_t |\lesssim \no{\chi_0^t u_t}_{m+n+1},
\end{split}
\end{eqnarray}
where the last inequality follows from Sobolev Lemma since $\chi_0^tu_t$ is smooth by Theorem \ref{t6.4g}. 
We use now Theorem~\ref{t6.1} and observe that $\chi_1^t\Box u_t=0$ (by property 2 of $\chi^t_1$). It follows
$$ \no{\chi_0^tu_t}^2_{m+n+1}\lesssim g^*(t^{-1})^{2(m+n+2)}\no{u_t}^2.$$
Since $\Om$ is bounded pseudoconvex, then
$$\no{u_t}^2\lesssim \no{\Box u_t}^2=\no{\alpha_t}^2\lesssim 1.$$
This completes the proof of Proposition~\ref{p6.6}. 
\\

$\hfill\Box$

{\it End of Proof of Theorem~\ref{t6.3b}.}
We return to the proof of Theorem~\ref{t6.3b}. Combining \eqref{lower*} with \eqref{upper*} and \eqref{upperbis*}, we get the upper estimate
$$ t^{2k}F(t)^{-(m+\frac{1}{2})}\le C t^{2k-1} g^*(t^{-1})^{m+n+2}.$$  
It follows $$t^{\frac{1}{m}}F(t)^{-1-\frac{1}{2m}}\le C^{\frac{1}{m}} g^*(t^{-1})^{1+\frac{n+2}{m}}.$$
Since this holds for all positive $m$, we may take the limit as $m\to \infty$ and end up with
$$F(t)^{-1}\le g^*(t^{-1}).$$
This concludes the proof of Theorem~\ref{t6.3b}.

$\hfill\Box$

\section{Mc Neal's proof of lower bounds for the Bergman metric}

Let $\Om$ be a bounded domain in $\C^n$. Denote $H(\Om)$ by the space of holomorphic functions on $\Om$ and set $A^2(\Om)=H(\Om)\cap L^2(\Om)$. The space $A^2(\Om)$ is a Hilbert space under the inner product 
$$(\varphi,\psi)=\int_{\Om} \varphi\bar \psi dV.$$

For each point $a\in \Om$, the  evaluation map  
$$e_a: A^2(\Om)\ni \varphi \mapsto \varphi(a)\in \C$$
 is a bounded linear functional on $A^2(\Om)$. Therefore,  by the Riesz Representation Theorem, there is  a unique element in  $A^2(\Om) $, that we  denote by $K_\Om(\cdot, a)$, such that 
\begin{eqnarray}\Label{6.18c}
\varphi(a)=e_a(\varphi)=( \varphi ,K(\cdot, a))=\int \varphi(\zeta)\overline{K(\zeta, a)}dV(\zeta)
\end{eqnarray}
for all $\varphi\in A^2(\Om)$. The function 
$$K_\Om:\Om\times\Om\to \C $$
obtained in this way  is called the {\it Bergman kernel}  for  $\Om$. \index{Bergman ! kernel, $K_\Om(z,w)$}

For the choice $\varphi=K_\Om(\cdot, z)$ and $a=\zeta\in \Om$, \eqref{6.18c} yields the proof of the following  fundamental symmetry property of $K$.

\begin{lemma} The Bergman kernel $K$ satisfies
\begin{eqnarray}\Label{6.19c}
K(\zeta, z)=\overline{K(z, \zeta)}~~~~~~\T{~~~for all } \zeta, z\in \Om.
\end{eqnarray}
\end{lemma}

Notice that by \eqref{6.19c}, we can rewrite \eqref{6.18c} in the form
\begin{eqnarray}
\varphi(z)=\int \varphi(\zeta)K(z, \zeta )dV(\zeta)~~~\T{ for all } \varphi \in A^2(\Om) \T{ and } z\in \Om. 
\end{eqnarray}

Since $A^2(\Om)$ is a closed subspace of $L_2(\Om)$, there is an orthogonal projection operator $P:=P_{\Om}$ \index{Bergman ! projection, $P, P_{k-1}$}
\begin{eqnarray}
P: L^2(\Om)\to A^2(\Om).
\end{eqnarray}
This operator is usually called  the {\it Bergman projection}. The Schwartz kernel theorem implies that the Bergman projection is represented by integration against the Bergman kernel function, that is, 
\begin{eqnarray}\Label{6.22c}
P\varphi(z)=\int \varphi(\zeta)K(z, \zeta )dV(\zeta)\T{    for all } \varphi \in L_2(\Om) \T{ and } z\in \Om.
\end{eqnarray}
The relation between the Bergman projection $P$ and the $\dib$-Neumann operator $N$,  is expressed by Kohn's formula 
\begin{eqnarray}\Label{6.23c}
P=I-\dib^*N\dib.
\end{eqnarray}

From the Bergman kernel $K_\Om$,  one obtains a very interesting metric on $\C T(\Om)$.  Let 
\begin{eqnarray}
g_{ij}(z)=\frac{\di^2}{\di z_i\di \bar z_j} \log K(z,z)
\end{eqnarray}
 and,  for $X=\sum^n_{i=j}a_j\frac{\di}{\di z_j}$, define 
\begin{eqnarray}
B_\Om(z, X)=\Big(\sum_{ij=1}^n g_{ij} a_i \bar a_j\Big)^\frac{1}{2}.
\end{eqnarray}\index{Bergman ! metric, $B_\Om{z,X}$}
This differential metric, the Bergman metric, is primarily interesting because of its holomorphic invariance properties. If we derive the set-theoretic metric on $\Om$ from $B_\Om$ in the standard way, then, with respect to this metric, distances are preserved under a biholomorphic transformation on $\Om$.\\

One can obtain the value of the Bergman kernel function on the diagonal of $\Om\times\Om$ and the length of tangent vector $X$ in the Bergman metric by solving  the following extreme value problems : 
\begin{eqnarray}\begin{split}\Label{6.26c}
K_{\Om}(z,z)=&\inf\{ \no{\varphi}^2 :\varphi \in H(\Om), \varphi(z)=1\}^{-1}\\
=&\sup\{|\varphi(z)|^2 : \varphi \in H(\Om), \no{\varphi}\le 1\}
\end{split}
\end{eqnarray}
and 
\begin{eqnarray}\begin{split}\Label{6.27c}
B_{\Om}(z,X)=&\frac{\inf\{ \no{\varphi} :\varphi \in H(\Om), \varphi(z)=0, X\varphi(z)=1\}^{-1}}{\sqrt{K_{\Om}(z,z)}}\\
=&\frac{\sup\{|X\varphi(z)| : \varphi \in H(\Om), \varphi(z)=0, \no{\varphi}\le 1\}}{\sqrt{K_{\Om}(z,z)}}.
\end{split}
\end{eqnarray}

The purpose of this section is to study the boundary behavior of $B_{\Om}(z,X)$ for $z$ near a point $z_0\in b\Om$, when a $(f\T-1)^1$ estimate for the $\dib$-Neumann problem holds. In this section, we show the second part of Theorem \ref{main3}, that is, 

\begin{theorem}\Label{t6.9}Let $\Om\subset\subset \C^n$ be   pseudoconvex with a smooth boundary in a neighborhood $U$ of $z_0$. Suppose that  $(f\T-1)^1$ holds in $U$ with $f\gg \log$. Given a constant $\eta>0$,  there exists $U'\subset  U$,  so that for all $z\in U'\cap \Om$ and any $X\in T^{1,0}_z\C^n$,
$$B_\Om(z,X)\simge g(\delta(z)^{-1+\eta})|X|$$
 where $\delta(z)$ is the distance of $z$ to $ b\Om$ and $g=\frac{f}{\log}$.    
\end{theorem}

Theorem \ref{t6.9} is a new version of Theorem 1 by McNeal \cite{McN92}. In Theorem \ref{t6.9}, we treat  the case in which the domain has non-finite type.  For the proof of Theorem \ref{t6.9},  we first use the classical results about locally comparable properties of Bergman kernel and Bergman metric, that is, 
\begin{proposition}\Label{p6.10}
Let $\Om_1, \Om_2$ be  bounded pseudoconvex domains in $\C^n$ such that a portion of $b\Om_1$ and $b\Om_2$ coincide. Then
\begin{eqnarray}
\Label{44f} K_{\Om_1}(z,z)&\cong& K_{\Om_2}(z,z);\\
\Label{45f}B_{\Om_1}(z,X)&\cong& B_{\Om_2}(z,X),
\end{eqnarray}
for $z$ near the coincidental portion of the two boundaries and $X\in T^{1,0}_z(\C^n)$.
 \end{proposition}
(cf. \cite{McN92} or \cite{DFH84}.)\\  

To apply Proposition \ref{p6.10}, we construct a smooth pseudoconvex domain $\tilde\Om$, contained in $\Om$, that shares a piece  of its boundary with $b\Om$ near $z_0$.  The crucial property that $\tilde\Om$ has, for our purpose, is the exact, global regularity of the $\dib$-Neumann operator. In fact, one can show that

\begin{proposition}\Label{p6.11}
Let $\Om$ be  a smooth, bounded, pseudoconvex domain in $\C^n$ and let $z_0\in b\Om$. Then, there exist a neighborhood $U$ of $z_0$ and a smooth, bounded, pseudoconvex domain $\tilde\Om$ satisfying the following properties:
\begin{enumerate}
  \item[(i)] $\tilde\Om\subset \Om\cap U$,
  \item[(ii)] $b\tilde\Om\cap b\Om$ contains a neighborhood of $z_0$ in $b\Om$,
  \item[(iii)] all points in $b\tilde\Om\setminus b\Om$ are points of strong pseudoconvexity.    
\end{enumerate}
 \end{proposition}
A proof of Proposition \ref{p6.11} can be found in \cite{McN92}. Let $z_0\in b\Om$ such that an $(f\T-1)^1$ estimate holds in a neighborhood of $z_0$. By Proposition  \ref{p6.11} one can get that the $\dib$-Neumann operator $N$ of $\tilde\Om$ is exactly, globally regular and that the  $(f\T-1)^1$ estimate holds in a neighborhood of $z_0$. Therefore, Theorem \ref{t6.9} is an immediate consequence of Proposition \ref{p6.10}, \ref{p6.11} and the following

\begin{theorem}\Label{t6.12}
Let $\Om$ be a smooth, bounded, pseudoconvex domain in $\C^n$ which has an exactly, globally regular $\dib$-Neumann operator. Let $z_0\in b\Om$ and $U$ be a neighborhood of $z_0$. Suppose that  $(f\T-1)^1$ holds in $U$ with $f\gg \log$. Given a constant $\eta>0$,  there exists $U'\subset  U$,  so that for all $z\in U'\cap \Om$ and any $X\in T^{1,0}_z\C^n$,
$$B_\Om(z,X)\simge g(\delta(z)^{-1+\eta})|X|,$$
 where $\delta(z)$ is the distance of $z\in b\Om$ to $g=\frac{f}{\log}$. 
  \end{theorem}
Before giving the proof of Theorem \ref{t6.12}, we need some preliminaries.\\

Let $\psi$ be the cutoff functions such that  
 $$\psi(z) =\begin{cases}0  &\T{~~~if } z\in B(z_0,1),\\ 
1&\T{~~~if } z\in \C^n\setminus B(z_0,2),\end{cases} $$
where $B(z_0, a)$ is the ball in $\C^n$ with center $z_0$ and radius $a$.  We also set  $\psi^t=\psi (\frac{z}{t})$. To prove Theorem \ref{t6.12}, we also need the following

\begin{proposition}\Label{p6.13} Suppose that $\Om$  satisfies the hypothesis of Theorem \ref{t6.12}. Let $\psi^t$ be defined as above. Then if $\alpha\in C^\infty_c(B(z_0, \frac{t}{8})\cap \bar\Om)^1$,  for any nonnegative integer $s_1, s_2$, we have 
\begin{equation}
\Label{newbis}
\no{\psi^tN\alpha }^2_{s_1}\lesssim  g^*(t^{-1})^{2(s_1+s_2+4)} \no{\alpha}^2_{-s_2}.
\end{equation}
\end{proposition}
{\it Proof. } We choose a triple of cutoff functions $\chi^t_0, \chi^t_1 $ and $\chi^t_2$ in Theorem \ref{t6.1}, such that $\chi^t_0\equiv 1$ on a neighborhood of the support of the derivative of $\psi^t$ and supp $\chi_2^t\subset  B(z_0, 3t)\setminus B(z_0, \frac{t}{2})$; hence $\chi_1^t\alpha=0$.   We notice that for $t$ sufficiently small, supp$\chi_j^t\subset \subset U$ for  $j=0,1,2$, so that we can apply Theorem \ref{t6.1}  to this triple of cutoff functions. Using the global regularity estimate and Theorem \ref{t6.1} for an arbitrary 1-form $u\in \A^{0,1}\cap \Dom(\Box)$, we have 
\begin{eqnarray}\Label{6.46}
\begin{split}
\no{\psi^t u}_{s_1}^2\lesssim& \no{\Box \psi^t u}_{s_1}^2\\
\lesssim& \no{\psi^t\Box  u}_{s_1}^2+\no{[\Box, \psi^t]u}_{s_1}^2\\
\lesssim& \no{\psi^t\Box  u}_{s_1}^2+t^{-2}\no{\chi^t_0 u}_{s_1+1}^2+t^{-4}\no{\chi^t_0 u}_{s_1}^2\\
\lesssim&  \no{\psi^t\Box u}_{s_1}^2+t^{-2(s_1+2)}\no{\chi^t_1\Box  u}_{s_1+1}^2 +g^*(t^{-1})^{2(s_1+3)}\no{\chi_2^t u}^2.
\end{split}
\end{eqnarray}

Recall that we are supposing that the $\dib$-Neumann operator is globally regular.  If $\alpha \in C^\infty(\bar\Om)^1$, then $N\alpha\in C^\infty(\bar \Om)^1\cap \T{Dom}(\Box)$. Substituting $u=N\alpha$ in \eqref{6.46} for $\alpha\in C^\infty_c(B(z_0, \frac{t}{8})\cap \bar\Om)^1$, we obtain 
\begin{eqnarray}\Label{6.47}
\no{\psi^t N\alpha }_{s_1}^2\lesssim g^*(t^{-1})^{2(s_1+3)}\no{\chi^t_2 N\alpha}^2.
\end{eqnarray}

However, 
$$\no{\chi^t_2 N\alpha}=\sup\{|(\chi_2^t N\alpha, \beta)| : \no{\beta}\le 1\},$$
and the self-adjointness of $N$ and the Cauchy-Schwartz inequality yield  

\begin{eqnarray}
\begin{split}
|(\chi_2^t N\alpha, \beta)|=&|(\alpha, N\chi_2^t \beta)|\\
=&|(\alpha, \tilde \chi^t_0 N\chi_2^t \beta )|\\
\lesssim& \no{\alpha}_{-s_2}\no{\tilde\chi^t_0 N\chi_2^t \beta }_{s_2},
\end{split}
\end{eqnarray}
where $\tilde \chi^t_0$  is a cutoff function such that  $\tilde\chi^t_0\equiv 1$ on supp $\alpha$. Let $\tilde \chi^t_1$ be a cutoff function such that $\tilde\chi^t_1=1$ on supp $\tilde \chi^t_0$ and $\T{supp~} \tilde \chi^t_1 \subset\subset  B(z_0, \frac{t}{4})$. Hence  $\T{supp~}\tilde \chi^t_1 \cap\T{supp~}\chi^t_2=\emptyset$. Let $\tilde\chi^t_2$ be a cutoff function such that $\tilde\chi^t_2=1$ on supp $\tilde\chi^t_1$.   Using again Theorem \ref{t6.1} for the triple of cutoff functions $\tilde \chi^t_0,\tilde \chi^t_1$ and $\chi_2^t$, we obtain 
\begin{eqnarray}
\begin{split}
\no{\tilde\chi^t_0  N\chi_2^t \beta}_{s_2}^2\lesssim& t^{-2s_2}\no{\tilde\chi^t_1 \chi_2^t \beta }_{s_2}^2+g^*(t^{-1})^{2(s_2+1)}\no{\tilde\chi^t_2 N \chi_2^t \beta}^2\\
\lesssim&g^*(t^{-1})^{2(s_2+1)}\no{\tilde\chi^t_2 N \chi_2^t \beta}^2\\
\lesssim&g^*(t^{-1})^{2(s_2+1)}\no{\beta}^2.
\end{split}
\end{eqnarray}
Taking supremum over all $\no{\beta}\le 1$, we get \eqref{newbis}.

$\hfill\Box$

{\it Proof of Theorem \ref{t6.12}.} Let $\eta=(\eta_1,...,\eta_n)$ the complex coordinates near $z_0$. $L_2$ norms and derivatives are computed with respect to the $\eta$ variables. For $z\in U'$, we can choose $\eta$ such that  $X(z)=\frac{\di}{\di \eta_1}$. We may also suppose that $\frac{\di r}{\di x_n}(z)>0$ where $x_n+iy_n=\eta_n$.  If $z\in U'$ and $z\not\in b\Om$, we define
$$h_z(\eta)=\frac{K_{\Om}(\eta, z)}{\sqrt{K_\Om(z,z)}};$$
hence $\no{h_z}=1$ and $\frac{|h_z(z)|}{\sqrt{K_{\Om}(z,z)}}=1$. 
We set $$\gamma_z(\eta)=R(z)(\eta_1-z_1)h_z(\eta)$$
where $R(z)$ will be  chosen later. It is an obvious fact that $\gamma_z \in H(\Om)$ and $\gamma_z(z)=0$. We  show that for a specified  choice of $R(z)$, the norm $\no{\gamma_z}\le 1$. Then by \eqref{6.27c}, we get 
\begin{eqnarray}\Label{6.50}
B_\Om(z,X)\ge \frac{ |X\gamma_z(z)|}{\sqrt{K_{\Om}(z,z)}}=\frac{|R(z)h_z(z)|}{\sqrt{K_{\Om}(z,z)}}=|R(z)|.
\end{eqnarray}

In all what follows,  $z$ is  fixed in $U'$.  Let $t=\theta(\delta(z))$, for a function $\theta$ to be chosen later.  For $\psi^t$, as in Proposition \ref{p6.13}, we define $\psi^t_z(\eta)=\psi^t(\eta-z)$.  We write
 \begin{eqnarray}\Label{6.51f}
\gamma_z(\eta)=\psi_z^t(\eta)\gamma_z(\eta)+(1-\psi_z^t(\eta))\gamma_z(\eta).
\end{eqnarray}

A simple estimation of the second term gives
 \begin{eqnarray}\Label{6.52f}
\no{(1-\psi_z^t)\gamma_z}\lesssim |R(z)|t = |R(z)|\theta(\delta(z)).
\end{eqnarray}
 
In order to estimate the first term in right hand side of \eqref{6.51f}, we rewrite this term  as
 \begin{eqnarray}
\psi^t_z(\eta)\gamma_z(\eta)=\psi^t_z (\eta)R(z)(\eta_1-z_1)\frac{D^m\bar G(z)}{D^m\bar G(z)} h_z(\eta),
\end{eqnarray}
where $G$ is  the holomorphic function defined in Lemma~\ref{l6.5b}, and $D^m$ denote the operator $\frac{\di^m}{\di \bar\eta_n^m}$ for $m\in \4N$. 
Notice that by \eqref{6.26c} and Lemma~\ref{l6.5b} for $m=0$, we get 
\begin{equation}
\Label{newter}
K_\Om(z,z)\simge |G(z)|^2\simge \delta^{-1}(z).
\end{equation}
Therefore,
 \begin{eqnarray}\Label{6.31c}
\begin{split}
\no{\psi_z^t\gamma_z}\le& |R(z)|\cdot\T{diam}\Om \cdot\frac{1}{|D^m\bar G(x)|}\cdot\frac{1}{\sqrt{K_\Om(z,z)}}\cdot\no{\psi_z^tK_\Om(\cdot, z)D^m\bar G(z)}\\
\lesssim& |R(z)| \delta^{m+1}\no{\psi_z^tK_\Om(\cdot, z)D^m\bar G(z)}.
\end{split}
\end{eqnarray}
 Now, we consider the norm in last line of \eqref{6.31c}.  Let $\phi\in C^\infty_c(B(0,1))$ be a nonnegative radial function with $\int\phi=1$. If $z\in U'$, set 
$$\phi^t_z(\zeta)=\Big(\frac{t}{10}\Big)^{-2n}\phi\Big(\frac{\zeta-z}{t/10}\Big).$$
From the mean value theorem for antiholomorphic functions, \eqref{6.22c} and \eqref{6.23c}, we have
 \begin{eqnarray}\begin{split}
K_\Om(\zeta, z)D^m\bar G(z)
=&\int_\Om K(\zeta, w)\big(D^m\bar G(w)\big)\phi^t_z(w)dw\\
=&P\Big(\big( D^m\bar G(\zeta)\big) \phi^t_z(\zeta)\Big)\\
=&\big( D^m\bar G(\zeta)\big) \phi^t_z(\zeta)-\dib^* N\dib\Big (\big( D^m\bar G(\zeta)\big) \phi^t_z(\zeta)\Big).
\end{split}
\end{eqnarray}

Notice that the supports of $\psi^t_0$ and $\phi^t_z$ are disjoint,  and that supp $\dib\Big (\big( D^m\bar G(\zeta)\big) \phi^t_z(\zeta)\Big)$ is contained in $B(z, \frac{t}{8})$ for all $z\in U'$. We may therefore apply Proposition \ref{p6.13} for $z_0$ replaced by $z$, with $s_1=1$, and obtain

 \begin{eqnarray}\begin{split}
\no{\psi^t_z K_\Om(\cdot , z)D^m\bar G(z)}^2=&\no{\psi^t_z \dib^* N\dib\Big (\big( D^m\bar G\big) \phi_z^t\Big)}^2\\
\lesssim&\no{\psi^t_z N\dib\Big (\big( D^m\bar G\big) \phi^t_z\Big)}_1^2+\no{[\psi_z^t, \dib^*] N\dib\Big (\big( D^m\bar G\big) \phi^t_z\Big)}^2\\
\lesssim & g^*(t^{-1})^{2(1+s_2)}\no{ \dib\Big (\big( D^m\bar G\big) \phi^t_z\Big)}_{-s_2}^2\\
\lesssim &  g^*(t^{-1})^{2(2+s_2)}\no{\big( D^m\bar G\big) \phi^t_z}_{-s_2+1}^2\\
\lesssim &  g^*(t^{-1})^{2(2+s_2)}\no{ D^m\bar G}_{-m}\no{ \phi^t_z}_{-s_2+m+1},
\end{split}
\end{eqnarray}
where the last inequality follows from the Cauchy-Schwartz inequality. We notice that $\no{ D^m\bar G}_{-m}\lesssim\no{\bar G}\le 1 $ and, if $s_2-m-1>n$,  the Sobolev's lemma implies that 
 \begin{eqnarray}\begin{split}
\no{ \phi^t_z}^2_{-s_2+m+1}=&\sup\{(|(\phi^t_z, h)| : h\in C^\infty_c, \no{h}_{s_2-m-1}\le 1\}\\
\lesssim& \no{\phi^t_z}\\
=&1.
\end{split}
\end{eqnarray}
Therefore, we get 
 \begin{eqnarray}\Label{6.35c}
\begin{split}
\no{\psi^t_zK_\Om(\cdot , z)D^m\bar G(z)}^2\lesssim g^*(\theta(\delta(z))^{-1})^{2(m+n+4)}.
\end{split}
\end{eqnarray}
Combining \eqref{6.31c} and \eqref{6.35c}, we obtain
\begin{eqnarray}
\begin{split}
\no{\psi^t_z\gamma_z}\lesssim& |R(z)| \delta(z)^{m+1}g^*(\theta(\delta(z))^{-1})^{m+n+4}.
\end{split}
\end{eqnarray}
Choose $\theta=g(\delta^{-1+\eta})^{-1}$, so that  $g^*(\theta(\delta)^{-1})= \delta^{-1+\eta}$ for $\eta> 0$. It follows
\begin{eqnarray}\Label{6.60}
\begin{split}
\no{\psi^t_z\gamma_z}\lesssim& |R(z)| \delta(z)^{\eta m-(1-\eta)n-3+4\eta}.
\end{split}
\end{eqnarray}

In order that the right-hand side of \eqref{6.52f} be $\lesssim 1$, we must have $|R(z)|\lesssim \theta(z)^{-1}$. Letting $R(z)=\theta(z)^{-1}$ in formula \eqref{6.60}, then forces
\begin{eqnarray}\Label{6.61}
g(\delta^{-1+\eta})\delta(z)^{\eta m-(1-\eta)n-3+4\eta}\lesssim 1
\end{eqnarray}
in order that $\no{\psi^t_0\gamma_z}\lesssim 1$. Since $g(\delta^{-1+\eta})=\frac{f}{\log}(\delta^{-1+\eta})\lesssim \delta^{-1} $,  \eqref{6.61} follows by letting $m\to +\infty$. We thus conclude that $\no{\gamma_z}\lesssim 1$, and then from \eqref{6.50}  we get  
$$B_\Om(z, X)\simge |R(z)|=g(\delta^{-1+\eta}).$$
Theorem \ref{t6.12} is proved.

$\hfill\Box$

 \section{How the  Property $(f\T-\M\T-P)^k$ follows from the $(f\T-\M)^k$ estimate}
It is classical that $\epsilon$-subelliptic estimates imply Property $(\tilde{f}\T-1\T-P)^k$ with  $\tilde{f}(\delta^{-1})=\delta^{-\tilde\epsilon}$ for a new $\tilde\epsilon<\epsilon$. This is highly non-trivial and requires a big amount of Catlin's theory. First one proves that $\epsilon$-subelliptic estimates imply D'Angelo's finite type $D$ according to our discussion in Section 6.4.3. 
Here the link between $\epsilon$ and $D$ is accurate and one has $D\le\frac1\epsilon$. \index{Type !  finite type } \index{Notation ! $D$} Next, according to \cite{C87}, we have the converse. The proof goes through the construction of a family of weights, that is, Property $(\tilde f\T-1\T-P)^k$.  Combining with the former step, we get that $\epsilon$-subelliptic estimates imply  Property $(\tilde f\T-1\T-P)^k$; however, the relation of $\tilde\epsilon$ with $D$ and, hence with the initial $\epsilon$, is very rough. This makes very interesting to have at disposal a direct proof of the implication of Property $(\tilde f\T-1\T-P)^k$ from $\epsilon$-subelliptic estimates. In this implication, we also search for a better accuracy about indices. The possible existence of a direct implication was a question raised also by Straube in Vienna in December 2009. Here is the answer. 

\begin{theorem}
Let $\Om\subset\subset \C^n$ be a pseudoconvex domain with a smooth boundary in a neighborhood $U$ of $z_0\in b\Om$. Suppose that $(f\T-1)^1$ holds in $U$ with $f\gg \log$. Then, given a constant $\eta >0$, there exists $U'\subset U$, such that Property $(\tilde f\T-1\T-P)^1$ holds in $U'$ with $\tilde f(\delta^{-1})=\frac{f}{\log^{3/2+\eta}}(\delta^{-1+\eta})$.  
\end{theorem}

{\it Proof.} The notations $K_\Om(z,z$, $\delta(z)$, $U'$ and $ \eta$ are the same as in Theorem \ref{t6.9}. Define
\begin{eqnarray}
\phi(z)=\frac{\log K_\Om(z,z)}{\big(\log (\delta^{-1}(z))\big)^{1+2\eta}}-\frac{1}{\big(\log (\delta^{-1}(z))\big)^{\eta }}
\end{eqnarray}
for $z\in U'$. 
Recall that $K_\Omega(z,z)\simgeq \delta^{-1}(z)$ by \eqref{newter} whereas the opposite estimate is obvious because $\Om$ contains an osculating ball at any boundary point. Thus $\phi(z)\to0$ as $\delta(z)\to0$ (and in particular, $\phi$ is bounded).  
  Moreover,  we have 
\begin{eqnarray}\Label{6.66f}\begin{split}
\di\dib \phi(z)(u^\tau) =&\frac{\di\dib \log K_\Om(z,z) (u^\tau)}{\big(\log (\delta^{-1}(z))\big)^{1+2\eta}}+(1+2\eta)\frac{\log K_\Om(z,z)\cdot \di\dib\delta(z)(u^\tau)}{\delta(z)\big(\log (\delta^{-1}(z))\big)^{2+\eta}}\\
&\hskip3.5cm -\eta\frac{\di\dib \delta(z)(u^\tau)}{\delta(z)\big(\log (\delta^{-1}(z))\big)^{1+\eta }}\\
 =&\frac{\di\dib \log K_\Om(z,z) (u^\tau)}{\big(\log (\delta^{-1}(z))\big)^{1+2\eta}}+\frac{ \di\dib\delta(z)(u^\tau)}{\delta(z)\big(\log (\delta^{-1}(z))\big)^{1+\eta}}\times\\
&\hskip2.5cm\times \left ((1+2\eta)\frac{\log K_\Om(z,z)}{\log \delta^{-1}(z)}-\eta \big(\log \delta^{-1}(z)\big)^{\eta }\right).
\end{split}
\end{eqnarray}
Here, the  last line between brackets  is negative when $z$ approaches $b\Om$ because  its first term stays bounded whereas the second diverges to $-\infty$.
Since $\Om$ is pseudoconvex at $z_0$, then $\di\dib\delta(z)(u^\tau)\le0$. Combining with Theorem \ref{t6.9}, we obtain
\begin{eqnarray}\Label{6.67f}\begin{split}
\di\dib \phi(z)(u^\tau) \ge &\frac{B_\Om(z,u^\tau)^2}{\log(\delta^{-1}(z))^{1+2\eta}}\\
\simge &\frac{(f(\delta^{-1+\eta}(z)))^2}{(\log \delta^{-1+\eta}(z))^2\cdot \log(\delta^{-1}(z))^{1+2\eta}}|u^\tau|^2\\
=&\left(\frac{f}{\log^{\frac{3}{2}+\eta}}\Big(\delta^{-1+\eta}(z)\Big)\right)^2|u^\tau|^2
\end{split}
\end{eqnarray}
for all $z$ near $b\Om$. The inequality \eqref{6.67f} implies the proof of the theorem.\\

$\hfill\Box$

\begin{corollary}
Let $\Om\subset\subset \C^n$ be a pseudoconvex domain with a smooth boundary in a neighborhood of  $z_0\in b\Om$. Suppose that $\epsilon$-subelliptic estimate holds in a neighborhood of $z_0$.  Then,  Property $(f\T-1\T-P)^1$ also holds with $f(\delta^{-1})=\delta^{-\epsilon+\eta}$ for any $\eta>0$.  
\end{corollary}


\chapter{ Global  regularity} 
\label{Chapter7}

\section{Compactness estimate and global regularity}
In this section, we outline the proof of the fact that a compactness estimate implies global regularity. \\

 A global compactness estimate is defined as follows. For every positive number $M$ and for any $u\in C^\infty(\Om)^k\cap \T{Dom}(\dib^*)$, we have
\begin{eqnarray}\Label{7.1}
M\no{u}^2\lesssim Q(u,u) +C_M\no{u}^2_{-1}.\end{eqnarray}
\\
It is well-known that a compactness estimate  implies global regularity.
The idea of the proof is very simple. In order to prove global regularity, it suffices to prove 
\begin{eqnarray}\Label{7.1bis}
\no{u}_s\lesssim \no{\Box u }_s,
\end{eqnarray}
for any $u\in C^\infty(\bar\Om)^k\cap \T{Dom}(\Box)$ and for any integer $s$. 
In fact, by elliptic regularization, \eqref{7.1bis} passes from $C^\infty $- to $H_s$-forms.
When \eqref{7.1bis} is satisfied for any $s$, we say that $N_k$, the inverse to $\Box$ on $k$-forms, is {\it exactly regular}.
 We start by remarking that, since the operator $\Box$ is non-characteristic with respect to the boundary, then 
\begin{eqnarray}\Label{7.2}\no{u}^2_s\lesssim \no{\Box u}_{s-2}^2+||\La^{s-1} Du||^2,\end{eqnarray}
where $D$ is any differential operator of order 1 and $\La^{s-1}$ is the tangential differential operator of order $s-1$. By Lemma \ref{l3.38}, the estimate \eqref{7.1} implies that
\begin{eqnarray}\Label{7.3}M\no{D\La^{-1}u}^2\le Q(u,u)+C_M\no{u}^2_{-1}.\end{eqnarray}
Here $D$ stands for $D_r$ or $\La$ and the conclusion follows from the non-characteristicity of $D_r$.  \\
\index{Regularity ! exactly}

We now  estimate the last term of \eqref{7.2}; we have 
\begin{eqnarray}\begin{split}
M\no{\La^{s-1} Du}^2\lesssim &M\no{D\La^{-1}\La^s u}^2+M\no{u}^2_{s-1}\\
\le & Q(\La^s u,\La^s u)+C_M\no{u}^2_{s-1}\\
\lesssim & (\La^s\Box  u,\La^su)+\no{[\dib, \La^s] u}^2+\no{[\dib, \La^s] u}^2\\
&+\no{[\dib^*,[\dib,\La^s] u}^2+\no{[\dib,[\dib^*,\La^s] u}^2+C_M\no{u}^2_{s-1}\\
\lesssim & \no{\La^s\Box  u}^2+\no{\La^{s-1}Du}^2+\no{\La^{s-2}D^2u}^2+C_M\no{u}^2_{s-1}\\
\lesssim & \no{\Box  u}_s^2+\no{\La^{s-1}Du}^2+C_M\no{u}^2_{s-1},
\end{split}
\end{eqnarray}
where the second inequality follows from \eqref{7.3}. Now, the term $\no{\La^{s-1}Du}^2$ can be absorbed in the left-hand side for $M$  sufficiently large. By induction, we obtain the estimate \eqref{7.1bis}. \\

In Chapter 1, we have introduced a local compactness estimate in a neighborhood $U$ of $z_0\in b\Om$ as
\begin{eqnarray}\Label{7.6}
M\no{u}^2\lesssim Q(u,u) +C_M\no{u}^2_{-1}
\end{eqnarray}
 for any $u\in C_c^\infty(U\cap \Om)^k\cap \T{Dom}(\dib^*)$. The following lemma yields a global estimate from a local estimate.
  
\begin{lemma} Let $\Om$ be a bounded domain. Assume that the estimate \eqref{7.6} holds at any boundary point. Then \eqref{7.1} also holds. 
\end{lemma}
{\it Proof. } Let $\{\zeta_j\}^N_{j=0}$ be a partition of the unity such that $\zeta_0\in C^\infty_c(\Om)$, $\zeta_j\in C^\infty_c(U_j),$ $j=1,...,N$ and
$$\sum_{j=0}^N\zeta_j^2=1 \qquad \T{ on  }\bar\Om,$$
where $\{U_j\}_{j=1,...,N}$ is a covering of $b\Om$.

For $u\in C^\infty(\bar\Om)^k\cap \T{Dom}(\dib^*)$, we wish to show \eqref{7.1}. 
From the interior elliptic regularity of $Q$ we have
$$\no{\zeta_0 u}^2_1\lesssim Q(\zeta_0u,\zeta_0u).$$
On the other hand, by the Cauchy-Schwarz inequality, we have
$$M\no{\zeta_0u}^2\lesssim \no{\zeta_0u}_1^2+C_M\no{\zeta_0u}_{-1}^2.$$
It follows
\begin{eqnarray}\begin{split}
M\no{\zeta_0u}^2\lesssim& Q(\zeta_0u,\zeta_0u)+C_M\no{\zeta_0u}_{-1}^2\\
 \lesssim&Q(u,u)+C_M\no{u}_{-1}^2.
\end{split}
\end{eqnarray}
Similarly, for $j=1,...,N$, using the hypothesis, we have
\begin{eqnarray}\begin{split}
M\no{\zeta_ju}^2\lesssim& Q(\zeta_ju,\zeta_ju)+C_M\no{\zeta_ju}_{-1}^2\\
 \lesssim&Q(u,u)+C_M\no{u}_{-1}^2.
\end{split}
\end{eqnarray}
Summing up over $j$, we get the proof of the lemma.

$\hfill\Box$
\\

We have thus proved that local compactness estimate over a covering of $b\Om$ implies global compactness. In turn, this latter yields \eqref{7.1bis} and hence global regularity. We collect in a single statement the conclusions of this section.

\begin{theorem}
\begin{enumerate}
  \item   Let $\Om\subset\C^n$ be a smoothly  bounded $q$-pseudoconvex (resp. $q$-pseudoconcave) domain  at any boundary point. Assume that a compactness estimate  holds on $(0,k)$-forms  with $q\le k\le n-1$ (resp. $1 \le k\le q $) in a neighborhood of any boundary point . Then the $\dib$-Neumann operator $N_k$ on $(0,k)$-forms is exactly, globally regular.
  \item   Let $\Om\subset\C^n$, $n\ge 3$, be a smooth,  bounded  annulus  defined by $\Om=\Om_1\setminus \Om_2$ where $\bar\Om_2\subset\Om_1$ and $\Om_1$ is $p$-pseudoconvex and $\Om_2$ is $(n-q-1)$-pseudoconvex. Assume that a compactness estimate  holds on $(0,k)$-forms  with $p \le k\le q $  in a neighborhood of any boundary point. Then the $\dib$-Neumann operator $N_k$ on $(0,k)$-forms, $p \le k\le q $, is exactly, globally regular.
\item Let $M\subset\C^n$, $n\ge3$,  be a $C^\infty$, compact,  $q$-pseudoconvex hypersurface. Assume that a compactness estimate  holds on $(0,k)$-forms  with $q\le k\le n-1-q$  in a neighborhood of any point. Then the Green operator $G_k$ ($G_k:=\Box_b^{-1}$) on $(0,k)$-form is exactly, globally regular.
\end{enumerate}
\end{theorem}

\section{Failure of compactness estimate in presence of global regularity}
Is compactness estimate necessary for global regularity? Answer is no. We start from  the following statement about failure of compactness estimate.
\begin{proposition}
\Label{reinhardt}
Let $\Om$ be a smooth bounded pseudoconvex domain of $\C^n$ with a ``$n-1$-Reinhardt flat"  piece of boundary. This means that, in some choice of coordinates
$$
b\Om\supset b\Bbb D\times \Bbb D^{n-1}_\epsilon,
$$
where  $\Bbb D$ is the unit polydisc in $\C$ and $\Bbb D^{n-1}_\epsilon$ the $\epsilon$-polydisc in $\C^{n-1}$. Then, compactness of $N_k$ does not hold. \index{Notation ! \Bbb D^k}
\end{proposition}


 \begin{remark}
 This result generalizes the one by Krantz in \cite{Kr88} and is close to the further development contained in \cite{BS99}.
 In the original statement, Reinhardt domains having a flat portion of the boundary, are considered. Here, it is not required that the domain is fully Reinhardt but, we still ask for a flat Reinhardt portion of the boundary.
 \end{remark}
{\it Proof of Proposition~\ref{reinhardt}.} We prove the proposition for the case $n=2$ and $k=1$, since the general proof is identical. 
 Let $\psi\in C^\infty_c(\R)$ satisfy $\psi(t)=\begin{cases}1 \T{~if~} t\le \frac{1}{4} \\ 0 \T{~if~} t>\frac{1}{2}.\end{cases}$  For  $m\in \Z^+$,  set 
$$u_m(z_1,z_2)=\sqrt{2(m+1)}z_1^{m}\psi(|z_2|^2)d\bar z_2.$$
   
Then,  we have $r_{z_1}(u_m)_1+r_{z_2}(u_m)_2=0$ on $b\Om$, that is, $u_m\in C^\infty(\bar\Om)^1\cap \T{Dom}(\dib^*)$. Moreover, we have
\begin{eqnarray}\Label{7.10f}
\no{u_m}^2=2(m+1)\int_{\Om} |z_1|^{2m}\psi(|z_2|^2)dV\ge 2 (m+1)\int_{D(0, 1)\times D(0,\frac{1}{2})}|z_1|^{2m} dV \simeq1.
\end{eqnarray}
On the other hand, one checks readily that $\dib u_m=0$ and $\dib^* u_m=\sqrt{2(m+1)}z_1^m\partial_{z_2}\chi$; this yields $\no{\dib^* u_m}^2\simleq 1$. In conclusion,
\begin{equation}
\Label{newquater}
||u_m||_0\simeq Q(u_m,u_m).
\end{equation}
We remark now that, by $H_0$-boundedness, there must exist a convergent subsequence $\{u_{m_k}\}$. Since $\Om$ is circular with respect to $z_1$, and since integration of $z_1^{m}\bar z_1^{\tilde m}$ over the circle is $0$ unless $m=\tilde m$, then the sequence is orthogonal and therefore the limit must be $0$. By the compactness of the embedding $H_{-1}\hookrightarrow H_0$, we conclude $||u_{m_k}||_{-1}\to0$. This, in combination with \eqref{newquater}, violates compactness.

$\hfill\Box$
\\

As for global regularity, there are several criteria which do not require compactness. The first is the so-called Condition (T) (cf. \cite{CS01} p. 129). Under the choice of a defining function $r$ of $\Om$ and of a normal vector field $L_n=\underset{j=1}{\overset n\sum}r_{\bar z_j}\partial_{z_j}$, this is expressed as follows. For any $\epsilon$ there is a vector field $T=T_\epsilon$, tangent to $b\Om$ and whose component along $L_n-\bar L_n$ has uniformly positive lower bound, such that
\begin{equation}
\Label{(T)}
\left|\langle [T,S],L_n\rangle|_{b\Om}\right|<\epsilon
\end{equation}
for any $S\in\C\otimes T\C^n$. We refer to \cite{CS01} Theorem 6.2.1 for the proof that Condition (T) implies global regularity. An easy application of this result is
\begin{proposition}
\Label{reinhardtcomplete}
Let $\Omega$ be a Reinhardt complete, smooth bounded, pseudoconvex domain of $\C^n$.  Then, the $\dib$-Neumann operator $N_k$ on $(0, k)$-forms is exactly regular in Sobolev norms, that is
\begin{eqnarray}
\no{N_k\alpha }_s \le  C_s \no{\alpha}_s ,
\end{eqnarray}
for $s \ge 0$ and all $\alpha \in H_s(\Om)^k$.
\end{proposition}
{\it Proof.}
The proof goes through Condition (T). Since $\Om$ is Reinhardt, then in particular it is circular, that is, invariant under multiplication by $e^{i\theta}\in S^1$ and therefore $T:=i\underset{j=1}{\overset n\sum}z_j\partial_{z_j}-i\underset{j=1}{\overset n\sum}\bar z_j\partial_{\bar z_j}$, when restricted to $b\Om$, is tangent to $b\Om$. It is an easy exercice to check that in order that Condition (T) is fulfilled, it suffices to show that $T$ is not complex tangential, that is,
\begin{equation}
\Label{complextangential}
\underset{j=1}{\overset n\sum}z_j\partial_{z_j}(r)|_{b\Om}\neq0.
\end{equation}
If, instead, we have $\underset{j=1}{\overset n\sum}z_j\partial_{z_j}(r)|_{b\Om}= 0$ at some $z^o\in b\Om$, then the vector $z^o-0$ is orthogonal to $\partial r(z^o)$ and therefore, there are other points $z^1\in \Om$ such that $|z^1_j|>|z^o_j|$ for any $j$. Since $\Om$ is Reinhardt complete, this is a contradiction which proves \eqref{complextangential}. Thus Condition (T) is verified and the proposition is proved.

$\hfill\Box$
\\

The most general criterion for global regularity is due to Straube \cite{S08}.

\begin{theorem}\Label{t7.4e}{\bf  [Straube 08]}
Let $\Om$  be a smooth bounded pseudoconvex domain in $\C^n$, $r$ a defining
function for  $\Om$. Let $1 \le k \le n$. Assume that there is a constant $C$ such that for all $\epsilon > 0$  there exist a defining function $r^\epsilon$ for $\Om$
 and a constant $C_\epsilon$  with
\begin{eqnarray}
 1/C \le | \nabla r^\epsilon| \le C ~~~\T{~~on~~} b\Om ,
\end{eqnarray}
and
\begin{eqnarray}\Label{6.13g}
\sumK \left\|  \sum_{ij=1}^n \overline{r^\epsilon_{z_i\bar z_j}r_{\bar z_i}}u_{jK}   \right\|^2\le \epsilon(\no{\dib u}^2+\no{\dib^* u}^2)+C_\epsilon\no{u}_{-1}^2
\end{eqnarray}
for all $u \in  C^\infty(\bar\Om)^k\cap \T{Dom}(\dib^*)$. Then the $\dib$-Neumann operator $N_k$ on $(0, k)$-forms is exactly regular in Sobolev norms.
\end{theorem}

We  postpone to the next section a modified proof of Straube's theorem but first discuss of some examples. A general class of domains which enjoy global regularity are convex domains. They belong to the larger class of domains endowed with a plurisubharmonic defining function which also have global regularity (cf. e.g. \cite{BS91}). Historically, this fact is proved by the aid of Condition (T). For the selfcontainedness of our exposition, we prefere to offer a proof which uses Straube's Theorem.  
\begin{theorem}\Label{t7.3e}{\bf [Boas-Straube 91]} Let $\Om$ be a smooth bounded domain in $\C^n$ admitting a defining function $r$ that is plurisubharmonic on the boundary, that is, $\di\dib r(z)\geq0$ for any $z\in b\Om$. Then the $\dib$-Neumann operator $N_k$ is exactly, globally regular for any $k\ge1$.
\end{theorem}

{\it Proof.} By Theorems \ref{t3.12} and \ref{t7.4e}, we only need to show that for any $\epsilon>0$, there exists a family $\{\Phi^\epsilon\}$ such that   
\begin{eqnarray}
\begin{cases}
|\Phi^\epsilon|\lesssim 1\\
\di\dib\Phi^\epsilon(u^\tau, u^\tau) \simge \dfrac{1}{\epsilon}\big| \overset{n}{\underset{ij=1}{\sum}}\overline{r_{z_i\bar z_j} r_{\bar z_i}}u^\tau_{j}\big |^2
\end{cases}~~~~\T{~~on~~} b\Om
\end{eqnarray}
for any $u=\sum u_jd\bar z_j\in C^\infty(\bar\Om)^1\cap \T{Dom}(\dib^*)$. \\

We define $ \Phi^\epsilon :=\dfrac{r}{\epsilon}$. Since $\di\dib r|_{b\Om}\ge0$, then by Cauchy-Schwartz
\begin{eqnarray}
\di\dib r(u^\tau, u^\tau) \di\dib r(\dib r, \dib r)\ge \left|\overline{\di\dib r (
\dib r,u^\tau)}\right|^2=\left|\overset{n}{\underset{ij=1}{\sum}}\overline{r_{z_i\bar z_j} r_{\bar z_i}}u^\tau_{j} \right|^2\,\,\T{ on $b\Om$}.
\end{eqnarray}
Since $\di\dib (\dib r, \dib r)$  is bounded on $b\Om$, then

\begin{eqnarray}\begin{split}
\di\dib\Phi^\epsilon(u^\tau, u^\tau)=& \dfrac{\di\dib r(u^\tau,u^\tau)}{\epsilon}\\
\simgeq& \dfrac{\di\dib r(u^\tau,u^\tau) \di\dib r(\dib r, \dib r)}{\epsilon }\\
\simge&\dfrac{1}{\epsilon}\left|\overset{n}{\underset{ij=1}{\sum}}\overline{r_{z_i\bar z_j} r_{\bar z_i}}u^\tau_{j} \right|^2.
\end{split}
\end{eqnarray}
This concludes the proof of Theorem \ref{t7.3e}.

$\hfill\Box$
\begin{example}
Let $\Om$ be the domain in $\C^2$ defined by $|z_1|^2+\zeta(z_1,z_2)-1<0$ where
$$
\zeta=\begin{cases}
0\quad\T{if $|z_2|<\frac12$},
\\
1\quad\T{if $|z_2|\ge1$},
\end{cases}
$$
where we also assume, as a general fact, that $\Om$ is pseudoconvex. In this situation, according to Theorem~\ref{reinhardt}, there is no compactness. However, global regularity can be achieved in several cases as for instance for
\begin{itemize}
\item[(i)] $\zeta=\zeta(|z_1|,|z_2|)$ (so that $\Om$ is Reinhardt) and $\zeta$ is decreasing with respect to  both $|z_1|$ and $|z_2|$ (so that $\Om$ is complete),
\item[(ii)] $\zeta$ is general 
but, for any vector $w\in\C^2$, it satisfies $|w_1|^2+\di\dib\zeta(w,w)\ge0$ (so that $\Om$ admits a plurisubharmonic defining function).
\end{itemize}
\end{example}

\section{Weak compactness estimate and global regularity}

The purpose of this section is to show that  global regularity follows from an estimate which is weaker than compactness.  The proof that we give is a variation of a Straube's \cite{S08} result.

\begin{theorem}\Label{t7.4}
Let $\Om$ be a smooth, bounded, pseudoconvex domain in $\C^n$. Assume that for any positive constant $\epsilon>0$ there exists a defining function $r^\epsilon$ of $\Om$ with $\underset{k}{\sum}|r^\epsilon_{z_j}|^2\cong 1$ on $b\Om$ such that 
\begin{eqnarray}\Label{6.18g}
\sumK \Big\|\sum^n_{ij} r^\epsilon_{z_i\bar z_j }r^\epsilon_{\bar z_i} \bar u_{jK}\Big\|^2\le \epsilon\left( \no{\dib u}^2+\no{\dib^* u}^2\right)+C_\epsilon\no{u}^2_{-1},
\end{eqnarray}
for any $u\in C^\infty(\bar \Om)^k\cap \T{Dom}(\dib^*)$. Then the Bergman projection $P_{l-1}$ is exactly, globally regular for any $l\ge k$.
\end{theorem}

\begin{remark}
We remark that Theorem \ref{t7.4e} and \ref{t7.4} are equivalent. In fact, it is easy to check that \eqref{6.13g} and \eqref{6.18g} are equivalent;  and by Boas-Strauble \cite{BS90}, $N_k$ is exactly , globally regular if and only if $P_{k-1}, P_k, P_{k+1}$ are exactly, globally regular. 
\end{remark}

Before starting the proof, we need to fix some notations and state some preliminary result whose proof is immediate.
Let $U$ be a neighborhood of $b\Om$. For  any $\epsilon>0$, we may assume that the defining function $r:=r^\epsilon$ of $\Om$ satisfies $\sum_{k=1}^{n}|r_{z_k}|^2\not=0$ on $U$. We define a system of (1,0) vector fields by
$$N=\frac{1}{\sum_{k=1}^n |r_{z_k}|^2}\sum_{k=1}^n r_{\bar z_k}\frac{\di }{\di z_k},\quad T=N-\bar N \quad \T{and } L_j=\frac{\di }{\di z_j}-r_{z_j}N$$  
for $j=1,...,n$. Notice that $T$ and $L_j, 1\le j\le n$, are tangential and that  $\bar T=-T$.  

First, we consider $u\in C_c^\infty(U\cap \bar\Om)^k$. Using integration by parts, we get 
\begin{eqnarray}
\begin{split}\no { L_j u}^2\lesssim ([\bar L_j, L_j]u,u)+\no{\bar L_j u}^2+\no{u}^2\\
\lesssim \no{u}_1.\no{u}+\no{\bar L_j u}^2.\\
\end{split}
\end{eqnarray}
Denote $\S=\T{span}\{L_1,...,L_n, \frac{\di }{\di \bar z_1},...,\frac{\di }{\di \bar z_1}, Id\}$. 

\begin{proposition} For $s\in \Bbb N$ and $S\in \S$,  we have
\begin{eqnarray}\Label{7.11}
\no{S u}_{s-1}^2\lesssim \no{\dib u}^2_{s-1}+\no{\dib^* u}^2_{s-1}+\no{u}_{s-1}\cdot\no{u}_s,\end{eqnarray}
and
 \begin{eqnarray}\Label{7.12}
\no{ u}_s^2\lesssim \no{\dib u}^2_{s-1}+\no{\dib^* u}^2_{s-1}+\no{u}_{s-1}.\no{u}_s +\no{T^su}^2,\end{eqnarray}
for any $u\in C_c^\infty(U\cap\bar \Om)^k\cap \T{Dom}(\dib^*)$. Moreover, for $u\in C^\infty_c(\Om)^k\cap \T{Dom}(\dib^*)$, we have
\begin{eqnarray}\Label{7.13}
\no{u}^2_{s+1}\lesssim \no{\dib u}^2_{s}+\no{\dib^* u}^2_s+\no{u}^2_{s-1}.\end{eqnarray}
\end{proposition}
{\it Proof. } The proof of \eqref{7.11} follows from the fact that 
\begin{eqnarray}\Label{7.14b}
\no{\frac{\di u}{\di \bar z_j}}_{s-1}^2\lesssim \no{\dib u}^2_{s-1}+\no{\dib^* u}^2_{s-1}+\no{u}_{s-1}^2\end{eqnarray}
for $j=1, ...,n$, and 
\begin{eqnarray}\Label{7.15b}
\no{L_j u}_{s-1}^2\lesssim \no{\dib u}^2_{s-1}+\no{\dib^* u}^2_{s-1}+\no{u}_{s-1}\cdot\no{u}_s,\end{eqnarray}
for $j=1,...,n-1$. (The inequalities \eqref{7.14b} and \eqref{7.15b} may be found in \cite{BS91}, p.83 or in \cite{CS01}, Section 6.2.) \\

The proof of \eqref{7.12} follows by the induction in $j$ from
 \begin{eqnarray}\Label{7.15bis}
\no{ T^ju}_{s-j}^2\lesssim \no{\dib u}^2_{s-1}+\no{\dib^* u}^2_{s-1}+\no{u}_{s-1}\no{u}_s +\no{T^{j+1}u}_{s-j-1}^2.
\end{eqnarray}  The last estimate \eqref{7.13} follows from ellipticity in the interior.\\

$\hfill\Box$

\begin{lemma}
\begin{enumerate}
  \item[(i)] We have
\begin{eqnarray}
\begin{split}
[\frac{\di}{\di\bar   z_j}, T]=& \theta_j T+S_{j},
\end{split}
\end{eqnarray}
where $\theta_j=-\frac{1 }{\underset{k}{\sum}|r_{z_k}|^2}\underset{i}{\sum} r_{z_i\bar z_j}r_{\bar z_i}$, and $S_{j},\in \S$.
\item[(ii)] We also have
\begin{eqnarray}
\begin{split}
[\frac{\di}{\di  \bar z_j}, T^{2s}]=& 2s\theta_j T^{s+1}+ A^{2s-1}_jS_{j\,s},\\
\end{split}
\end{eqnarray}
where $A^{2s-1}_j$ is a tangential differential  operator  of order $2s-1$, and $S_{j\,s}\in \S$.
\end{enumerate}
\end{lemma}
{\it Proof. } By differentiating and using formula $\frac{\di}{\di z_j}=r_{z_j}T+L_j+r_{z_j}\bar L_n$, we get the proof of the lemma. 
\\

$\hfill\Box$ 
\\

{\it Proof of Theorem \ref{t7.4}.} For any $s\in \Bbb N$, we want to  show by  induction that $P_l$ is exactly continuous on $H_s(\Om)^l$ for $l=n-1,..., k-1$.\\

We start from the remark that, when $u$ has degree $n$, the condition $u\in\T{Dom}(\dib^*)$ coincides with the Dirichlet condition $u|_{b\Om}\equiv0$. 
Thus $Q$ is elliptic on $n$-forms and hence 
 $P_{n-1}=I-\dib^* N_n\dib$ is  exactly continuous on $H_s(\Om)^{n-1}$. The induction hypothesis is that $P_k$ is exactly continuous on $H_s(\Om)^k$.  We need to show that $P_{k-1}$ is exactly continuous on $H_s(\Om)^k$ for any $s$.\\

We  first prove that 
\begin{eqnarray}\Label{7.16}
\no{P_{k-1}\alpha}^2_s\lesssim \no{\alpha}^2_s+\no{P_{k-1}\alpha}^2_{s-1}\end{eqnarray}
 for any $\alpha\in C_c^\infty(U\cap \bar\Om)^{k-1}$. 
We have  $$(I-P_{k-1})\alpha =\dib^*N_{k}\dib\alpha\in C^\infty_c(U\cap \bar\Om)^{k-1}\cap \T{Dom}(\dib^*).$$ 
Using \eqref{7.12}, for $u$ replaced by $(I-P_{k-1})\alpha$, we get

\begin{eqnarray}
\begin{split}
\no{(I-P_{k-1})\alpha}^2_s\lesssim &\no{\dib \dib^*N_k \dib \alpha }_{s-1}^2+ \no{(I-P_{k-1})\alpha}^2_{s-1}+\no{T^s(I-P_{k-1})\alpha}^2\\
\lesssim &\no{ \alpha }_{s}^2+ \no{P_{k-1}\alpha}^2_{s-1}+\no{T^sP_{k-1}\alpha}^2.
\end{split}
\end{eqnarray}
Hence
\begin{eqnarray}\Label{7.18}
\begin{split}
\no{P_{k-1}\alpha}^2_s\lesssim &\no{ \alpha }_{s}^2+ \no{P_{k-1}\alpha}^2_{s-1}+\no{T^sP_{k-1}\alpha}^2.
\end{split}
\end{eqnarray}
Similarly, by using \eqref{7.11}, we get 
\begin{eqnarray}\Label{7.19}
\begin{split}
\no{SP_{k-1}\alpha}^2_{s-1}\lesssim &\no{ \alpha }_{s}^2+\no{P_{k-1}\alpha}_s\no{P_{k-1}\alpha}_{s-1}
\end{split}
\end{eqnarray}
for any $S\in \S$. \\

Now, we estimate the last term of \eqref{7.18}. We have
\begin{eqnarray}\Label{7.20}
\begin{split}
\no{T^sP_{k-1}\alpha}^2= &(T^sP_{k-1}\alpha, T^s\alpha)-(T^sP_{k-1}\alpha, T^s \dib^*N_k\dib\alpha)\\
= &(T^sP_{k-1}\alpha, T^s\alpha)-((T^s)^*T^s\dib P_{k-1}\alpha, N_k\dib\alpha)\\
&-([\dib, (T^s)^*T^s]P_{k-1}\alpha,N_k\dib\alpha).
\end{split}
\end{eqnarray}
We set $\bar\theta_j=\sum_ir_{z_j\bar z_i}r_{\bar z_i}$ and define an operator $\bar\Theta:\,\A^{0,k-1}\to \A^{0,k}$ by
\begin{equation*}
\bar\Theta u=\underset{|K|=k-2}{{\sum}'}\underset{i<j}\sum\left(\bar\theta_iu_{jK}-\bar\theta_ju_{iK}\right)d\bar z_i\wedge d\bar z_j\wedge d\bar z_K.
\end{equation*}
We also define operators $\bar I_j:\A^{0,k-1}\to \A^{0,k}$ by
 $$\bar I_ju=\sum_{|J|=k-1} u_J d\bar z_j\we d\bar z_J.$$
We have
\begin{eqnarray}
\begin{split}
[\dib, (T^s)^*T^s ]=&[\dib, T^{2s} ]+[\dib, ((T^s)^*-T^s)T^s ]\\
=&2s\bar\Theta T^{2s}+ \sum_jA_{j}^{2s}S_j\bar I_j\\
=&2s(T^s)^*\bar\Theta T^{2s}+ \sum_j(C^{s}_j)^*B_{j}^{s-1}S_j\bar I_j.
\end{split}
\end{eqnarray}
where $A_j^{2s-1}, B_j^{s-1}$ and $C_j^{s}$ are  tangential operators of order $2s-1$, $s-1$, and $s$ respectively.

Thus, we can continue  \eqref{7.20} by
\begin{eqnarray}\Label{7.22}
\begin{split}
=&(T^sP_{k-1}\alpha, T^s\alpha)-(\bar\Theta 2sT^sP_{k-1}\alpha,T^s N_k\dib\alpha)+\sum_{j}(B^{s-1}_jS_j \bar I_jP_{k-1}\alpha,C^{s}_j X_j^sN_k\dib\alpha)\\
\le&  \epsilon \Big(\no{P_{k-1}\alpha}_s^2+\no{N_k\dib\alpha}^2_s\Big)+\frac{1}{\epsilon}\Big(\no{\alpha}_s^2+\sum_j\no{S_j  P_{k-1}\alpha}^2_{s-1}+\no{\bar\Theta^* T^sN_k\dib\alpha}^2\Big).
\end{split}
\end{eqnarray}

Using the hypothesis of the theorem,  we get
\begin{eqnarray}\Label{7.23}
\begin{split}
\no{\bar\Theta^* T^sN_k\dib\alpha}^2\lesssim&\epsilon( \no{\dib  T^sN_k\dib\alpha}^2+\no{\dib^* T^sN_k\dib\alpha}^2)+C_\epsilon \no{T^sN_k\dib\alpha}^2_{-1}\\
\lesssim& \epsilon(\no{ T^s\dib N_k\dib\alpha}^2+\no{[\dib,  T^s]N_k\dib\alpha}^2+\no{T^s\dib^*N_k\dib\alpha }^2
\\
&\qquad+\no{[\dib^* ,T^s]N_k\dib\alpha}^2)+C_\epsilon \no{T^sN_k\dib\alpha}^2_{-1}
\\
\lesssim& \epsilon(\no{N_k\dib\alpha}^2_s+\no{P_{k-1}\alpha }_s^2+\no{\alpha}_s^2)+C_\epsilon \no{N_k\dib\alpha}^2_{s-1},
\end{split}
\end{eqnarray}
where we are using, among other things, that $  \dib N_k\dib\alpha=0$.
Combining \eqref{7.18}, \eqref{7.22} and \eqref{7.23}, we obtain
\begin{eqnarray}\Label{7.24}
\begin{split}
\no{P_{k-1}\alpha}^2_s\lesssim &\no{ \alpha }_{s}^2+\no{P_{k-1}\alpha}^2_{s-1}+\sqrt{\epsilon}\no{N_k\dib\alpha}^2_s+C_\epsilon\no{N_k\dib\alpha}^2_{s-1}.
\end{split}
\end{eqnarray}

We use now   Boas-Straube's formula in \cite{BS91}. For $\alpha\in H^s(\Om)^{k-1}$ we have 
$$N_k\dib\alpha=P_kw_tN_{t,k}\dib w_{-t}(I-P_{k-1})\alpha,$$
where $N_{t,k}$ is the solution operator to the weighted $\dib$-Neumann problem with weight $w_t(z)=\exp(-t|z|^2)$. By the result of Kohn \cite{K73}, we have that $N_{t,k}\dib$ is also exactly continuous on $H_s(\Om)^{k-1}$. It follows
\begin{eqnarray}\Label{7.25}
\begin{split}
\no{N_k\dib\alpha}_s\lesssim \no{\alpha}_s+\no{P_{k-1}\alpha}_s,
\end{split}
\end{eqnarray}
for any $\alpha\in H^s(\Om)^{k-1}$. From \eqref{7.24} and \eqref{7.25}, we get \eqref{7.16}.\\

Finally, for any $\alpha\in C^\infty(\bar\Om)^k$ we can write $\alpha=\chi_1\alpha+\chi_2\alpha$ where $\chi_1\in C_c^\infty(U\cap \bar\Om)^k$ and $\chi_2\in C_c^\infty(\Om)^k$ where $U$ is a neighborhood of $b\Om$. 

We have 
\begin{eqnarray}\Label{7.26}
\begin{split}
\no{P_{k-1}\alpha}^2_s\lesssim & \no{P_{k-1}(\chi_1\alpha)}^2_s+\no{P_{k-1}(\chi_2\alpha)}^2_s\\
\lesssim &\no{\chi_1\alpha}^2_s+\no{P_{k-1}(\chi_1\alpha)}^2_{s-1}+ \no{\chi_2\alpha}_s+\no{P_{k-1}(\chi_2\alpha)}^2_{s-1}\\
\lesssim &\no{\alpha}^2_s+\no{P_{k-1}(\chi_1\alpha)}^2_{s-1}+\no{P_{k-1}(\chi_2\alpha)}^2_{s-1},\\
\end{split}
\end{eqnarray}
for any $\alpha\in H^s(\Om)^{k-1}$.\\

Since 
$\no{\dib^*N_k\dib\alpha}^2=(\dib^*N_k\dib\alpha,\dib^*N_k\dib\alpha)=(\dib \alpha, N_k\dib\alpha)=(\alpha,\dib^*N_k\dib\alpha)
$, then $\no{P_{k-1}\alpha}\lesssim \no{\alpha}$.
We assume  by induction that $\no{P_{k-1}\alpha}_{s-1}^2\lesssim\no{\alpha}_{s-1}^2$ for any $\alpha\in C^\infty(\bar\Om)^k$; then by \eqref{7.26}, we get 
\begin{eqnarray}\Label{7.27}
\begin{split}
\no{P_{k-1}\alpha}^2_s\lesssim & \no{\alpha}_s^2,
\end{split}
\end{eqnarray}
for any $\alpha\in C^\infty(\bar\Om)^k$. Using the method of elliptic regularization as in \cite{KN65} and \cite{FK72} the a priori estimate \eqref{7.27} becomes an actual estimate which yields the conclusion of the proof of Theorem~\ref{t7.4}.

$\hfill\Box$



\chapter{The Property $(f\T-\M\T-P)^k$ in some classes of domains} 
\label{Chapter5}
In this chapter, we rivisit classical estimates in the framework of the method of the weights and exhibit some new classes of domains which enjoy $(f\T-\M\T-P)^k$ Property. We also discuss about optimal subelliptic estimates.

\section{Domains which satisfy the $Z(k)$ condition}
In this section, we consider a domain satisfying the $Z(k)$ condition at a boundary point. This class of domains is probably the simplest of non-pseudoconvex domains which enjoy $\frac12$-subelliptic estimates. We refer to \cite{FK72}  for the definition of this property which is also recalled here in \ref{e2.2}. 
\begin{theorem}\Label{t5.1} A domain $\Om\subset\C^n$ satisfies $Z(k)$ condition at $z_0\in b\Om$ if and only if in a neighborhood $U$ of $z_0$ the $\frac12$-subelliptic estimate holds, that is, 
$$|||u|||^2_{1/2}\lesssim Q(u,u)$$
 for any $u\in C_c^\infty(U \cap\bar\Om )^k\cap \T{Dom}(\dib^*)$.
\end{theorem}
This is a classical result about non-pseudoconvex domains that can be found in  \cite{H65}, \cite{FK72} and others. In this section, we give a new way to get $\frac{1}{2}$-subelliptic estimates by constructing a family of weights as required by Property $(f\T-\M\T-P)^k$ under the choice $\M=1$ . \\

{\it Proof. } Let $\Om$ satisfy $Z(k)$ condition at $z_0\in b\Om$. Then $\Om$ is strongly $k$-pseudoconvex or strongly $k$-pseudoconcave. There are a number $q_o\not=k$ and a neighborhood $U$ of $z_0$ such that
\begin{eqnarray}
\sumK\sum_{ij=1}^{n-1}r_{ij}u_{iK}\bar u_{jK}-\sum_{j=1}^{q_o}r_{jj}|u|^2\simge |u|^2 \T{~on~} U\cap \bar\Om
\end{eqnarray}
for any $u\in C_0^\infty(U\cap \bar\Om)^k\cap \T{Dom}(\dib^*)$.  So, what is left  to show, is that  $\Om$ satisfies Property $(f\T-1\T-P)^k$ at $z_0$ with $f(\delta^{-1})=\delta^{-1/2}$. \\

We define $\Phi^\delta=-\frac{r}{\delta}$. For $z\in S_\delta$, we see that $\Phi^\delta$ is absolutely bounded and   
\begin{eqnarray}
\begin{split}
H^{\Phi^\delta}_{q_o}(u^\tau)\simge&\delta^{-1}\left(\sumK\sum_{ij=1}^{n-1}r_{ij}u^\tau_{iK}\bar u^\tau_{jK}-\sum_{j=1}^{q_o}r_{jj}|u^\tau|^2\right)\\
\simge&\delta^{-1}|u^\tau|^2
\end{split}
\end{eqnarray}
 Moreover, $L_j(\Phi^\delta)=0$ for any $j\leq n-1$.  Then the family $\{\Phi^\delta\}_{\delta>0}$ satisfies Property $(f\T-1\T-P)^k$.

$\hfill\Box$\\

\begin{corollary} Let $M$ be a smooth hypersurface in $\C^n$. Assume that $M$ satisfies $Z(k)$ and $Z(n-1-k)$ condition at $z_0$. Then there is a neigborhood $U$ of $z_0$ such that
$$\no{u}_{b, 1/2}^2\lesssim Q_b(u,u)$$
for any $u\in C^\infty_c(U\cap M)^k$.
\end{corollary}
The proof of Corollary 5.2 follows from Theorem 5.1 and Corolary \ref{c4.12}. 

\section{$q$-decoupled-pseudoconvex/concave domains}

 Let $\Om \subset\C^n $ be defined in a neighborhood of $z_0$ by
\begin{eqnarray}\Label{q1q2}
r=2\T{Re}z_{n}-{a}(z_1,\dots, z_{q_o})+{b}(z_{q_o+1},\dots, z_{n-1})<0
\end{eqnarray}
where $a$ and $b$ are real functions such that $\partial\bar\partial a$ and $\partial\bar\partial b$ are $\geq0$. 

\begin{proposition} \Label{p5.3} The domain $\Om$, defined by \eqref{q1q2},  is $(q_o+1)$-pseudoconvex  and  $(q_o-1)$-pseudoconcave at $z_0$.  
\end{proposition}
{\it Proof.}  We consider the basis of vector fields 
$$L_i=\frac{\di}{\di z_j}-r_{z_j}\frac{\di}{\di z_n}, j=1...n-1 \T{ and } L_n=\frac{\di}{\di z_n}.$$  Let $\om_1,...,\om_n=\di r$ be the dual (1,0) forms of these vector fields. We may choose the Hermitian metric in which $\om_1,...,\om_n$ are orthonormal. Then 
$$\di\dib r=-\sum_{ij=1}^{q_o} a_{ij}\om_i\we \bar\om_j+\sum_{ij=q_o+1}^{n-1} b_{ij}\om_i\we \bar\om_j$$
where $a_{ij}$ and $b_{ij}$ are the Levi matrices of $a$ and $b$ respectively. It follows, for any k-form $u$
\begin{eqnarray}
H_{q_o}^{r}(u^\tau)=\sum_{j=1}^{q_o}a_{jj}|u^\tau|^2-\sumK\sum_{ij=1}^{q_o} a_{ij}u^\tau_{iK}\bar u^\tau_{jK}+ \sumK\sum_{ij=q_o+1}^{n-1}b_{ij}u^\tau_{iK}\bar u^\tau_{jK}.
\end{eqnarray}
We prove now that $H^r_{q_o+1}(u^\tau)\ge0$ for $k\ge q_o+1$, that is, that $\Om$ is $(q_o+1)$-pseudoconvex. (The proof of the $(q_o-1)$pseudoconcavity of $\Om$ is similar.) In fact, extending  $(a_{ij})$ from a $q_o\times q_o$ to an $(n-1)\times(n-1)$ matrix by adding zeroes, and using the positivity of the resulting 2-form, we conclude
$$
\sumK\sum_{i,j=1}^{q_o}a_{ij}u^\tau_{iK}\bar u^\tau_{jK}\le\sum_{j=1}^{q_o}a_{jj}|u^\tau|^2.
$$
It follows
\begin{equation}
\Label{5.4bis}
H^r_{q_o}(u^\tau)\geq  \sumK\sum_{i,j=q_o+1}^{n-1}b_{ij} u^\tau_{iK}\bar u^\tau_{jK}.
\end{equation}
Note that, $b\Om$ being {\it rigid}, \eqref{5.4bis} holds not only on $b\Om$, but also on $\bar\Omega$ near $b\Om$.

$\hfill\Box$

\begin{definition} $\Om$ is said to be {\it $q_2$-decoupled-pseudoconvex}  (resp. {\it $q_1$-decoupled-pseudoconcave}) at $z_0$ if there are  functions $h_j(z_j)$  such that \index{Domain ! $q$-decoupled-pseudoconvex/concave}
${b}(z_{q_o+1},\dots, z_{n-1})=\sum_{j=q_2}^{n-1}h_j(z_j)$
(resp. $ {a}(z_1,\dots, z_{q_o})=\sum_{j=1}^{q_o}h_j(z_j)$)
\end{definition}
\begin{remark}Since $\di\dib a$ and $\di\dib b$ are $\ge 0$, then $\frac{\di^2  h_j}{\di z_j\di\bar z_j}\ge 0$ for any $j$.
\end{remark}
\begin{remark}If $\Om$ is $(q_o+1)$-decoupled-pseudoconvex  then from \eqref{5.4bis} we get
\begin{eqnarray}\Label{5.7}
H^r_{q_o}(u^\tau)\simge \sum_{j=k}^{n-1}\frac{\di^2 h_j}{\di z_j\di\bar z_j} \sumK |u^\tau_{jK}|^2\ge 0 \T{  on  } U\cap \bar\Om
\end{eqnarray}
 for any $u^\tau$ of degree  $k\ge q_o+1$. On the other hand, if $\Om$ is $(q_o-1)$-decoupled-pseudoconcave,  then, owing to $(b_{ij})\ge0$, we have
\begin{eqnarray}\Label{5.8}
H^r_{q_o}(u^\tau)\simge\sum_{j=1}^{k+1}\frac{\di^2 h_j}{\di z_j\di\bar z_j} \big( |u^\tau|^2-\sumK |u^\tau_{jK}|^2\big) \T{  on  } U\cap \bar\Om,
\end{eqnarray} for any $u^\tau$ of degree  $k\le q_o-1$.
\end{remark}

\begin{theorem}\Label{t5.7} Let $\Om$ be $(q_o+1)$-decoupled-pseudoconvex at $z_0=0$. 
Write coordinates in $\C^n$ as $z=x+iy$ and suppose that there are  invertible functions $F_j$, $j=q_o+1,...,n-1$, with $\frac{F_j(|t|)}{|t|^2}$ increasing with respect to  $t$ near $0$ such that 
\begin{eqnarray}\Label{5.9}
\T{~~~either~~~}\frac{\di^2 h_j}{\di z_j\di\bar z_j}(z_j)\simge \frac{F_j(|x_j|)}{x_j^2}, \T{~~~or~~~ }\frac{\di^2 h_j}{\di z_j\di\bar z_j}(z_j)\simge \frac{F_j(|y_j|)}{y_j^2}
\end{eqnarray}
Also, assume by reordering, that the $Fj$'s are increasing, that is, $...F_j\lesssim F_{j+1}...$.
Then $(f\T-1)^k$ estimates hold in degree $k\ge q_2$  with $f(\delta^{-1})=(F^{*}_k(\delta))^{-1}$ where $F_j^{*}$ are the inverse functions to the $F_j$'s .

Similarly, let $\Om$ be $(q_o-1)$-decoupled-pseudoconcave and assume that there are $F_j$ with $\frac{F_j(|t|)}{|t|^2}$ increasing in  $t$ and satisfying \eqref{5.9} for $j=1,...,q_o$; by reordering, further assume that $...F_j\simgeq F_{j+1}...$. Then $(f\T-1)^k$ estimates hold in degree $k\le q_o-1$ for $f(\delta^{-1})=(F^{*}_k(\delta))^{-1}$.
\end{theorem}

\begin{example}\Label{e5.1}\rm If $h_j(z_j)=|z_j|^{2m_j}$ or $|x_j|^{2m_j}$ with $m_j\ge m_{j+1}$, then we get $\epsilon$-subelliptic estimates with $\epsilon=\frac{1}{2m_k}$. If $h_j(z_j)=\exp(-\frac{1}{|z_j|^{m_j}})$ or $\exp(-\frac{1}{|x_j|^{m_j}})$ with $m_j\ge m_{j+1}$, then we get $(f\T-1)^k$ estimates with  $f(t)=(\log t)^\frac{1}{m_k}$. 
\end{example}

\begin{example}  \rm Let  $$r=2\T{Re}z_{n}-a(z_1,...,z_{q_o})+\sum_{j=q_o+1}^{n-1}h_j(z_j).$$
where $\di\dib a\ge 0$ and $h_j(z_j)$ is defined in Example  \ref{e5.1}. Then we get $(f\T-1)^k$ estimates at $z_0$ for  the domain $\Om^+=\{+r<0\}$ (resp. $\Om^-=\{-r< 0\}$) for any form of degree $k\ge (q_o+1)$ (resp. $k\le (n-q_o-2) $). By Theorem~\ref{main2} , $(f\T-1)^k$ estimates for the system $(\dib_b, \dib_b^*)$ on $M=\{r=0\}$ hold in any degree $k$ such that $q_o+1\le k\le n-q_o-2$. 
\end{example}

{\it Proof of Theorem \ref{t5.7}. } We may assume that $\frac{\di^2 h_j}{\di z_j\di\bar z_j}(z_j)\simge \frac{F_j(|x_j|)}{x_j^2}$ for any $j$ because,  if $\frac{\di^2 h_j}{\di z_j\di\bar z_j}(z_j)\simge \frac{F_j(|y_j|)}{y_j^2}$, we change coordinates by $z_j:=iz_j$. Let $C$ be a postive constant such that $C\frac{\di^2 h_j}{\di z_j\di\bar z_j}(z_j)\ge \frac{F_j(|x_j|)}{x_j^2}$.\\

{\it\underline{The case $(q_o+1)$-decoupled-pseudoconvex.}}
For any form of degre $k\geq q_o+1$, we define the family of weights $\Phi=\Phi^\delta_k$ by 
\begin{eqnarray}
\Phi =C{\frac{r}{\delta}}-2\sum_{j=k}^{n-1} \exp({-\frac{x_j^2}{4F_k^{*}(\delta)^2}}).
\end{eqnarray} 
The weights $\Phi$'s are absolutely bounded on $S_\delta$. Computation of the Levi form $\di\dib \Phi$ shows that
$$\di\dib \Phi=C\frac{\di\dib r}{\delta}+\frac{1}{F_k^*(\delta)^2}\sum_{j=k}^{n-1}\Big(1-\frac{x_j^2}{2F_k^*(\delta)^2}\Big)\exp(-\frac{x_j^2}{4F_k^{*}(\delta)^2})\om_j\we \bar\om_j.$$

Then
\begin{eqnarray}\Label{5.10}
\begin{split}
H^\Phi_{q_o}(u^\tau)=&\frac{C}{\delta}H^r_{q_o} (u^\tau)\\
&+\frac{1}{F_k^*(\delta)^2}\sum_{j=k}^{n-1}\Big(1-\frac{x_j^2}{2F_k^*(\delta)^2}\Big)\exp(-\frac{x_j^2}{4F_k^{*}(\delta)^2})\sumK|u^\tau_{jK}|^2.
\end{split}
\end{eqnarray}
Using \eqref{5.7} and the hypothesis \eqref{5.9}, we get the estimate
$$
\frac C8H^r_{q_o}(u^\tau)\ge\sum_{j=k}^{n-1}\left(\frac1\delta\frac{F_j(|x_j|)}{x_j^2}\right)\sumK|u_{jK}|^2.
$$
We use the notations

 $$A_j=\frac{1}{\delta}\frac{F(|x_j|)}{x_j^2}\T{ ~~~and~~~ } B_j=\frac{1}{F_k^*(\delta)^2}\Big(1-\frac{x_j^2}{2F_k^*(\delta)^2}\Big)\exp(-\frac{x_j^2}{4F_k^{*}(\delta)^2}).$$
where we notice  that $A_j\ge 0$ for any $j$.  
Thus \eqref{5.10} reads as
\begin{equation}
\Label{newb}
H^\Phi_{q_o}(u^\tau)\ge \sum_{j=k}^{n-1}(A_j+B_j)\sumK|u^\tau_{jK}|^2.
\end{equation}

 For each $j$ wit $k\le j\le n-1$,  we consider two cases 
\\
{\bf Case 1.} If $|x_j|\le F_k^*(\delta)$, we have
$$B_j\ge \frac{1}{2F_k^*(\delta)^2}e^{-1/4}\ge cF_k^*(\delta)^{-2};$$
and hence $$A_j+ B_j\simge F_k^*(\delta)^{-2}.$$
{\bf Case 2.} Otherwise, assume  $|x_j|\ge F_k^*(\delta)$. Using our assumption that $\frac{F(|x_j|)}{x_j^2}$ is increasing in $|x_j|$, we get
$$A_j=\frac{1}{\delta}\frac{F_j(|x_j|)}{x_j^2}\ge \frac{1}{\delta}\frac{F_k(F_k^*(\delta))}{F_k^*(\delta)^2}=\frac{1}{\delta}\frac{\delta}{F_k^*(\delta)^2}=F_k^*(\delta)^{-2}.$$
In this case, $B_j$ can get negative values; however, by using the fact that $\underset{t\ge \frac{1}{2}}\min\Big\{(1-t)e^{-t/2}\}=-2e^{-3/2}$ for $t=\frac{x_j^2}{2F_k^*(\delta)^2}\ge\frac{1}{2}$ we have
$$B_j\ge -2{e^{-3/2}}F^*_k(\delta)^{-2}.$$ 
This implies
$$A_j+B_j\simge F^*_k(\delta)^{-2}.$$

Therefore, continuing our estimate  \eqref{newb}, we obtain
\begin{eqnarray}
\begin{split}
H^\Phi_{q_o}(u^\tau)
\simge&\sum_{j=k}^{n-1}F_k^*(\delta)^{-2}\sumK|u^\tau_{jK}|^2\\
\simge&f(\delta)^2|u^\tau|^2,
\end{split}
\end{eqnarray}
where the last inequality follows from $\overset{n-1}{\underset{j=k}{\sum}}\sumK|u^\tau_{jK}|^2\ge\sumJ|u^\tau_J|^2= |u^\tau|^2$. Moreover, we see that $\Phi_j=\frac{r_j}{\delta}$ and therefore $\Phi_j=0$ for any $j\le q_o$. Hence $\sum_{j=1}^{q_1+1}|\Phi_ju^\tau|^2=0$.
 Thus the weights $\Phi=\Phi_k^\delta$ satisfy Property $(f\T-1\T-P)^k$ . Applying Theorem \ref{main1}, we get
$$\no{f(\Lambda)u}^2\lesssim Q(u,u)$$
for any $u\in C_c^\infty(U\cap\bar\Om)^k\cap \T{Dom}(\dib^*)$ with $f(\delta^{-1})=F^*_k(\delta)^{-1}$.\\

{\it\underline{The case $q$-decoupled-pseudoconcave.}}
For each $k\le q_o-1$, we define the family of weights $\Phi=\Phi^\delta_k$ by 
\begin{eqnarray}
\Phi =C{\frac{r}{\delta}}+2\sum_{j=1}^{k+1} \exp({-\frac{x_j^2}{4F_{k+1}^{*}(\delta)^2}}).
\end{eqnarray} 

In the same way as before, we see that the $\Phi_k^\delta$'s are absolutely bounded on $S_\delta$ and their Levi forms satisfy
\begin{eqnarray}\Label{5.15}
\begin{split}
H^\Phi_{q_o}(u^\tau)
\simge&F_{k+1}^*(\delta)^{-2}\sum_{j=1}^{k+1}\Big(|u^\tau|^2-\sumK|u^\tau_{jK}|^2\Big)\\
\simge&F_{k+1}^*(\delta)^{-2}\Big((k+1)|u^\tau|^2-\sum_{j=1}^{k+1}\sumK|u^\tau_{jK}|^2\Big)\\
\simge&F_{k+1}^*(\delta)^{-2}|u^\tau|^2.
\end{split}
\end{eqnarray}

Now, we need to show that
\begin{eqnarray}\Label{5.16}
\begin{split}
H^\Phi_{q_o}(u^\tau)\simge \Big(\sum_{j=1}^{q_o}|\Phi_j|^2\Big) |u^\tau|^2.
\end{split}
\end{eqnarray}

We note that $\Phi_j=0$ for $k+2\le j\le n-1 $ and that
\begin{eqnarray}\Label{5.17}
|\Phi_j|^2=4F^*_j(\delta)^{-2} \Big(\frac{x_j^2}{2F^*_j(\delta)^2}\exp(-\frac{x_j^2}{2F^*_j(\delta)^2})\Big)\lesssim F^*_j(\delta)^{-2} 
\end{eqnarray}
for $1\le j\le k+1$. 
Thus, from \eqref{5.15} and \eqref{5.17} we get \eqref{5.16}.\\

$\hfill\Box$

\section{Subelliptic estimates for regular coordinate domains}
We state precise subelliptic estimates for the $\bar\partial$-Neumann problem over the class of {\it regular coordinate domains} of $\C^{n+1}$.\\

\index{Domain ! regular coordinate}

We consider a domain $\Omega\subset\C^{n+1}$ defined by
\begin{equation}
\Label{1}
2\T{Re} z_{n+1}+\underset{j=1}{\overset{n}\sum}|f_j(z)|^2<0,
\T{ with $f_j$ holomorphic, $f_j=f_j(z_1,...,z_j)$ and $\partial_{z_j}^{m_j}f_j\neq0$}.
\end{equation}
 This is called a { \it regular coordinate domain}; the inequality which defines $\Omega$ is denoted by $r<0$.
We discuss { \it subelliptic} estimates for the $\bar\partial$-Neumann problem on $\Omega$ 
\begin{equation}
\Label{subelliptic}
|||u|||^2_\epsilon\simleq ||\bar\partial u||^2_0+||\bar\partial^*u||^2_0+||u||^2_0\quad
\T{for any $u\in C^\infty_c(U\cap \bar\Om)^1\cap\T{Dom}(\bar\partial^*)$}.
\end{equation}
Our problem is to find the { \it optimal} $\epsilon$.
There are several relevant numbers related to $\Omega$ 
\begin{itemize}
\item
$ m=m_1\cdot...\cdot m_{n}$ the {\it multiplicity}.
\item
$D$ the {\it D'Angelo type} defined as the maximal order of contact of a complex curve with $b\Omega$. Note here that since $\sum_j|f_j|^2\simgeq |z|^{2m}$ then necessarily $D\leq 2m$.
\item
$\epsilon$ the (optimal) index of subelliptic estimates. It satisfies
$$
\epsilon\leq \frac1D\T{ Catlin 1983 \cite{C83}},\qquad \epsilon\overset ?\geq \frac1{2m}
\T{ D'Angelo conjecture 1993 \cite{D93}}.
$$
\end{itemize}
We define a new number $\gamma$. For this, we write
$$
f_j=g_j(z_1,...,z_j)+z_j^{m_j}+O(z_j^{m_j+1})\quad\T{ for $g_j=O(z_1^{{\lambda^1_j}},...,z_{j}^{\lambda^{j-1}_j})$}.
$$
Let $j_o$, $j_o\geq2$,  $j\geq2$, be the first index
with the property that $g_j$ is independent of $z_j$ for any $j\leq j_o-1$ and write  $l_j^i$ for the minimum between $\lambda^i_j$ and $m_i$ (resp. $m_i-\eta$ for any $\eta>0$) when $i\leq j_o-1$ and $j\leq j_o-1$ (resp. $j\geq j_o$), and otherwise put $l_j^i=1$. Define
\begin{equation}
\Label{gamma}
\gamma_j=\underset{i\leq j-1}\min\frac{l_j^i}{m_j}\gamma_i,\qquad \gamma=\underset j\min \gamma_j.
\end{equation}
Note that $\frac1{2m}\leq\frac\gamma2\leq\frac1D$.
Here is our main result of \cite{KZ08c} (which  is also presented in \cite{CC08} with $\frac\gamma2$ replaced by $\frac1{2m}$).
\begin{theorem}
\Label{t1}
Let $\Omega $ be a regular coordinate domain and let $\gamma$ be the number defined by \eqref{gamma}; then we have $\epsilon$-subelliptic estimates for $\epsilon=\frac\gamma2$.
\end{theorem}
\begin{example}\rm
\Label{e1}
For the domain of $\C^{n+1}$ defined by
$$
2\T{Re} z_{n+1}+|z_1^{m_1}|^2+\underset {j=2}{\overset{n}\sum}|z_j^{m_j}-z_{j-1}^{l_j}|^2<0,\quad l_j\leq m_{j-1}\leq m_j,
$$
the number $\gamma$ is given by  $\gamma=\frac{l_2\cdot...\cdot l_{n}}{m_1\cdot...\cdot m_{n}}$. We claim that $\epsilon=\frac\gamma2=\frac1D$. The first equality follows from Theorem~\ref{t1}. As for the second, we can easily find the  { \it critical curve} $\Gamma$ with the maximal contact with $\partial\Omega$; this is  parameterized over $\tau\in \Delta$ by
$$
\tau\mapsto (\tau^{\frac1{m_1\gamma}},\tau^{\frac{l_2}{m_2m_1\gamma}},...,\tau,0).
$$
In fact, we have
$
r|_\Gamma=|\tau^{\frac1\gamma}|^2+|\not{\tau^{\frac{l_2}{m_1\gamma}}}-\not{\tau^{\frac{l_2}{m_1\gamma}}}|+...=|\tau^{\frac1\gamma}|^2$
and therefore $D\geq \frac2\gamma$. On the other hand  $\epsilon=\frac\gamma2<\frac1D$  by Catlin 1983 \cite{C83}.
\end{example}
\begin{example}\rm
For the domain in $\C^3$ defined by
$$
2\T{Re} z_3+|z_1^4|^2+|z_2^6+\underset{j=0}{\overset5\sum}c_jz_2^jz_1^{\alpha_j}+O(z_2^7)|^2<0,
$$
with $\alpha_j\geq 3$ for any $j$, we have $\gamma_1=\frac14,\,\,\gamma_2=\frac3{6\cdot 4}$ and $\gamma=\gamma_2$.
\end{example}
\begin{example}\rm
\Label{e2}
Let us consider in $\C^4$ the domain defined by
$$
2\T{Re}  z_4+|z_1^6|^2+|z_2^4+z_1^3|^2+|z_3^4+z_3z_1^a+z_3z_2^b|^2<0.
$$
Here $\gamma_1=\frac16$, $\gamma_2=\frac3{6\cdot 4}$ and  $\gamma_3=\min(\frac{3\cdot b}{6\cdot 4\cdot 4},\frac a{6\cdot 4},\frac {3}{6\cdot 4})$ and $\gamma=\gamma_3$; in particular, if $a\geq3$, and $b\geq4$, we have $\gamma=\frac {3}{6\cdot 4}$.
\end{example}

\noindent
{\it Proof of Theorem~\ref{t1}.}

In order to establish \eqref{subelliptic} it suffices to find a family of bounded weights $\{\Phi^\delta\}$  whose Levi form satisfies over the strip $S_\delta:=\{z\in\Omega:\,-r(z)<\delta\}$, the estimate $\partial\bar\partial\Phi^\delta(z)(u,\bar u)\simgeq \delta^{-\gamma}|u|^2$ for $u\in\C^n$. We choose $\alpha\geq 1$, put $\alpha_j=\alpha (m_{j+1}...m_n)$ and  choose a smooth cut off function $\chi$ with $\chi\equiv1$ in $[0,1]$ and $\chi\equiv0$ in $[2,+\infty)$. We define
\begin{equation}
\Label{10}
\begin{split}
\Phi^\delta=-\log(\frac{-r+\delta}\delta)+&\underset{j=j_o}{\overset n\sum}\underset{h=1}{\overset{m_j-1}\sum}\frac1{|\log\,*|}\log\left(|\partial_{z_j}^hf_j|^2+\frac{\delta^{(m_j-h)\gamma_j}}{|\log\delta|^{(m_j-h)\alpha_j}}\right)
\\
&+c\underset{j=1}{\overset n\sum}\chi(\frac{|z_j|^2}{\delta^{\gamma_j}})\log\left(\frac{|z_j|^2+\delta^{\gamma_j}}{\delta^{\gamma_j}}\right),
\end{split}
\end{equation}
where $*=\frac{\delta^{\gamma_j(m_j-h)}}{|\log\delta|^{(m_j-h)\alpha_j}}$; notice  that $\log *\sim \log \delta$. The weights $\Phi^\delta$ that we have defined are bounded in the strip $S_\delta$.
When calculating the Levi form, we observe that $\partial\bar\partial r(u,\bar u)=\underset j\sum|\partial f_j\cdot u|^2$ and $\partial\bar\partial|\partial^h_{z_j}f_j|^2(u,\bar u)=|\partial\partial^h_{z_j}f_j\cdot u|^2$ and finally $\partial\bar\partial |z_j|^2(u,\bar u)=|u_j|^2.$ Thus, we have got a decomposition
\begin{equation}
\Label{10,5}
\partial\bar\partial \Phi^\delta(u,\bar u)=\underset {j=1}{\overset n\sum} A_j+\underset {j=j_o}{\overset n\sum}\underset {h=1}{\overset {m_j-1}\sum } B^h_j+\underset {j=1}{\overset n\sum}  C_j,
\end{equation}
with the estimates
\begin{equation*}
\begin{cases}
A_j\simgeq\delta^{-1}|\partial f_j\cdot u|^2 &\T{ for any $z\in S_\delta$},
\\
B_j^h\simgeq {\delta^{-(m_j-h)\gamma_j}}\frac{|\log\delta|^{(m_j-h)\gamma_j}}{|\log*|}|\partial\partial_{z_j}^hf_j\cdot u|^2&\T{ if $|\partial_{z_j}^hf_j|^2<\frac{\delta^{(m_j-h)\gamma_j}}{|\log\delta|^{(m_j-h)\alpha_j}}$},
\\
C_j\simgeq c\delta^{-\gamma_j}|u_j|^2&\T{ if $|z_j|^2<\delta^{\gamma_j}$}.
\end{cases}
\end{equation*}
We note that the $A_j$'s and $B_j^h$'s are positive for any $z$ but, instead, the $ C_j$'s can take negative values when $|z_j|>\delta^{\gamma_j}$; however, $|C_j|\simleq c\delta^{-\gamma_j}|u_j|^2$ and thus the $C_j$'s are controlled by the $A_j$'s and $B_j^h$'s. We define
$$
 D_j=A_j+C_j\quad\T{for $j\leq j_o-1$},\qquad D_j=A_j+\underset{h\leq m_j-1}\sum B_j^h+ C_j\quad\T{for $j\geq j_o$}.
$$
We wish to start by proving that, when $j\leq j_o-1$, then
\begin{equation}
\Label{2.12}
\underset{i\leq j}\sum D_i\simgeq
\underset{i\leq j}\sum
\delta^{-s_i\gamma_i}|z_i|^{2(s_i-1)}|u_i|^2\T{ for any $s_i\leq m_i$}.
\end{equation}
We use induction and show how to pass from step $j-1$ to step $j$ (the step $j=1$ being elementary). We fix our choice $s_i=l^i_j$
and remark that
\begin{equation}
\Label{2.12bis}
\begin{split}
A_j+\underset{i\leq j-1}\sum D_i&\simgeq \delta^{-1}\left| \partial f_j\cdot u\right|^2
+\underset{i\leq j-1}\sum\delta^{-l_j^i\gamma_i}|z_i|^{2(l_j^i-1)}|u_i|^2
\\
&\simgeq \underset{i\leq j-1}\sum\delta^{-l_j^i\gamma_i}\left[|z_j|^{2m_j-1}|u_j|^2-{|z_i|^{2(l_j^i-1)}|u_i|^2}\right]+\underset{i\leq j-1}\sum\delta^{-l_j^i\gamma_i}{|z_i|^{2(l_j^i-1)}|u_i|^2}.
\end{split}
\end{equation}
This proves \eqref{2.12} for $s=m_j$. On the other hand, we have
\begin{equation}
\Label{2.12ter}
C_j+\delta^{-m_j\gamma_j}|z_j|^{2(m_j-1)}\simgeq (\delta^{-\gamma_j}+\delta^{-m_j\gamma_j}|z_j|^{2(m_j-1)})|u_j|^2.
\end{equation}
This is clear for $|z_j|^2\leq \delta^{\gamma_j}$; otherwise,   $C_j$ gets negative
but it is controlled by the second term in the left of \eqref{2.12ter} for small $c$. By combining \eqref{2.12bis} with \eqref{2.12ter} we get \eqref{2.12} for $s=m_j$ and $s=1$ and thus also for any  $1\leq s\leq m_j$.
 This concludes the proof of our claim \eqref{2.12}.
 
We pass to treat the terms
$
D_j=A_j+\underset h\sum B_j^h+C_j\quad\T{for $j\geq j_o$}.
$
We begin by an auxiliary statement: if for some $i$ with $j_o\leq i\leq j-1$ and for any $1\leq s_{i'}\leq m_{i'}$, we have
\begin{equation}
\Label{20bis}
\underset{i'\leq {i}}\sum D_{i'}\simgeq \delta^{-\gamma_{{i}}}|\log\delta|^{\alpha_{{i}}-1}|u_{{i}}|^2+\underset{j_o\leq i'\leq {i}-1}\sum\delta^{-\gamma_{i'}}|u_{i'}|^2+\underset{i'\leq j_o-1}\sum \delta^{-s_{i'}\gamma_{i'}}|z_{i'}|^{2(s_{i'}-1)}|u_{i'}|^2,
\end{equation}
then we also have
\begin{equation}
\Label{*}
\underset{i\leq j}\sum D_i\simgeq \delta^{-\gamma_j}|\log\delta|^{\alpha_j-1}|u_j|^2.
\end{equation}
We  prove the implication from step $j-1$ to $j$.
By choosing $s_{i'}=l_j^{i'}$ in \eqref{20bis}, and observing that, if $i\leq j-2$, then $\delta^{-\gamma_i}\geq \delta^{-\gamma_{j-1}}|\log\delta|^a$ for any $a$, we get
\begin{equation}
\Label{**}
\begin{split}
A_j+\underset{i\leq j-1}\sum D_i&\simgeq\underset{i\leq j_o-1}\sum\delta^{-\gamma_il^i_j}\left(|\partial_{z_j}f_j|^2|u_j|^2-|z_i|^{2(l_j^i-1)}|u_i|^2\right)
\\&+\underset{i\geq j_o}\sum\delta^{-\gamma_{j-1}}|\log\delta|^{\alpha_{j-1}-1}\left(|\partial_{z_j}f_j|^2|u_j|^2-|u_i|^2\right)
+\underset{i\leq j-1}\sum D_i
\\
&\simgeq \delta^{-\gamma_{j-1}}|\log\delta|^{\alpha_{j-1}-1}|\partial_{z_j}f_j|^2|u_j|^2.
\end{split}
\end{equation}
If now  $|\partial_{z_j}f_j|^2\geq \frac{\delta^{(m_j-1)\gamma_j}}{|\log\delta|^{(m_j-1)\alpha_j}}$, then \eqref{**} can be continued by
\begin{equation*}
\begin{split}
{}&
\geq \delta^{-\gamma_{j-1}+(m_j-1)\gamma_j}|\log\delta|^{(\alpha_{j-1}-1)-(m_j-1)\alpha_j}|u_j|^2
\\
&\geq \delta^{-\gamma_j}|\log\delta|^{\alpha_j-1}|u_j|^2.
\end{split}
\end{equation*}
If not, we pass to use $B_j^1$ instead of $A_j$. In this way we jump from $\partial_{z_j}^{h}f_j$ to $\partial_{z_j}^{h+1}f_j$ until we reach $B_j^{m_j-1}$; since $|z_j|^2=|\partial_{z_j}^{m_j-1}f_j|^2$ is smaller than $\frac{\delta^{\gamma_j}}{|\log\delta|^{\alpha_j}}$ (otherwise we would have used the former term $B_j^{m_j-2}$), then $B_j^{m_j-1}$ is bigger than the right side of \eqref{*}. This concludes the proof of the auxiliary statement.

We show that, for any value of $|z_j|^2$
\begin{equation}
\Label{11}
\underset{i\leq j}\sum D_i\simgeq \delta^{-\gamma_j}|u_j|^2,
\end{equation}
whereas, when $|z_j|^2\geq \delta^{\gamma_j}$
\begin{equation}
\Label{12}
\underset{i\leq j}\sum D_i\simgeq
\begin{cases}
\T{\rm either} &\delta^{-\gamma_j}|\log\delta|^{\alpha_j-1}|u_j|^2
\\
\T{\rm or}&\delta^{-\gamma_{j-1}}|z_j|^{2(m_j-1)}|u_j|^2.
\end{cases}
\end{equation}
For $j\leq j_o-1$, the claim has already been proved in \eqref{2.12}: the second alternative in \eqref{12} holds. If $|z_j|^2\leq \delta^{\gamma_j}$, then $C_j\simgeq c\delta^{-\gamma_j}|u_j|^2$. Assume therefore $|z_j|^2\geq \delta^{\gamma_j}$ and suppose \eqref{12} true up to step $j-1$; we prove that it also holds for $j$ (which  implies \eqref{11}).
First, if in the inductive statement it is the first of \eqref{12} which is fulfilled at some step $i$ with $j_o\leq i\leq j-1$, the first is also fulfilled at step $j$; this follows from the auxiliary statement. Otherwise, assume we have the second for any $i\leq j-1$; (we surely do have for any $i\leq j_o-1$).
Now, if for some $i\leq j_o-1$, we have $|z_i|^{l_j^i}\geq |z_j|^{m_j-1}$, or, for some $j_o\leq i\leq j-1$, we have $|z_i|\geq |z_j|^{m_j-1}$, then, owing to $|z_j|\geq \delta^{\gamma_j}$, we have
\begin{equation}
\Label{12bis}
\begin{cases}
\delta^{-1}|z_i|^{2(m_i-1)}\geq \delta^{-\gamma_il^i_j-\eta}|z_i|^{2(l_j^i-1)},&i\leq j_o-1,
\\
\delta^{-\gamma_{i-1}}|z_i|^{2(m_i-1)}\geq \delta^{-\gamma_i-\eta},&i\geq j_o.
\end{cases}
\end{equation}
To prove the second of \eqref{12bis}, it suffices to notice that $$\delta^{-\gamma_{i-1}}|z_i|^{2(m_i-1)}\geq\delta^{-\gamma_{i-1}+\gamma_j(m_j-1)(m_i-1)}\geq \delta^{-\gamma_i-\eta}.$$ As for the first, we notice that
\begin{equation*}
\begin{split}
\delta^{-1}|z_i|^{2m_i-1}&\geq \delta^{-1+(m_i-l_j^i)\frac{m_j-1}{l_j^i}\gamma_j}|z_i|^{2(l_j^i-1)}
\\
&\geq\delta^{-l^i_j\gamma_i-(\frac{m_i}{l^i_j}-1)\gamma_j}|z_i|^{2(l_j^i-1)}\geq\delta^{-l^i_j-\eta}|z_i|^{2(l_j^i-1)},
\end{split}
\end{equation*}
(because $l_j^i\leq m_i-\eta$).
This proves \eqref{12bis}.
By \eqref{12bis}, the second of \eqref{12} is converted into the first in the inductive statement for $i\leq j-1$ (and thus also for $j$ owing to the auxiliary statement). Thus the only critical case occurs when both the inequalities
\begin{equation}
\Label{13}
\begin{cases}
|z_i|^{l_j^i}\leq |z_j|^{m_j-1},&i\leq j_o-1,
\\
|z_i|\leq |z_j|^{m_j-1},&i\geq j_o,
\end{cases}
\end{equation}
are fulfilled. But we have in this situation
\begin{equation*}
\begin{split}
|\partial_{z_j}f_j|^2&\geq |z_j|^{2(m_j-1)}-\frac12\left(\underset{i=1}{\overset{j_o-1}\sum}|z_i|^{l_j^i-1}+\underset{i=j_o}{\overset{j-1}\sum}|z_i|^2\right)
\\
&\geq |z_{j}|^{2(m_j-1)},
\end{split}
\end{equation*}
which implies $A_j+\sum_i D_i\simgeq \delta^{-\gamma_{j-1}}|z_j|^{2(m_j-1)}|u_j|^2$. This yields the second of \eqref{12}.
Thus induction works and  brings us to step $j=n$. At this point we can disregard \eqref{12} (though it did a great job for the inductive argument): \eqref{11} for any $j\leq n$ yields the conclusion of the proof.\\

$\hfill\Box$


\label{Bibliography}

\bibliographystyle{alphanum}  


\chapter*{Vita}
\addcontentsline{toc}{chapter}{Vita}

\begin{center}
		TRAN VU KHANH\\
(Email : khanh@math.unipd.it; vukhanh83@yahoo.com)
\end{center}
\begin{enumerate}
  \item[Sept. 1983~: ]  Born, Ca Mau, Viet Nam.
  \item[2001-2005~: ]  Undergraduate student at Hochiminh University of Natural Sciences, Vietnam.
\item[~~~~2006~~~~ : ]  Teaching assistant at Department of Mathematics and Computer
Sciences, Hochiminh University of Natural Sciences, Vietnam.
\item[2007-2009 ~: ] Ph.D. student at University of Padova, Italy.
\end{enumerate}

\end{document}